\documentclass[preprint,sort&compress,12pt]{elsarticle}
\usepackage{bm}
\usepackage{amsmath}
\usepackage{mathrsfs}
\usepackage{stmaryrd}
\usepackage{amsfonts}
\usepackage{amssymb}
\usepackage{amsthm}
\usepackage{xcolor}

\newcommand{\mm}{\mathrm}
\newcommand{\ml}{\mathcal}

\newcommand{\be}{\begin{equation}}
\newcommand{\bea}{\begin{equation}\begin{aligned}}
\newcommand{\beas}{\begin{equation*}\begin{aligned}}
\newcommand{\eeas}{\end{aligned}\end{equation*}}
\newcommand{\eea}{\end{aligned}\end{equation}}
\newcommand{\ee}{\end{equation}}

\makeatletter
\def\ps@pprintTitle{%
	\let\@evenfoot\@oddfoot
}
\makeatother
\begin{document}
\begin{frontmatter}
\title{Global Solvability for the Compressible Hookean Viscoelastic Fluids\\ with a Free Boundary in Some Classes of Large Data}

\author[aa,sJ]{Fei Jiang}
\ead{jiangfei0591@163.com}
\author[aa,sJ]{Youyi Zhao\corref{cor1}}
\ead{zhaoyouyi957@163.com}
\cortext[cor1]{Corresponding author. }
\address[aa]{School of Mathematics and Statistics, Fuzhou University, Fuzhou, 350108, China.}
\address[sJ]{Key Laboratory of Operations Research and Control of Universities in Fujian, Fuzhou 350108, China.}
\begin{abstract}
Recently Jiang--Jiang established a global (in time) existence result for unique strong solutions of the two-dimensional (2D) free-boundary problem of an incompressible Hookean viscoelastic fluid, the rest state of which is defined in a slab, in some classes of large data \cite{JFJSGS}. In particular, Jiang--Jiang's mathematical result shows that, if the initial free boundary is flat, the way the elastic deformation under the large  elasticity coefficient $\kappa$ acts on the free boundary prevents the natural tendency of the fluid to form singularities, even when the initial velocity is properly large. However it is not clear whether their result can be extended to the corresponding \emph{3D case}. In this paper, we further find a similar result in the \emph{3D} stratified (immiscible) compressible Hookean viscoelastic fluids in an infinite slab with two restrictive conditions: that the elasticity coefficients of two fluids are equal, and that the initial density functions satisfy the asymptotic stability condition in Lagrangian coordinates. These two restrictive conditions in the compressible case contribute us to avoid the essential obstacles that would be faced in the extension of Jiang--Jiang's result from two dimensions to our 3D case. In addition, we can further obtain a new result regarding the vanishing phenomena of the nonlinear interactions of solutions with the fixed initial velocity and the initial zero perturbation deformation. Such a new result roughly presents that  the solutions of the problem considered by us can be approximated by the ones of a linear problem for sufficiently large $\kappa$.
\end{abstract}
\begin{keyword}
Compressible viscoelastic fluids; Stratified fluids;
Exponential time-decay; Large initial velocity; Vanishing phenomena of nonlinear interactions.
\end{keyword}
\end{frontmatter}


\newtheorem{thm}{Theorem}[section]
\newtheorem{lem}{Lemma}[section]
\newtheorem{pro}{Proposition}[section]
\newtheorem{cor}{Corollary}[section]
\newproof{pf}{\emph{Proof}}
\newdefinition{rem}{Remark}[section]
\newtheorem{definition}{Definition}[section]

\section{Introduction}\label{sec:01}
\numberwithin{equation}{section}

Viscoelastic materials include a wide range of fluids with elastic properties, as
well as solids with fluid properties. The models of viscoelastic fluids formulated by Oldroyd, in particular the classical Oldroyd-B model,
have been studied by many authors. In this paper, we consider the following incompressible Oldroyd model which includes a viscous stress component
and a stress component for a neo-Hookean solid \cite{hu2010local,QJZZPWZZF}:
\begin{equation}\label{0101}
\begin{cases}
\rho v_{t}+\mm{div}\left(\rho v\right) =0,\\[0.5mm]
\rho v_{t}+\rho v\cdot\nabla v
+\mm{div}\mathcal{S}(\rho,v,U)=0, \\[0.5mm]
 {U}_{t}+ v\cdot\nabla {U}= \nabla v{U},
\end{cases}
\end{equation}
where the tensor $\mathcal{S}(\rho,v,U)$ is given as follows:
\begin{align}
\label{2020202109232002205}
\mathcal{S}(\rho,v,U):=P(\rho)\mathbb{I}-\mathcal{V}(v)-\kappa\left(UU^\top/\mm{det}U-\mathbb{I}\right).
\end{align}
We call \eqref{0101}$_1$ the continuity equation,
\eqref{0101}$_2$ the momentum equations and \eqref{0101}$_3$ the deformation equations.
Next we shall introduce the mathematical notations in the above motion equations.

The unknowns ${\rho}:={\rho}(x,t)$, ${v}:={v}(x,t)$ and ${U}:={U}(x,t)$
represent the density, velocity and deformation tensor (a 3$\times$3 matrix valued function) of the viscoelastic fluid, resp.
The pressure function $ P(\tau)\in C^{4}(\mathbb{R}^+) $ in \eqref{0101} is always assumed to be positive and strictly
increasing with respect to $\tau$,
and the viscosity (stress) tensor $\mathcal{V}(v)$ is given by
$$\mathcal{V}(v):=\mu\mathbb{D}v+\lambda\mm{div}v\mathbb{I}\;\mbox{ with }\;
\mathbb{D}v:=\nabla v+ \nabla v^\top,$$
where $\lambda:= \varsigma-2\mu/3 $, and the parameters $\mu>0$ and $\varsigma\geqslant 0$ represent the shear and bulk viscosity coefficients, resp.;
$\kappa$ denotes the elasticity coefficient (i.e., the ratio between the kinetic and elastic energies)  \cite{LFHSM}, and the term $\kappa\mm{div}\left(UU^\top/\mm{det}U-\mathbb{I}\right)$  represents the elasticity, where $\mm{det}U$ means the determinant of $U$.
In addition the letter $\mathbb{I}$ and the superscript $\top$ in \eqref{2020202109232002205} are the $3\times3$ identity matrix and  the transpose, resp.

For strong solutions of  both the Cauchy and the initial-boundary value problems of \eqref{0101}, the authors in \cite{HXWDHTTIVPTCVF2,QJZZPWZZF,HXWDHTTIVPTCVF1,hu2012strong,HXPWGC} have established the global (in time) existence of solutions in various functional spaces
whenever the initial data is a small perturbation around the rest state $(\bar{\rho}, 0, \mathbb{I})$, where the constant $\bar{\rho}>0$ denotes the equilibrium density.
The global existence of weak solutions to \eqref{0101} with small perturbations near the rest state is established by Hu  \cite{HXPGETDCMFJDE}.
We also refer to  \cite{MR4435915,MR4350195} in relevant progresses.
However, it is still a longstanding open problem whether a global solution of the equations of compressible viscoelastic fluids  exists for any general large initial data, even in the two-dimensional case.  At present, other mathematical topics for viscoelastic fluids
have also been widely investigated, such as the decay-in-time of solutions with small perturbations \cite{MR3927502,MR4152207,HXPWGC}, the incompressible limit \cite{LZZYGE},
the corresponding incompressible case \cite{HXPLFHG,LFHLCZPO2,LFZPGCon2,chen2006global,LZLCZYGARMA}
and so on.
We also refer the readers to \cite{MR4335131,MR4064197,MNIMG,MR4059368,LPLMNCAMN,MR2990047} and the references cited therein for the local/global existence of solutions to the other closely related models in viscoelastic fluids.

Viscoelasticity is a material property that exhibits both viscous and elastic characteristics under deformation. In particular, an elastic fluid strains when it is stretched and quickly returns to its original state once the stress is removed. This means that the elasticity will has a stabilizing effect on the motion of viscoelastic fluids.   In fact Jiang et al. verified that the elasticity can inhibit the Rayleigh--Taylor instability and thermal instability for a properly large elasticity coefficient $\kappa$ in Hookean viscoelastic fluids \cite{JFJWGCOSdd,FJWGCZXOE,jiang2020,JFJSCVPDE1}. Recently, Ishigaki also found that the elasticity in viscoelastic fluids can accelerate the decay-in-tim of the $L^p$ estimates with $p>2$ \cite{MR4152207}.  This stabilizing effect can be also observed in the incompressible elastic fluids (i.e., in the absence of viscosity). For example, in \cite{STCBT1,STCBT2}, the authors established an interesting global existence results of global (in time) classical solutions  with  small initial data for the incompressible  elastic fluids. Note that such a result is not known for Euler equations, where the elastic effect is not present. Later Lei and Wang further independently obtained the global (in time) classical solutions with the stability of lower-order derivatives for the 2D elastic fluid system \cite{LZSTCZY2,WXCG}. In addition, the stabilizing effect of elasticity on the local-in-time motion of elastic fluids can be found in the free-boundary case, interested readers can refer to \cite{CCRMHJLWDH} for the vortex sheet problem and \cite{LHWWZZF} for the Rayleigh--Taylor problem. We mention that recently the well-posedness of the free-boundary problems of incompressible/compressible elastic/viscoelastic fluids have been widely investigated, see \cite{GL2022,ZJY2022,TrakhininJDE,GuMei2023,XLZPZZFGAR} and the references cited therein.

Motivated by the inhibition phenomenon of the Rayleigh--Taylor instability  by elasticity in viscoelastic fluids for a properly large elasticity coefficient
$\kappa$ in  \cite{JFJWGCOSdd},  Jiang--Jiang obtained a global existence result of strong solutions to the initial value problem of an incompressible/compressible Hookean viscoelastic fluid defined in a periodic domain, when the initial velocity is \emph{relatively} smaller than the elasticity coefficient \cite{JFJSsdafsafJMFMOSERT}. This means that the strong elasticity can prevent the development of singularities even when the initial velocity is large, thus playing a similar role to viscosity in preventing the formation of singularities in viscous flows. Recently, Jiang--Jiang
further established a similar result in the system of incompressible viscoelastic fluids with a free boundary by developing new ideas.
More precisely, for any given initial velocity $v|_{t=0}$ perturbing the rest state of an incompressible viscoelastic fluid with a free boundary, the elasticity
under the relatively larger  $\kappa$  can prevent the formation of singularities on the free boundary for the 2D case, in which the fluid domain of the rest state is a slab with a fixed lower boundary and an upper free boundary. Jiang--Jiang's mathematical result shows that, if the initial free boundary is flat, the way the elastic deformation under the relatively larger $\kappa$ acts on the free boundary prevents the natural tendency of the fluid to form singularities, even when the initial velocity is properly large.
We mention  that the self-overlap singularities of the free boundary may form if $\kappa$ is relatively smaller than
the initial perturbation velocity, please refer to Di Iorio--Marcati--Spirito's result for  the formation of splash  singularities on the  free boundary of a 2D incompressible viscoelastic fluid \cite{di2020splash}, which extends the result obtained
by Castro, C\'ordoba, Fefferman, Gancedo and G\'omez-Serrano for the Navier--Stokes equations in \cite{castro2019splash}.

Unfortunately it seems to be difficult to further extend Jiang--Jiang's result for the 2D free-boundary problem of incompressible viscoelastic fluids in \cite{JFJSGS} to the corresponding 3D case due to some essential obstacles arising from the free boundary and the pressure function. In this paper, however, we find a similar result in the stratified (immiscible) \emph{compressible} viscoelastic fluids (with an internal free interface) in a slab, where both  elasticity coefficients of the two fluids shall be equal (i.e., the condition of uniform elasticity coefficients) and the initial density functions satisfy the asymptotic stability condition of density (see \eqref{20202109250111564}) in Lagrangian coordinates (see Theorem \ref{thm:2109n} for details).  In addition, we further derive the vanishing phenomena of the nonlinear interactions (i.e., the solution of the nonlinear problem can be approximately by the one of a linear problem) in the Lagrangian coordinates as $\kappa\to \infty$ (see Theorem \ref{thm:210902} for details). We mention that such a conclusion can not be expected in the incompressible case with a free boundary, as discussed in \cite{JFJSGS}.

\section{Main results}
Before further stating our main results in details,  we shall first formulate our problem.
\subsection{A Model of stratified Hookean viscoelastic fluids}

We consider two distinct, immiscible, Hookean viscoelastic fluids
evolving in a moving domain $\Omega(t)=\Omega_+(t)\cup\Omega_{-}(t)$ for time $t\geqslant 0$. The upper fluid fills the upper domain
\begin{align}
\label{20224102301821}
 \Omega_+(t):=\{(x_{\mm{h}},x_3)^\top~|~x_{\mm{h}}:=(x_1,x_2)^\top\in \mathbb{R}^2 ,\ d(x_{\mm{h}},t)< x_3<h_+\},
\end{align}
and the lower fluid fills the lower domain
\begin{align}
\label{202241023018211}
\Omega_-(t):=\{(x_{\mm{h}},x_3)^\top~|~x_{\mm{h}}\in \mathbb{R}^2 ,\ h_-< x_3<d(x_{\mm{h}},t)\}.
\end{align}
We assume that $h_+$ and $h_-$ are two fixed and given constants satisfying $h_-<h_+$, but
the internal surface function $d:=d(x_{\mm{h}}, t)$ is free and unknown.
The internal surface
\begin{align}
\label{202410301926}
\Sigma(t):=\{x_3=d\}
\end{align} moves between the two Hookean viscoelastic   fluids,
and
$\Sigma_\pm:=\{x_3=h_\pm\}$ are the fixed upper and lower boundaries of $\Omega(t)$, resp.

We use the equations in \eqref{0101}  to describe the motion of the stratified compressible Hookean viscoelastic fluids, and
 add the subscript $_+$, resp. $_-$ to the notations of the known  physical parameters, pressure functions and other unknown functions
 in \eqref{0101} for the upper, resp. lower fluids. Thus  the motion equations of the stratified compressible Hookean viscoelastic fluids read  as follows:
\begin{equation}\label{0101f1}
\begin{cases}
\partial_t\rho_\pm+\mm{div}\left(\rho_{\pm} v_\pm\right)=0&\mbox{in } \Omega_\pm(t),\\[0.5mm]
\rho_\pm \partial_t v_\pm+\rho_\pm v_\pm\cdot\nabla v_\pm+\mm{div}\mathcal{S}_\pm(\rho_\pm, v_\pm, U_\pm)=0
&\mbox{in } \Omega_\pm(t),\\[0.5mm]
\partial_t U_\pm+v_\pm\cdot \nabla U_\pm=\nabla v_\pm U_\pm
&\mbox{in } \Omega_\pm(t).
\end{cases}
\end{equation}
Here the notation $\mathcal{S}_\pm$ are defined by \eqref{2020202109232002205} with $(v_\pm, U_\pm,P_\pm, \mu_\pm,\varsigma_\pm,\kappa_\pm)$ in place of $(v, U,P, \mu,\varsigma,\kappa)$, $P_\pm:=P_\pm(\tau)|_{\tau=\rho_\pm}$ denote the two different pressure functions of the upper and lower Hookean viscoelastic fluids.

 Motivated by both models of the stratified compressible viscous fluids in \cite{JJTIWYJTC} and the incompressible Hookean viscoelastic fluid with a free boundary \cite{XLZPZZFGAR},
for two Hookean viscoelastic fluids   meeting at a free boundary, the standard assumptions are that the velocity is continuous across the interface and that the jump in the normal stress is zero under \emph{ignoring the internal surface tension}.
This requires us to enforce the jump conditions
\begin{align}\label{201612262131}
\llbracket  v \rrbracket =0\mbox{ and }
\llbracket  \mathcal{S}  \rrbracket \vec{\nu} = 0 \; \mbox{ on } \Sigma(t),
\end{align}
where we have written the normal vector to the boundary $\Sigma(t)$ as $\vec{\nu}$, denoted the interfacial jump by
$\llbracket  f \rrbracket :=f_+|_{\Sigma(t)} -f_-|_{\Sigma(t)}$, and $f_\pm|_{\Sigma(t)}$ are the traces of the functions $f_\pm$ on $\Sigma(t)$. We mention that our results in Theorems \ref{thm:2109n}--\ref{thm:210902} can be extended to the case with internal surface tension, please refer to \cite{MR4762618} for details on how to further estimate the internal surface tension.
We will also enforce the condition that the fluid velocity vanishes at the fixed boundaries; we implement this via the boundary conditions
\begin{align}\label{201612262117}
v_\pm=0  \; \mbox{ on } \Sigma_\pm.
\end{align}

To simplify the representation of \eqref{0101f1} and \eqref{201612262117}, we introduce the indicator functions $\chi_{\Omega_\pm(t)}$ and denote
\begin{align}
& \rho:=\rho_+\chi_{\Omega_+(t)} +\rho_-\chi_{\Omega_-(t)},\  v:=v_+\chi_{\Omega_+(t)} +v_-\chi_{\Omega_-(t)},\
U=U_+\chi_{\Omega_+(t)} +U_-\chi_{\Omega_-(t)}, \label{201611092055} \\
& \mu:=\mu_+\chi_{\Omega_+(t)} +\mu_-\chi_{\Omega_-(t)},\
 \varsigma:=\varsigma_+\chi_{\Omega_+(t)} +\varsigma_-\chi_{\Omega_-(t)},\
\kappa:=\kappa_+\chi_{\Omega_+(t)} +\kappa_-\chi_{\Omega_-(t)} .\label{20161109205512}
\end{align}
Thus  one has
\begin{equation}\label{201611091004}
\begin{cases}
  \rho_t+\mm{div}(\rho{  v})=0& \mbox{in } \Omega(t),\\[0.5mm]
\rho v_t+\rho v\cdot\nabla v+\mm{div}\mathcal{S} =0\;\;&\mbox{in } \Omega(t),\\[0.5mm]
 {U}_{t}+ v\cdot\nabla {U}= \nabla v{U}&  \mbox{in } \Omega (t)
\end{cases}\end{equation}
and
\begin{align}\label{201611091004n}
v=0\; \mbox{ on } \Sigma_-^+,
\end{align}
where $ \mathcal{S}  $ is defined by \eqref{2020202109232002205}
with $(\rho,v,U, P, \mu, \varsigma)$ given by \eqref{201611092055}--\eqref{20161109205512}, and $\Sigma_-^+:=\Sigma_+\cup\Sigma_-$.
Moreover, under the first jump condition in \eqref{201612262131}, the internal surface function is defined by $v$, i.e.,
\begin{align}\label{201612262217}
d_t+v_1(x_{\mm{h}},d) \partial_1d+v_2(x_{\mm{h}},d) \partial_2d=v_3(x_{\mm{h}},d)  \mbox{ on }\mathbb{R}^2.
\end{align}

Finally, we impose the initial data for $(\rho, v, U,d)$:
\begin{align}
\label{201612262216}
(\rho, v, U )|_{t=0}:=(\rho^0,v^0, U^0 ) \mbox{ in } \Omega\!\!\!\!\!-\setminus \Sigma(0) \mbox{ and }
d|_{t=0}=d^0 \mbox{ on } \mathbb{R}^2,
\end{align}
where $\Omega\!\!\!\!-:=\mathbb{R}^2\times (h_-,h_+)$ and $\Sigma(0)=\{x_3=d(x_\mm{h},0)\}$. Then \eqref{201612262131}  and \eqref{201611091004}--\eqref{201612262216}
constitute an initial-boundary value problem for stratified compressible  Hookean viscoelastic fluids with an free interface, which we call the SCVF model for simplicity.

Now let us consider a  rest state of the SCVF model. Choose a constant $\bar{d}\in (h_-,h_+)$, and consider positive density constants $\bar{\rho}_\pm$, which   satisfy the  equilibrium state
 \begin{align}
\label{201611051547}
 \begin{cases}
\nabla P_\pm(\bar{\rho}_\pm)=0 &\mbox{in } \Omega_\pm,\\[0.5mm]
    \llbracket  P(\bar{\rho}) \rrbracket \mathbf{e}^3=0  &\mbox{on }\Sigma,
  \end{cases}
\end{align}
where $\mathbf{e}^3:=(0,0,1)^\top$,
$$\Sigma:=\mathbb{R}^2\times\{\bar{d}\},\ \Omega_+:=\mathbb{R}^2\times\{\bar{d}<x_3<h_+\}\mbox{ and }
\Omega_-:=\mathbb{R}^2\times\{h_-<x_3<\bar{d}\}.$$

Let $\bar{\rho}:=\bar{\rho}_+\chi_{\Omega_+} +\bar{\rho}_-\chi_{\Omega_-}$. Then  $(\rho,v,U)=(\bar{\rho},0,\mathbb{I})$ with $d=\bar{d}$ is
 an rest state (solution) of the SCVF model. Without loss of generality, we assume that $\bar{d}=0$ in this paper.
If $\bar{d}$ is not zero, we can adjust the $x_3$-coordinate to make $\bar{d}=0$.
Thus $h_-<0$, and $d$ can be called the displacement function of the point at the interface deviating from the plane $\Sigma$.

\subsection{Reformulation in Lagrangian coordinates}\label{sec:002}
It is well-known that the movement of the free interface $\Sigma(t)$ and the subsequent change of the
domains $\Omega_\pm(t)$ in Eulerian coordinates will result in severe mathematical difficulties.
Hence, we shall switch our verification to Lagrangian coordinates,
so that the interface and the domains are the fixed plane and the fixed domains, resp.
To this purpose, we take $\Omega_+$ and $\Omega_-$ to be the fixed Lagrangian domains, and assume that there exist
invertible mappings
\begin{equation*}
\zeta_\pm^0:\Omega_\pm\rightarrow \Omega_\pm(0),
\end{equation*}
such that $
\det(\nabla\zeta_\pm^0)\neq 0$,
\begin{align}\label{05261240}
\Sigma(0)=\zeta_\pm^0(\Sigma),\
\Sigma_+=\zeta_+^0(\Sigma_+)\mbox{ and }
\Sigma_-=\zeta_-^0(\Sigma_-).
\end{align}
The first condition in \eqref{05261240} means that the initial interface
$\Sigma(0)$ is parameterized by the mapping $\zeta_\pm^0$
restricted to $\Sigma$, while the latter two conditions in \eqref{05261240} mean that
$\zeta_\pm^0$ map the fixed upper and lower boundaries into
themselves. Define the flow maps $\zeta_\pm$ as the solutions to
\begin{align*}
 \begin{cases}
\partial_t \zeta_\pm(y,t)=v_\pm(\zeta_\pm(y,t),t)&\mbox{in }\Omega_\pm,
\\[0.5mm]
\zeta_\pm(y,0)=\zeta_\pm^0(y)&\mbox{in }\Omega_\pm.
                  \end{cases}
\end{align*}
We denote the Eulerian coordinates by $(x,t)$ with $x=\zeta(y,t)$,
whereas the fixed $(y,t)\in \Omega\times \mathbb{R}^+$ stand for the
Lagrangian coordinates.

In order to switch back and forth from Lagrangian to Eulerian coordinates, we assume that
$\zeta_\pm(\cdot ,t)$ are invertible and
$\Omega_{\pm}(t)=\zeta_{\pm}(\Omega_{\pm},t)$. Since $v_\pm$ and
$\zeta_\pm^0$ are all continuous across $\Sigma$, we have
$\Sigma(t)=\zeta_\pm(\Sigma,t)$, i.e.,
\begin{align}
\label{201701011211}
    \llbracket    \zeta     \rrbracket =0 \mbox{ on }\Sigma.
\end{align}
In other words, the Eulerian domains of upper and lower fluids are the images of $\Omega_\pm$
under the mappings $\zeta_\pm$, and the free interface is the image
of $\Sigma$ under the mappings $\zeta_\pm(\cdot,t)$.
In addition,  in view of the non-slip boundary condition $v_\pm|_{\Sigma_\pm}=0$,
we  have
\begin{align*}
y=\zeta_\pm(y, t)\mbox{ on }\Sigma_\pm.
\end{align*}
From now on, we define that $\zeta=\zeta_+\chi_{\Omega_+}+\zeta_-\chi_{\Omega_-}$, $\zeta^0:=\zeta_+^0\chi_{\Omega_+}+\zeta_-^0\chi_{\Omega_-}$ and $\eta:=\zeta-y$.

Next we introduce some notations involving $\eta$.  We define  $\mathcal{A}:=(\ml{A}_{ij})_{3\times 3}$ via
$\ml{A}^\top=(\nabla (\eta+y))^{-1}:=(\partial_j (\eta+y)_i)^{-1}_{3\times 3}$, and
the differential operators $\nabla_{\ml{A}}$, $\mm{div}_\ml{A}$ and $\Delta_\ml{A}$ are defined
as  follows:
\begin{align}
&\nabla_{\ml{A}}w:=(\nabla_{\ml{A}}w_1,\nabla_{\ml{A}}w_2,\nabla_{\ml{A}}w_3)^\top,\ \nabla_{\ml{A}}w_i:=(\ml{A}_{1k}\partial_kw_i,
\ml{A}_{2k}\partial_kw_i,\ml{A}_{3k}\partial_kw_i)^\top,\nonumber \\
\label{062209345}
&\mm{div}_{\ml{A}}(f_1,f_2,f_3)^\top=(\mm{div}_{\ml{A}}f_1,\mm{div}_{\ml{A}}f_2,\mm{div}_{\ml{A}}f_3)^\top,
\ \mm{div}_{\ml{A}}f_i:=\ml{A}_{lk}\partial_k f_{il},\\
& \Delta_{\mathcal{A}}w:= (\Delta_{\mathcal{A}}w_1,\Delta_{\mathcal{A}}w_2,\Delta_{\mathcal{A}}w_3)^\top \mbox{ and }
\Delta_{\mathcal{A}}w_i:=\mm{div}_{\ml{A}}\nabla_{\ml{A}}w_i\nonumber
  \end{align}
for  vector functions $w:=(w_1,w_2,w_3)^\top$ and $f_i:=(f_{i1},f_{i2},f_{i3})^\top$, where we have used the Einstein convention of summation over repeated indices and  $\partial_{k}$ denotes the partial derivative with respect to the $k$-th component of the variable $y$, i.e., $\partial_{k}:=\partial_{y_k}$.

Finally, we further introduce some properties of $\mathcal{A}$.
\begin{enumerate}[\quad\  (1)]
  \item In view of the definition of $\mathcal{A}$, one can deduce the following two important properties:
\begin{align}
\label{AklJ=0}
\partial_l (J\mathcal{A}_{kl})=0
\end{align}
and
\begin{align}\label{AklJdeta}
  \partial_i(\eta+y)_k\ml{A}_{kj}= \ml{A}_{ik}\partial_k(\eta+y)_j=\delta_{ij},
\end{align}
where $J=\det (\nabla (\eta+y))$, $\delta_{ij}=0$ for $i\neq j$, and $\delta_{ij}=1$ for $i=j$.
The relation \eqref{AklJ=0} is often called the geometric identity.
  \item We can evaluate that
$J\mathcal{A}\mathbf{e}^3 =\partial_1(\eta+y)\times \partial_2(\eta+y)$,
and thus, by \eqref{201701011211},
\begin{align}
\label{06051441}
    \llbracket  J\mathcal{A}\mathbf{e}^3     \rrbracket =0.
    \end{align}
  \item In view of \eqref{06051441}, the unit normal $\vec{n}$ to $\Sigma(t)=\zeta(\Sigma,t)$ can be written as follows:
\begin{align}
\label{05291021n}
&\vec{n}=  \tilde{n} |_{\Sigma} \mbox{ with }\tilde{n}:=\frac{J\mathcal{A}\mathbf{e}^3}{|J\mathcal{A}\mathbf{e}^3|}.
\end{align}
\end{enumerate}

In Lagrangian coordinates, the displacement gradient tensors $\tilde{U}_\pm(y,t)$ are  defined by a Jacobi matrix
of $\zeta_\pm(y,t)$:
\begin{align*}
\tilde{U}_\pm(y,t):=\nabla \zeta_\pm(y,t), \mbox{ i.e.},\ \tilde{U}_{ij}:=\partial_{j}(\zeta_\pm(y,t))_i.
\end{align*}
When we study this deformation tensor in Eulerian coordinates, we shall
denote it by
$$U_\pm(x,t):=\tilde{U}_\pm(\zeta^{-1}_\pm(x,t),t).$$
Moreover, by virtue of the chain rule, it is easy to check that $U_\pm(x,t)$ satisfy
the transport equation \eqref{0101f1}$_3$.
This means that the deformation tensor in Lagrangian coordinates
can be directly represented by $\zeta_{\pm}$, namely $U_{\pm}=\nabla\zeta_{\pm}|_{y=\zeta^{-1}_\pm(x,t)}$,  if the initial data $U^0$ also satisfies
\begin{align}\label{202108231510}
U_{\pm}^0=\nabla\zeta_{\pm}^0((\zeta^0_\pm)^{-1}(x,t),t).
\end{align}

Define  the Lagrangian unknowns
\begin{align*}
(\varrho,u)(y,t)=(\rho,v)(\zeta(y,t),t) \mbox{ for } (y,t)\in \Omega \times\mathbb{R}^+,
\end{align*}
then the transformed (stratified) viscoelastic problem
in Lagrangian coordinates reads as follows:
\begin{equation}\label{201911091004nn}
\begin{cases}
\zeta_t=u&\mbox{in } \Omega ,\\[0.5mm]
\varrho_t+\varrho\mm{div}_{\mathcal{A}}u=0&\mbox{in } \Omega ,\\[0.5mm]
\varrho u_t+\mm{div}_{\mathcal{A}}
 ({\mathcal{S}}_{\mathcal{A}}(\varrho,u )-\kappa(\nabla\zeta\nabla\zeta^\top/J-\mathbb{I}))=0
&\mbox{in } \Omega ,\\[0.5mm]
\llbracket  \zeta\rrbracket=\llbracket  u\rrbracket=0,\;\;\;
\llbracket ({\mathcal{S}}_{\mathcal{A}}(\varrho,u  )-\kappa(\nabla\zeta\nabla\zeta^\top/J-\mathbb{I}))\vec{n}\rrbracket=0
& \mbox{on } \Sigma,\\[0.5mm]
(\zeta, u)=(y,0) & \mbox{on }\Sigma_{-}^{+},\\[0.5mm]
 (\varrho,\zeta, u)|_{t=0}=(\varrho^0, \zeta^0, u^0) & \mbox{in }\Omega,
\end{cases}
\end{equation}
where we have defined that $\Omega:=\Omega_+\cup\Omega_-$,
\begin{align}
& {\mathcal{S}}_{\mathcal{A}} (\varrho,u )
:=P(\varrho)\mathbb{I}-\mathcal{V}_{\mathcal{A}}(u)
,\ \mathcal{V}_{\mathcal{A}}(u):=\mu\mathbb{D}_{\mathcal{A}}u+\lambda \mm{div}_{\mathcal{A}}u\mathbb{I},\ \lambda:=\varsigma-2\mu/3,\nonumber\\
&
\mathbb{D}_{\mathcal{A}}u:=\nabla_{\mathcal{A}} u+\nabla_{\mathcal{A}} u^\top
\;\mbox{ and }\;\vec{n}:=\frac{J\mathcal{A}\mathbf{e}^3}{|J\mathcal{A}\mathbf{e}^3|}
 \mbox{ represents the unit normal to } \Sigma(t) . \nonumber
\end{align}

Next, we proceed to rewrite the density and the elasticity in Lagrangian coordinates.
Viscoelasticity is a material property that exhibits both viscous and elastic characteristics when undergoing deformation.
In particular, a viscoelastic fluid strains when stretched, and quickly returns to its rest state
once the stress is removed. Therefore, we naturally have the following asymptotic behaviors after the rest state of the viscoelastic fluid being perturbed:
\begin{align}\label{202108231612}
\zeta\rightarrow y\mbox{ and }\varrho(y,t)\rightarrow\bar{\rho} \mbox{ as } t\rightarrow\infty.
\end{align}
It follows from \eqref{201911091004nn}$_1$ that
$J_t=J\mm{div}_{\mathcal{A}}u$,
which, together with \eqref{201911091004nn}$_2$, yields
\begin{align*}
\partial_t(\varrho J)=0.
\end{align*}
Thus we deduce from the above identity and the asymptotic behavior \eqref{202108231612} that
\begin{align}
\label{20202109250111564}
\varrho J=\varrho^0 J^0=\bar{\rho},
\end{align} which
implies $\varrho=\bar{\rho} J^{-1}$,
provided the initial data $(\varrho^0,J^0)$ satisfies
\begin{align}\label{202108231715}
\varrho^0=\bar{\rho} J_0^{-1},
\end{align}
where $\varrho^0$ and $J^0$ are the initial data of $\varrho$ and $J$, resp.
We call \eqref{20202109250111564} \emph{the asymptotic stability condition of density} in  Lagrangian coordinates.

It is worth noting that by \eqref{AklJ=0} and the relation \eqref{AklJdeta}, we have
\begin{align}\label{202108231732}
\mm{div}_{\mathcal{A}}(\nabla \zeta\nabla \zeta^\top/J)
=J^{-1}\mm{div}\left(J{\mathcal{A}}^\top\nabla \zeta\nabla \zeta^\top/J\right)
=J^{-1}\Delta\eta.
\end{align}
Similarly, by virtue of \eqref{AklJdeta}--\eqref{05291021n}, we get
\begin{align}\label{202108231740}
\llbracket \kappa(\nabla\zeta\nabla\zeta^\top/J-\mathbb{I})\vec{n}\rrbracket
=\llbracket \kappa(\nabla\zeta \mathbf{e}^3-J\mathcal{A}\mathbf{e}^3)\rrbracket/|J\mathcal{A}\mathbf{e}^3| \mbox{ on } \Sigma.
\end{align}

Under the assumptions \eqref{202108231510} and  \eqref{202108231715},
we can use the relations \eqref{20202109250111564}, \eqref{202108231732}  and \eqref{202108231740} to rewrite \eqref{201911091004nn} as follows:
\begin{equation}\label{20210824nn}
\begin{cases}
\eta_t=u&\mbox{in } \Omega ,\\
\bar{\rho} u_t+J\mm{div}_{\mathcal{A}}
{\mathcal{S}}_{\mathcal{A}}(\bar{\rho} J^{-1},u)=\kappa\Delta\eta
&\mbox{in } \Omega ,\\
\llbracket  \eta\rrbracket=\llbracket  u\rrbracket=0,\;\;
\llbracket{\mathcal{S}}_{\mathcal{A}}(\bar{\rho} J^{-1},u)J\mathcal{A}\mathbf{e}^3\rrbracket
=\llbracket \kappa((\nabla\eta+\mathbb{I}) \mathbf{e}^3-J\mathcal{A}\mathbf{e}^3)\rrbracket
& \mbox{on } \Sigma,\\
(\eta, u)=(0,0) & \mbox{on }\Sigma_{-}^{+},\\
(\eta, u)|_{t=0}=(\eta^0, u^0) & \mbox{in }\Omega,
\end{cases}
\end{equation}
where $\eta=\zeta-y $ and  ${\mathcal{S}}_{\mathcal{A}}(\bar{\rho} J^{-1},u):=P(\bar{\rho} J^{-1})\mathbb{I}-\mathcal{V}_{\mathcal{A}}(u)$. If the initial-boundary value problem \eqref{20210824nn} admits a solution $(\eta,u)$, then we further get the both functions of density and the deformation tensor in Lagrangian coordinates by the relations
\begin{align}\label{202108240820}
\varrho=\bar{\rho} J^{-1}  \mbox{ and } \tilde{U}=\nabla(\eta+y).
\end{align}

Due to the mathematical difficulty arising from $\kappa_{+}\neq\kappa_{-}$ on the interface, we also consider the case of $\kappa$ being uniform as in \cite{CCRMHJLWDH,LHWWZZF,chen2020linear},
i.e.,
\begin{align}
\kappa_{+}=\kappa_{-}.
\end{align}
Under such case, thanks to \eqref{06051441}, the  initial-boundary value problem \eqref{20210824nn} reduces to
\begin{equation}\label{20210824nnm}
\begin{cases}
\eta_t=u&\mbox{in } \Omega ,\\
\bar{\rho} u_t+J\mm{div}_{\mathcal{A}}
{\mathcal{S}}_{\mathcal{A}}(\bar{\rho} J^{-1},u)=\kappa\Delta\eta
&\mbox{in } \Omega ,\\[0.5mm]
\llbracket  \eta\rrbracket=\llbracket  u\rrbracket=0,\;\;
\llbracket{\mathcal{S}}_{\mathcal{A}}(\bar{\rho} J^{-1},u)J\mathcal{A}\mathbf{e}^3-\kappa\nabla\eta \mathbf{e}^3\rrbracket
=0
& \mbox{on } \Sigma,\\
(\eta, u)=(0,0) & \mbox{on }\Sigma_{-}^{+},\\
(\eta, u)|_{t=0}=(\eta^0, u^0) & \mbox{in }\Omega.
\end{cases}
\end{equation}

\subsection{Notations}\label{202411061605}
To conveniently stating our results and presenting their proofs, we shall introduce some simplified notations:

(1) Basic notations:
$\Omega\!\!\!\!-:=\mathbb{R}^2\times(h_{-},h_{+})$,
$\int:=\int_{\Omega}$, $\int_{\Sigma}:=\int_{\mathbb{R}^2}$.
$\mathbb{R}^+_0:=[0,\infty)$,
$a\lesssim b$ means that $a\leqslant cb$ for some positive constant $c$, which is independent of $\kappa$ but may depend on the domain $\Omega$ and the other known physical parameters, such as $\bar{\rho}_\pm$, $\mu_\pm$ and $\varsigma_\pm$, and may vary from line to line or from place to place. However sometimes we also renew to denote the constant $c$ by $C_i$ for $1\leqslant i\leqslant 3$ (or $c_j$ for $1\leqslant j\leqslant 5$), under such case, $C_i$ (or $c_j$) are still independent of $\kappa$ but do not vary from line to line or from place to place.  The notation $\partial_{\mm{h}}^{\alpha}$ denotes $\partial_1^{\alpha_1}\partial_2^{\alpha_2}$ for some multi-index of order
$\alpha:=(\alpha_1, \alpha_2)$.

(2) Simplified Banach spaces, norms and semi-norms:
\begin{equation*}
\begin{aligned}
&L^p:=L^p(\Omega), \ H^i :=W^{i,2}(\Omega),\ \|\cdot\|_{L^{p}}:=\|\cdot\|_{L^{p}(\Omega)} , \ \|\cdot\|_{i}:=\|\cdot\|_{H^{i}},\\
&\|\cdot\|_{i,j}^2:=\sum_{|\alpha|=i} \|\partial_{\mm{h}}^{\alpha}\cdot\|_{j}^2,\,\; \|\cdot\|^2_{\underline{i},j}:=\sum_{k=0}^i\|\cdot\|_{k,j}^2,\\
&|\cdot|_{s}:=\|\cdot\|_{H^{s}(\mathbb{R}^2)},\
|\cdot|_{L^{p}}:=\|\cdot\|_{L^{p}(\mathbb{R}^2)},
\end{aligned}
\end{equation*}
where $1<p\leqslant\infty$, $s$ is a real number, and $i$, $j$ are non-negative integers.
It should be noted that
$$|f|_{i+1/2}:=|f|_{\Omega_{+}}|_{i+1/2}+|f|_{\Omega_{-}}|_{i+1/2}, \;\mbox{ if } f\in H^{i+1}\mbox{ with }i\geqslant 0.$$
In addition, we also use the simplified norm
$\|\aleph(f_1,\ldots,f_n)\|_{X}:=
\sqrt{\sum_{1\leqslant i\leqslant n}\|\aleph f_i\|_{X}^2}$, where the number of components of the vector functions $f_1$, $\ldots$, $f_n$ may be different to each other, and $\aleph$ represents a differential operator.

(3) Simplified functional classes: for non-negative integer $i\geqslant1$,
\begin{equation*}
\begin{aligned}
&H^{1,i}:=\{f\in H^{i}~|~\partial_{\mm{h}}^{\alpha}f\in H^{i}\;\mbox{for}\;0\leqslant|\alpha|\leqslant1\},\\[0.5mm]
&H^{1,3}_{*}:=\{f\in H^{1,3}~|~\zeta:=f(y,t)+y:\Omega\rightarrow\Omega\;\mbox{is a $C^1$-homeomorphism},\\[0.5mm]
&\qquad \qquad \mbox{mapping, and satisfies }1/2\leqslant\mm{det}\nabla \zeta\leqslant3/2\},\\[0.5mm]
&H^{1,i}_{0}:=\{f\in H^{1,i}~|~f|_{\Sigma_{-}^{+}}=0\;\mbox{and}\; \llbracket f\rrbracket=0\},\quad
H^{1,3}_{*,0}:=H^{1,3}_{*}\cap H^{1,1}_{0}.
\end{aligned}
\end{equation*}

(4)  Energy/dissipation functionals:
\begin{align}
{\mathcal{E}}(t):=&\kappa^{-2}(\|u\|_{\underline{1},2}^2+\|u_{t}\|_{\underline{1},0}^2)+
{\kappa}^{-1}(\|\eta\|_{\underline{1},3}^2+\|u\|_{\underline{1},1}^2)+ \|\eta \|_{\underline{1},2}^2
  \nonumber\\[1mm]
&+\|(\mm{div}\eta,u)\|_{\underline{3},0}^2+(1+{\kappa})\|\eta\|_{\underline{3},1}^2,
\nonumber\\[1mm]
{\mathcal{D}}(t):=&
\kappa^{-2}(\|u\|_{\underline{1},3}^2+\|u_{t}\|_{\underline{1},1}^2)+ {\kappa}^{-1}\|u\|_{\underline{1},2}^2
+\|\eta\|_{\underline{1},3}^2+\kappa^{{1}/{2}}\|\eta \|_{\underline{1},2}^2
\nonumber\\[1mm]&
+\|u\|_{\underline{3},1}^2+\kappa\|\eta\|_{\underline{3},1}^2\nonumber.
\end{align}

\subsection{Global solvability in some classes of large data}\label{subsec:03}
Our first result concerns the existence of global solutions for the initial-boundary value problem \eqref{20210824nnm} in some classes of large data.
\begin{thm}\label{thm:2109n}
Let $(\eta^0, u^0)\in H^{1,3}_{*,0}\times H^{1,2}_{0}$. We further assume that $\kappa$ is a constant, i.e., $\kappa_-=\kappa_+$.
There exist positive constants $c_1$ (sufficiently small) and $c_2$, such that if initial data $(\eta^0, u^0)$ satisfies
the compatibility condition
\begin{align}\label{202409091654nn}
\llbracket{\mathcal{S}}_{\mathcal{A}^0}(\bar{\rho} J_0^{-1},u^0)J^0\mathcal{A}^0\mathbf{e}^3-\kappa\nabla\eta^0 \mathbf{e}^3\rrbracket=0
 \mbox{ on } \Sigma
\end{align}
and the condition of large elastic coefficient
\begin{align}\label{202109091654}
\kappa^{-1}\max\big\{ \big(4c_2(1+{\kappa}^{-8}){E}^0 \big)^{1/14}, \big(4c_2(1+{\kappa}^{-8}){E}^0 \big)^{8}\big\}\leqslant  {c_1} ,
\end{align}
then  the initial-boundary value  problem \eqref{20210824nnm} admits a unique strong solution
$(\eta, u)\in C^{0}(\mathbb{R}^+_0,$ $H^{1,3}_{*,0} \times H^{1,2}_{0})$;
moreover, it enjoys the stability estimate
\begin{align}
& \mathcal{E} (t) + \int_{0}^{t} {\mathcal{D}}(\tau)e^{\frac{(\tau-t)}{c_3 (1+\kappa^{-1})}}\mm{d}\tau
\lesssim e^{-\frac{t}{c_3 (1+\kappa^{-1})}}(1+\kappa^{-8}){E}^0  \label{202109090810}
\end{align}
for some positive constant $c_3>0$,
where we have defined that
\begin{align*}
E^0:=&\kappa^{-2}\|u^0\|_{\underline{1},2}^2+
\kappa^{-1}(\|\eta^0\|_{\underline{1},3}^2+\|u^0\|_{\underline{1},1}^2)+ \|\eta^0 \|_{\underline{1},2}^2
  \nonumber\\[1mm] &
+\|(\mm{div}\eta^0,u^0)\|_{\underline{3},0}^2+(1+{\kappa})\|\eta^0\|_{\underline{3},1}^2.
\end{align*}
\end{thm}
\begin{rem}
 Recalling that the solution $\eta$ in Theorem \ref{thm:2109n} satisfies
$\zeta(y,t):=\eta(y,t)+y\in H_*^{1,3} $  for each $t\in \mathbb{R}_0^+$,
we can easily obtain an existence result of strong solutions to the original problem, i.e., the SCVF model,
by an inverse transformation of Lagrangian coordinates and taking  $c_1$ sufficiently small, where
$\Omega_+(t)$, $\Omega_-(t)$ and  $\Sigma(t) $ are defined by \eqref{20224102301821}--\eqref{202410301926} with $d:=\zeta_3(\zeta_{\mm{h}}^{-1}(x_{\mm{h}},t), 0,t)$, resp. \cite{JFJHJS2023}. In addition, let
$$\alpha= \kappa^{-2}\|u^0\|_{\underline{1},2}^2+{\kappa}^{-1}\| u^0\|_{\underline{1},1}^2+\|u^0\|_{\underline{3},0}^2
,$$
then the condition of large elasticity coefficient in \eqref{202109091654} with $\eta^0=0$ reduces to
\begin{align}
&\max \bigg\{ \big(4c_2(1+\kappa^{-8})\alpha \big)^{1/14}, \big(4c_2(1+{\kappa}^{-8})\alpha \big) ^8\bigg\}
 \leqslant c_1\kappa,
 \end{align}
which shows that the singularities do not form on the internal free boundary $\Sigma(t)$ of the compressible Hookean viscoelastic fluid,
if the elasticity coefficient is \emph{relatively} larger than the initial perturbation velocity.
\end{rem}

The basic idea in the proof for Theorem \ref{thm:2109n} can be found in  \cite{JFJSGS}, where Jiang--Jiang investigated the well-posedness of the \emph{2D} free-boundary problem of incompressible viscoelastic fluids in some classes of large data. Next we briefly sketch the proof idea of Theorem \ref{thm:2109n} and explain why Jiang--Jiang's result in \cite{JFJSGS} can be extended to the case of \emph{3D} stratified compressible viscoelastic fluids in this paper.

Taking the inner product of \eqref{20210824nnm}$_2$ with $u$ in $L^2$,
we obtain the basic energy identity of the initial-boundary value problem \eqref{20210824nnm}:
\begin{align}
&\frac{\mm{d}}{\mm{d}t}\left( \frac{1}{2}\| \bar{\rho} u\|_0^2+
 {\bar{\rho}}\int
\int_{\bar{\rho}/4}^{\bar{\rho}J^{-1}}\frac{P(z)}{z^2}\mm{d}z\mm{d}y
+\frac{\kappa}{2}\|\nabla \eta\|_0^2\right)\nonumber \\
&\label{202108asda301420}
+\frac{1}{2}\|\sqrt{\mu J}\mathbb{D}_{\mathcal{A}}u\|_0^2+\|\sqrt{\lambda J}\mm{div}_{\mathcal{A}}u\|_0^2
=0 .
\end{align}
Integrating the above identity over $(0,t)$ yields that
\begin{align}
&  \| \bar{\rho} u\|_0^2+2
 {\bar{\rho}}\int
\int_{\bar{\rho}/4}^{\bar{\rho}J^{-1}}\frac{P(z)}{z^2}\mm{d}z\mm{d}y
+ {\kappa} \|\nabla \eta\|_0^2 + \|\sqrt{\mu J}\mathbb{D}_{\mathcal{A}}u\|_0^2+2\|\sqrt{\lambda J}\mm{div}_{\mathcal{A}}u\|_0^2\nonumber \\
&\label{202108asda301422}
= \| \bar{\rho} u^0\|_0^2+2
 {\bar{\rho}}\int
\int_{\bar{\rho}/4}^{\bar{\rho}J^{-1}_0}\frac{P(z)}{z^2}\mm{d}z\mm{d}y
+ {\kappa} \|\nabla \eta^0\|_0^2  =:2I^0 .
\end{align}
In particular, we have
\begin{align}
\|\nabla \eta\|_0^2 \leqslant 2 I^0/\kappa.
\label{2022410302122}
\end{align}
We call $I^0$ the initial mechanical energy. It is easy to see that, for the given initial mechanical energy $I^0$,
$\|\nabla \eta\|_0^2 $ is enough small as $\kappa$ is sufficiently large, as shown by \eqref{2022410302122}.

Such a smallness property also exists in Jiang--Jiang's result for the free-boundary problem of  \emph{2D} incompressible viscoelastic fluids  in \cite{JFJSGS} and also plays a key role in establishing the existence of global solutions with some classes of large data. However the norm of pressure in \cite{JFJSGS} increases with respect to $\kappa$. To overcome this growth behavior, Jiang--Jiang used the three-layer energy method, which includes anisotropic energy/dissipation functionals with the weight of elasticity coefficient, to establish
\emph{a priori} stability estimates under proper \emph{a priori} assumptions. In our compressible case, the perturbation pressure does not increase with respect to $\kappa$, see the expression of pressure $P(\bar{\rho} J^{-1})$ in \eqref{201903120816}. Based on this fact and Jiang--Jiang's  analysis process for the 2D incompressible problem, we roughly observe that the problem  \eqref{20210824nnm}
may be approximated by the corresponding linear problem (see \eqref{20210824nnmmn}) with the smallness property
$$
\|\nabla\eta\|_{L^{\infty}}\ll 1\mbox{ for sufficiently large }\kappa.
 $$

Motivated by the approximate idea mentioned above and Jiang--Jiang's proof method in \cite{JFJSGS}, we carefully analyze the corresponding linear problem and  the structure of the nonlinear terms in
\eqref{20210824nnmm}. Based on this analysis, we choose the following \emph{a priori} assumptions:
\begin{align}\label{prio1}
\;\;\sup_{t\in[0,T]}\mathcal{E}(t)\leqslant K^2 \mbox{ for some } T>0
\end{align}
and
\begin{align}\label{prio2}
\kappa^{-1}\max\left\{K^{1/7}, K^{16}\right\}\in(0,\delta]\mbox{ with sufficiently small }\delta.
\end{align}The smallness of $\delta$ may only depends on the parameters $\mu$, $\lambda$ and $\bar{\rho}$.
Consequently, we can utilize an energy method to prove that there exist two constants $K$ and $\delta$, such that
the solution of \eqref{20210824nnm}  with \emph{a priori} assumptions \eqref{prio1} and \eqref{prio2} enjoys the \emph{a priori} estimate
\begin{align}\label{pr}
\sup_{t\in[0,T]}\mathcal{E}(t)\leqslant K^2/4,
\end{align}
where the constant $K$ will be determined later by \eqref{20241101031319}.

Roughly speaking, our energy method can be divided into two steps: the first step is to establish the tangential estimates,
and the second step is to close the normal estimates for $(\eta,u)$ using the tangential estimates.
To achieve this goal,
we shall make full use of the anisotropic Sobolev inequalities (see \eqref{202408241430} and \eqref{202108261428nm}) to carefully estimate the nonlinear terms ${\mathcal{N}}^{\mu}_{u}$, ${\mathcal{N}}^{\lambda}_{u}$ and ${\mathcal{N}}^{P}_{\eta}$ in the inhomogeneous form \eqref{20210824nnmm}. Based on the \emph{a priori} estimate mentioned above, together with the existence of unique local solutions, we can immediately obtain Theorem \ref{thm:2109n}. The detailed proof will be presented in Section \ref{20201010052127}.

It should be noted that the key inequality \eqref{2022410302122} can not be excepted to the case $\kappa_-\neq \kappa_+$. In fact, under such case, we have
the basic energy identity:
\begin{align}
&
\frac{\mm{d}}{\mm{d}t}\left( \frac{1}{2}\| \bar{\rho} u\|_0^2+
 {\bar{\rho}}\int
\int_{\bar{\rho}/4}^{\bar{\rho}J^{-1}}\frac{P(z)}{z^2}\mm{d}z\mm{d}y
+\frac{1}{4}\big(\|\sqrt{\kappa}\mathbb{D} \eta\|_0^2-2\|\sqrt{\kappa}\mm{div}\eta\|_0^2\big)\right)\nonumber \\
&
+\frac{1 }{2}\|\sqrt{\mu J}\mathbb{D}_{\mathcal{A}}u\|_0^2+ \|\sqrt{\lambda J}\mm{div}_{\mathcal{A}}u\|_0^2= \int  \kappa  \tilde{\mathcal{B}}^{{N}}: \nabla u\mm{d}y ,\label{202108asdsafa301420}
\end{align}
see the definition of $\tilde{\mathcal{B}}^{{N}}$ in \eqref{202502062146}. The additional nonlinear term $\kappa\tilde{\mathcal{B}}^{{N}}$ results in the failure of our method to establish Theorem \ref{thm:2109n}. This is also the reason why we shall consider the stratified case in Theorem \ref{thm:2109n}.

\subsection{Vanishing phenomena of the nonlinear interactions}
Now we state the vanishing phenomena of the nonlinear interactions with respect to large $\kappa$.  Such phenomena present that the solution of the initial-boundary value problem \eqref{20210824nnm} in Theorem \ref{thm:2109n} with $\eta^0=0$  can be approximated by the solution  of a linear problem, if $\kappa$ is sufficiently large.
\begin{thm}\label{thm:210902}
Let $(\eta,u)$ be the solution established in Theorem \ref{thm:2109n}, then we have
the stability estimates of the nonlinear solution $(\eta,u)$  around the linear solution $(\eta^{\mm{l}},u^{\mm{l}})$:
\begin{align}\label{202109061353}
&\|(\mm{div}\eta^{\mm{d}}, u^{\mm{d}})(t)\|_{\underline{3},0}^2+(1+\kappa) \|\eta^{\mm{d}}(t)\|_{\underline{3},1}^2
+ \int_{0}^{t} \big( \| \mm{div}\eta^{\mm{d}}(\tau)\|_{\underline{3},0}^2+
\|(\sqrt{\kappa}\eta^{\mm{d}}, u^{\mm{d}})(\tau)\|_{\underline{3},1}^2\big)e^{\frac{(\tau-t)}{c_4(1+\kappa^{-1})}}\mm{d}\tau
\nonumber \\
&\lesssim (1+\kappa^{-17})(\kappa^{-1/4}+\kappa^{-5/4})  e^{-\frac{t}{c_4(1+\kappa^{-1})}}{E}^0({E}^0+\|u^0\|_{\underline{1},2}^2),
\end{align}
and
\begin{align}
&\|\eta^{\mm{d}}(t)\|_{\underline{1},2}^2
+\kappa \int_{0}^{t} \|\eta^{\mm{d}}(\tau)\|_{\underline{1},2}^2e^{\frac{(\tau-t)}{c_5(1+\kappa^{-1})}}\mm{d}\tau\nonumber \\
&\lesssim (1+\kappa^{-17})(\kappa^{-1/8}+\kappa^{-5/4})  e^{-\frac{t}{c_5(1+\kappa^{-1})}}(1+{E}^0)
({E}^0 +\|u^0\|_{\underline{1},2}^2)
\label{202409061553}
\end{align}
for positive constants $c_4$ and $c_5$, where we have defined that $(\eta^{\mm{d}},u^{\mm{d}}):=(\eta-\eta^{\mm{l}},u-u^{\mm{l}})$
with $(\eta^{\mm{l}},u^{\mm{l}})\in C^{0}(\mathbb{R}^+,H^{1,3}_{0} \times H^{1,2}_{0})$ is the solution to the following linear problem
\begin{equation}\label{20210824nnmmn}
\begin{cases}
\eta_t^{\mm{l}}=u^{\mm{l}}&\mbox{in } \Omega ,\\[0.5mm]
\bar{\rho} u_t^{\mm{l}}
-\mm{div}\left({P}'(\bar{\rho})\bar{\rho}\mm{div}\eta^{\mm{l}} \mathbb{I}+\mu\mathbb{D}u^{\mm{l}}+\lambda\mm{div}u^{\mm{l}}\mathbb{I}\right)
=\kappa\Delta\eta^{\mm{l}}
&\mbox{in } \Omega ,\\[0.5mm]
\llbracket  \eta^{\mm{l}}\rrbracket=\llbracket  u^{\mm{l}}\rrbracket=0,
& \mbox{on } \Sigma,\\[0.5mm]
\llbracket {P}'(\bar{\rho})\bar{\rho}\mm{div}\eta^{\mm{l}}\mathbb{I}+\mu\mathbb{D}(u^{\mm{l}})+\lambda\mm{div}u^{\mm{l}} \mathbb{I}+\kappa\nabla\eta^{\mm{l}} \rrbracket \mathbf{e}^3
=0
& \mbox{on } \Sigma,\\[0.5mm]
(\eta^{\mm{l}}, u^{\mm{l}})=(0,0) & \mbox{on }\Sigma_{-}^{+},\\[0.5mm]
(\eta^{\mm{l}}, u^{\mm{l}})|_{t=0}=(\eta^0, u^0+u^{\mathrm{r}}) & \mbox{in }\Omega
\end{cases}
\end{equation}
and $(\eta^{\mm{l}},u^{\mm{l}})$ enjoys
\begin{align}\label{240504122nm}
\mbox{the same regularity as well as the solution $(\eta,u)$ in Theorem \ref{thm:2109n}}.
\end{align}
Moreover,  the function $u^{\mathrm{r}}\in {H}_{0}^{1,2}$ enjoys
\begin{align}\label{202408191959}
\|u^{\mathrm{r}}\|_{\underline{1},2}^2
\lesssim\|\eta^0\|_{\underline{2},1}\|\eta^0\|_{\underline{1},2}
\|(\eta^0, u^0)\|_{\underline{1},2}^2.
\end{align}
\end{thm}
Here and in what follows, we call the strong solution of the initial-boundary value problem \eqref{20210824nnm} (abbr. the nonlinear problem) the \emph{nonlinear solution}, and the strong solution of  the linear problem \eqref{20210824nnmmn}  the \emph{linear solution}. We also briefly sketch the proof idea of Theorem \ref{thm:210902}. Since both the linear and nonlinear solutions may enjoy different compatibility conditions, we shall first use the stratified elliptic theory to modify the given initial data $(\eta^0,u^0)$ in Theorem \ref{thm:2109n} so that the obtained new initial data can be used to generate a linear solution. Then we subtract the linear problem from the nonlinear problem to obtain an error problem. Following the argument of Theorem \ref{thm:2109n} with slight modifications, we can easily derive the stability estimate of the error function $(\eta^{\mm{d}},u^{\mm{d}})$ from the error problem. The  detailed derivation will be  presented in Section \ref{asymptotic}.

\vspace{3mm}
The rest of this paper is organized as follows.
In Section \ref{20201010101708}, we will further rewrite \eqref{20210824nnm} as an inhomogeneous form,
which is convenient for establishing the \emph{a priori} estimates of solutions.
Section \ref{preliminaries} is devoted to the derivation of some preliminary estimates which will
be utilized in the energy evolution.
In Sections \ref{20201010052127} and \ref{asymptotic}, we provide the detailed proofs for Theorems \ref{thm:2109n} and \ref{thm:210902}, resp.

\section{Inhomogeneous forms}\label{20201010101708}
In this section, we rewrite the initial-boundary value problem \eqref{20210824nnm} as an inhomogeneous form, which will be used
in the derivations of \emph{a priori} estimates for the solutions of \eqref{20210824nnm}.

To begin with,
we let ${\mathcal{B}} =J{\mathcal{A}}$,
$\tilde{\mathcal{B}} ={\mathcal{B}}-\mathbb{I}$,
 $$  \tilde{\mathcal{B}}^{{L}} =\left(\begin{array}{ccc}
\partial_2\eta_2+\partial_3\eta_3 &
-\partial_1\eta_2 &
- \partial_1\eta_3 \\
-\partial_2\eta_1 &
\partial_1\eta_1+\partial_3\eta_3 &
- \partial_2\eta_3 \\
-\partial_3\eta_1 &
- \partial_3\eta_2 &
\partial_1\eta_1+\partial_2\eta_2
        \end{array}\right)
$$
and
\begin{align}\tilde{\mathcal{B}}^{{N}}=\left(\begin{array}{ccc}
\partial_2\eta_2\partial_3\eta_3-\partial_2\eta_3\partial_3\eta_2 &
\partial_1\eta_3\partial_3\eta_2-\partial_1\eta_2\partial_3\eta_3 &
\partial_1\eta_2\partial_2\eta_3-\partial_1\eta_3\partial_2\eta_2    \\
\partial_2\eta_3\partial_3\eta_1-\partial_2\eta_1\partial_3\eta_3 &
\partial_1\eta_1\partial_3\eta_3-\partial_1\eta_3\partial_3\eta_1 &
\partial_1\eta_3\partial_2\eta_1-\partial_1\eta_1\partial_2\eta_3    \\
\partial_2\eta_1\partial_3\eta_2-\partial_2\eta_2\partial_3\eta_1 &
\partial_1\eta_2\partial_3\eta_1-\partial_1\eta_1\partial_3\eta_2 &
\partial_1\eta_1\partial_2\eta_2-\partial_1\eta_2\partial_2\eta_1
        \end{array}\right).
        \label{202502062146}
\end{align}
Then
\begin{align}
\label{2020201009270202218}
\tilde{\mathcal{B}} =\tilde{\mathcal{B}}^{{L}}+\tilde{\mathcal{B}}^{{N}}\end{align}
  and
\begin{align*}\mathcal{A}=J^{-1}\mathcal{B}=J^{-1}(\tilde{\mathcal{B}}+\mathbb{I})
.\end{align*}

Using the above identity, we have
\begin{align}
&\label{202108241105}
\mathbb{D}_{\mathcal{A}}u=\mathbb{D}u+\tilde{\mathcal{N}}^{1}_{u},\\[1mm]
&\label{202108241112}
\mm{div}_{\mathcal{A}}u=\mm{div}u+\tilde{\mathcal{N}}^{2}_{u},
\end{align}
where we have defined that
\begin{align*}
&\tilde{\mathcal{N}}^{1}_{u}:=\left(\tilde{\mathcal{N}}^{1,u}_{i,j}\right)_{3\times3}\;\mbox{with}\;\;\tilde{\mathcal{N}}^{1, u}_{ij}:=J^{-1}(\tilde{\mathcal{B}}_{jl}\partial_lu_i
+\tilde{\mathcal{B}}_{il}\partial_{l}u_j)+(J^{-1}-1)\left(\partial_ju_i+\partial_{i}u_j\right),\\[1mm]
&\tilde{\mathcal{N}}^{2}_{u}:=J^{-1}\tilde{\mathcal{B}}_{ij}\partial_ju_i+(J^{-1}-1)\partial_iu_i.
\end{align*}

Exploiting the Taylor expansion, it is easy to check that
\begin{align} \label{201903120816}
P(\bar{\rho}J^{-1})=P(\bar{\rho})-{P}'(\bar{\rho})\bar{\rho}\mm{div}\eta+\tilde{\mathcal{N}}^{3}_{\eta},
\end{align}
where we have defined that
\begin{align}
&\tilde{\mathcal{N}}^{3}_{\eta}: ={P}'(\bar{\rho})\bar{\rho}(J^{-1}-1+\mm{div}\eta)
+\mathcal{R}^{\eta},\nonumber \\[1mm]
&\mathcal{R}^{\eta}:= \int_{0}^{\bar{\rho}(J^{-1}-1)}(\bar{\rho}
(J^{-1}-1)-z)\frac{\mm{d}^2}{\mm{d}z^2} P (\bar{\rho}+z)\mm{d}z. \label{20202022101102103}
\end{align}
Consequently,
\begin{align}
&\mathcal{S}_{\mathcal{A}}(\bar{\rho}J^{-1},u)J\mathcal{A} \nonumber \\[1mm]
&=\left(P(\bar{\rho})-{P}'(\bar{\rho})\bar{\rho}\mm{div}\eta+\tilde{\mathcal{N}}^{3}_{\eta}\right)\mathcal{B}
-\mu\left(\mathbb{D}u+\tilde{\mathcal{N}}^{1}_{u}\right)\mathcal{B}
-\lambda\left(\mm{div}u \mathbb{I}+\tilde{\mathcal{N}}^{2}_{u}\right)\mathcal{B} \nonumber \\[1mm]
&=P(\bar{\rho})\mathcal{B}
-\left({P}'(\bar{\rho})\bar{\rho}\mm{div}\eta \mathbb{I}+\mu\mathbb{D}u+\lambda\mm{div}u\mathbb{I}\right)
+{\mathcal{N}}^{P}_{\eta}-{\mathcal{N}}^{\mu}_{u}-{\mathcal{N}}^{\lambda}_{u},
\label{20201101020111209}
\end{align}
where we have defined that
\begin{align*}
&{\mathcal{N}}^{P}_{\eta}
:=  \tilde{\mathcal{N}}^{3}_{\eta} {\mathcal{B}} -{P}'(\bar{\rho})\bar{\rho}\mm{div}\eta\tilde{\mathcal{B}}
\\[1mm]
&{\mathcal{N}}^{\mu}_{u}:= \mathbb{D}u\tilde{\mathcal{B}}
+\tilde{\mathcal{N}}^{1}_{u} \mathcal{B} \mbox{ and }{\mathcal{N}}^{\lambda}_{u}:=
 \mm{div}u\tilde{\mathcal{B}} +
 \tilde{\mathcal{N}}^{2}_{u}{\mathcal{B}} .
\end{align*}
Making use of \eqref{201611051547}$_2$, \eqref{AklJ=0}, \eqref{06051441} and \eqref{20201101020111209}, we further have
\begin{align}
& J\mm{div}_{\mathcal{A}}{\mathcal{S}}_{\mathcal{A}}(\bar{\rho} J^{-1},u)
= \mm{div}\left( {\mathcal{N}}^{P}_{\eta}-{\mathcal{N}}^{\mu}_{u}
-{\mathcal{N}}^{\lambda}_{u} -({P}'(\bar{\rho})\bar{\rho}\mm{div}\eta \mathbb{I}+\mu\mathbb{D}u+\lambda\mm{div}u\mathbb{I})\right)
\label{202108241348}
\end{align}
and\begin{align}\label{202108241257}
\llbracket{\mathcal{S}}_{\mathcal{A}}(\bar{\rho} J^{-1},u)J\mathcal{A}\mathbf{e}^3\rrbracket
=
\llbracket \mathcal{N}^{P}_{\eta} -{\mathcal{N}}^{\mu}_{u}-{\mathcal{N}}^{\lambda}_{u}\rrbracket \mathbf{e}^3
-\llbracket({P}'(\bar{\rho})\bar{\rho}\mm{div}\eta\mathbb{I}+\mu\mathbb{D}u+\lambda\mm{div}u\mathbb{I})\rrbracket \mathbf{e}^3
.
\end{align}
Consequently, thanks to \eqref{202108241348} and \eqref{202108241257}, the initial-boundary value problem \eqref{20210824nnm} can be rewritten as the following inhomogeneous form:
\begin{equation}\label{20210824nnmm}
\begin{cases}
\eta_t=u&\mbox{in } \Omega ,\\[0.5mm]
\bar{\rho} u_t
-\mm{div}\left({P}'(\bar{\rho})\bar{\rho}\mm{div}\eta \mathbb{I}+\mu\mathbb{D}u+\lambda\mm{div}u\mathbb{I}
+\kappa\nabla\eta^\top\right)
=\mm{div}({\mathcal{N}}^{\mu}_{u}+{\mathcal{N}}^{\lambda}_{u}-{\mathcal{N}}^{P}_{\eta})
&\mbox{in } \Omega ,\\[0.5mm]
\llbracket  \eta\rrbracket=\llbracket  u\rrbracket=0,
& \mbox{on } \Sigma,\\ \llbracket
{P}'(\bar{\rho})\bar{\rho}\mm{div}\eta \mathbb{I}+
\mu\mathbb{D}u+ \lambda\mm{div}u
\mathbb{I}+\kappa\nabla\eta \rrbracket \mathbf{e}^3
=\llbracket{\mathcal{N}}^{P}_{\eta}-{\mathcal{N}}^{\mu}_{u}-{\mathcal{N}}^{\lambda}_{u}\rrbracket \mathbf{e}^3
& \mbox{on } \Sigma,\\[0.5mm]
(\eta, u)=(0,0) & \mbox{on }\Sigma_{-}^{+},\\[0.5mm]
(\eta, u)|_{t=0}=(\eta^0, u^0) & \mbox{in }\Omega.
\end{cases}
\end{equation}

\section{Preliminaries}\label{preliminaries}
Before the proof of Theorem \ref{thm:2109n},
we shall first derive some preliminary estimates and then establish the estimates of the nonlinear terms ${\mathcal{N}}^{\mu}_{u}$, ${\mathcal{N}}^{\lambda}_{u}$ and $\mathcal{N}^{P}_{\eta}$. These estimates of nonlinear terms will be utilized in the derivations of energy evolution.
To this purpose,
we let $(\eta,u )$ be a solution to the problem \eqref{20210824nnm} defined on $\Omega\!\!\!\!\!-\times [0,T]$ with $T>0$
 and \emph{a priori} assume that
\begin{align}
&\label{202108261502}
(1+\|\nabla\eta\|_{\underline{3},0}^{1/2}\|\nabla\eta\|_{\underline{1},2}^{1/2})
 \|\nabla\eta\|_{\underline{2},0}^{1/2}\|\nabla\eta\|_{\underline{1},1}^{1/4}\|\nabla\eta\|_{\underline{2},1}^{1/4}
\lesssim  {\delta}.
\end{align}
It should be noted that the smallness of $\delta$ depends on $\bar{\rho}$, $\mu$, $\lambda$ and $\Omega$, but not on $\kappa$.
Next we use the condition \eqref{202108261502} to establish some preliminary estimates.
\begin{lem}\label{lem:2108261504}
Let $(\eta,u)$ satisfy $\eta_t=u$. Under the condition \eqref{202108261502} with sufficiently small $\delta$, we have
\begin{enumerate}
\item[(1)] the estimates for $\mathcal{B}$:
\begin{align}
&\label{202109071530}
\|\tilde{\mathcal{B}}\|_{\underline{i},0}\lesssim \|\nabla\eta\|_{\underline{i},0},\\[0.5mm]
 &\label{202108270837}
\|\tilde{\mathcal{B}} \|_{\underline{j},1}\lesssim \|\nabla\eta\|_{\underline{j},1},\\[0.5mm]
 &\label{202108270845}
\|\tilde{\mathcal{B}} \|_{\underline{1},2}\lesssim \|\nabla\eta\|_{\underline{1},2},\\[0.5mm]
 &\label{202108270845nm}
\|\tilde{\mathcal{B}}_{t} \|_{\underline{1},0}\lesssim \|\nabla u\|_{\underline{1},0},\\[0.5mm]
 &\label{202108270845n12m}
\|\tilde{\mathcal{B}}_{t} \|_{L^{\infty}}\lesssim \|\nabla u\|_{L^{\infty}},
\end{align}
where    $0\leqslant i\leqslant3$ and $1\leqslant j\leqslant2$ and $\tilde{\mathcal{B}}$ is defined by \eqref{2020201009270202218}.
\item[(2)]
the estimates for $J^{-1}$ and\; $\mathcal{R}^{\eta}$:
\begin{align}\label{202109071535}
& 1/2\leqslant J^{-1}\leqslant 2 ,\\[0.5mm]
&\label{202109071532}
 \|J^{-1}-1\|_{\underline{i},0}\lesssim \|\nabla\eta\|_{\underline{i},0},\\[0.5mm]
 &\label{202108270852}
\|J^{-1}-1\|_{\underline{j},1}\lesssim \|\nabla\eta\|_{\underline{j},1},\\[0.5mm]
 &\label{202108270852nm}
\|J^{-1}-1\|_{\underline{1},2}\lesssim \|\nabla\eta\|_{\underline{1},2},\\[0.5mm]
 &\label{202108270852n12}
\|J^{-1}_{t} \|_{\underline{1},0}\lesssim \|\nabla u\|_{\underline{1},0},\\[0.5mm]
 &\label{202108270852n12nm}
\|J^{-1}_{t} \|_{L^{\infty}}\lesssim \|\nabla u\|_{L^{\infty}},\\[0.5mm]
&\label{202108270856}
 \|(J^{-1}-1+\mm{div}\eta,\mathcal{R}^{\eta})\|_{\underline{i},0}
\lesssim
\|\nabla\eta\|_{\underline{3},0}^{1/2}\|\nabla\eta\|_{\underline{1},1}^{1/4}\|\nabla\eta\|_{\underline{1},2}^{1/4}
\|\nabla\eta\|_{\underline{i},0},\\[0.5mm]
&\label{202108270854}
\|(J^{-1}-1+\mm{div}\eta,\mathcal{R}^{\eta})\|_{\underline{j},1}
\lesssim
\|\nabla\eta\|_{\underline{3},0}^{1/2}\|\nabla\eta\|_{\underline{1},1}^{1/4}\|\nabla\eta\|_{\underline{1},2}^{1/4}
\|\nabla\eta\|_{\underline{j},1},\\[0.5mm]
&\label{202108270855}
\|(J^{-1}-1+\mm{div}\eta,\mathcal{R}^{\eta})\|_{\underline{2},1}
\lesssim
\|\nabla\eta\|_{\underline{3},0}^{1/2}\|\nabla\eta\|_{\underline{1},1}^{1/4}\|\nabla\eta\|_{\underline{1},2}^{1/4}
\|\nabla\eta\|_{\underline{2},1},\\[0.5mm]
&\label{202108270855n12}
\|(J^{-1}-1+\mm{div}\eta,\mathcal{R}^{\eta})_{t}\|_{\underline{1},0}
\lesssim
\|\nabla\eta\|_{\underline{3},0}^{1/2}\|\nabla\eta\|_{\underline{1},1}^{1/4}\|\nabla\eta\|_{\underline{1},2}^{1/4}
\|\nabla u\|_{\underline{1},0},\\[0.5mm]
&\label{202108270855n12nm}
\|(J^{-1}-1+\mm{div}\eta,\mathcal{R}^{\eta})_{t}\|_{L^{\infty}}
\lesssim
\|\nabla\eta\|_{\underline{3},0}^{1/2}\|\nabla\eta\|_{\underline{1},1}^{1/4}\|\nabla\eta\|_{\underline{1},2}^{1/4}
\|\nabla u\|_{L^{\infty}},
\end{align}
where $0\leqslant i\leqslant3$, $1\leqslant j\leqslant2$, $J:=\det(\nabla \eta+I)$ and $\mathcal{R}^{\eta}$ is defined in \eqref{20202022101102103}.
\end{enumerate}
\end{lem}
\begin{pf}
(1)
Let $1\leqslant$ $k$, $l$, $m, n\leqslant 3$.
By \eqref{202108261428nm}, we have
\begin{align*}
&\|\partial_{m}\eta_{n}\partial_{k}\eta_{l}\|_{\underline{i},0}
\lesssim
\|\nabla\eta\|_{\underline{2},0}^{1/2}\|\nabla\eta\|_{\underline{1},1}^{1/4}
\|\nabla\eta\|_{\underline{2},1}^{1/4}\|\nabla\eta\|_{\underline{i},0}  \mbox{ for } 0\leqslant i\leqslant 3
\end{align*}
and
\begin{align*}
&\|\partial_{3}(\partial_{m}\eta_{n}\partial_{k}\eta_{l})\|_{\underline{j},0}\nonumber\\
&\lesssim \|\partial_{3}\partial_{m} \eta_{n}\partial_{k}\eta_{l}\|_{\underline{j},0}
+\|\partial_{m}\eta_{n}\partial_{3} \partial_{k}\eta_{l}\|_{\underline{j},0}\nonumber \\
&\lesssim
\|\nabla\eta\|_{\underline{2},0}^{1/2}\|\nabla\eta\|_{\underline{1},1}^{1/4}\|\nabla\eta\|_{\underline{2},1}^{1/4}
\|\nabla^2\eta\|_{\underline{j},0}
+\|\nabla^2\eta\|_{\underline{1},0}^{3/4}\|\nabla^2\eta\|_{\underline{2},0}^{1/4}
\|\nabla\eta\|_{\underline{j},0}^{1/2}\|\nabla\eta\|_{\underline{j},1}^{1/2}\nonumber\\
&\lesssim
\|\nabla\eta\|_{\underline{2},0}^{1/2}\|\nabla\eta\|_{\underline{1},1}^{1/4}\|\nabla\eta\|_{\underline{2},1}^{1/4}
\|\nabla\eta\|_{\underline{j},1} \mbox{ for } j=1,\ 2.
\end{align*}
Exploiting \eqref{202108261428nm} and the interpolation inequality, one obtains
\begin{align*}
&\|\partial_{3}^2(\partial_{m}\eta_{n}\partial_{k}\eta_{l})\|_{\underline{1},0}\nonumber\\
&\lesssim \|\partial_{3}^2\partial_{m}\eta_{n}\partial_{k}\eta_{l}\|_{\underline{1},0}
+\|\partial_{3}\partial_{m}\eta_{n}\partial_{3}\partial_{k}\eta_{l}\|_{\underline{1},0}
+\|\partial_{m}\eta_{n}\partial_{3}^2\partial_{k}\eta_{l}\|_{\underline{1},0}\nonumber \\
&\lesssim
\|\nabla\eta\|_{\underline{2},0}^{1/2}\|\nabla\eta\|_{\underline{1},1}^{1/4}\|\nabla\eta\|_{\underline{2},1}^{1/4}
\|\nabla^3\eta\|_{\underline{1},0}
+\|\nabla^2\eta\|_{\underline{1},0}^{3/4}\|\nabla^2\eta\|_{\underline{2},0}^{1/4}
\|\nabla^2\eta\|_{\underline{1},0}^{1/2}\|\nabla^2\eta\|_{\underline{1},1}^{1/2}\nonumber \\
&\lesssim\|\nabla\eta\|_{\underline{2},0}^{1/2}\|\nabla\eta\|_{\underline{1},1}^{1/4}\|\nabla\eta\|_{\underline{2},1}^{1/4}
\|\nabla\eta\|_{\underline{1},2}.
\end{align*}
Combining with the above three estimates then yields
\begin{align}
&\label{202108271400n}
\|\partial_{m}\eta_{n}\partial_{k}\eta_{l}\|_{\underline{i},0}
\lesssim\|\nabla\eta\|_{\underline{2},0}^{1/2}\|\nabla\eta\|_{\underline{1},1}^{1/4}\|\nabla\eta\|_{\underline{2},1}^{1/4}
\|\nabla\eta\|_{\underline{i},0}
 \mbox{ for }0\leqslant i\leqslant 3,\\[0.5mm]
&\label{202108271412n}
\|\partial_{m}\eta_{n}\partial_{k}\eta_{l}\|_{\underline{j},1}
\lesssim\|\nabla\eta\|_{\underline{2},0}^{1/2}\|\nabla\eta\|_{\underline{1},1}^{1/4}\|\nabla\eta\|_{\underline{2},1}^{1/4}
\|\nabla\eta\|_{\underline{j},1}
 \mbox{ for }j=1,\ 2,\\[0.5mm]
&\label{202108271400}
\|\partial_{m}\eta_{n}\partial_{k}\eta_{l}\|_{\underline{1},2}
\lesssim\|\nabla\eta\|_{\underline{2},0}^{1/2}\|\nabla\eta\|_{\underline{1},1}^{1/4}\|\nabla\eta\|_{\underline{2},1}^{1/4}
\|\nabla\eta\|_{\underline{1},2}.
\end{align}
Recalling the expressions of $\tilde{\mathcal{B}}$, $\tilde{\mathcal{B}}^{{L}}$ and $\tilde{\mathcal{B}}^{{N}}$,
we immediately deduce \eqref{202109071530}--\eqref{202108270845} from \eqref{202108271400n}--\eqref{202108271400} under the assumption \eqref{202108261502} with sufficiently small $\delta$.

Moreover,
utilizing \eqref{202108261428nm} and the fact that $\eta_{t}=u$, we arrive at
\begin{align*}
&
\|\tilde{\ml{B}_{t}}\|_{\underline{1},0}
\lesssim(1+\|\nabla\eta\|_{\underline{2},0}^{1/2}\|\nabla\eta\|_{\underline{1},1}^{1/4}\|\nabla\eta\|_{\underline{2},1}^{1/4})
\|\nabla u\|_{\underline{1},0},\\[1mm]
&
\|\tilde{\ml{B}_{t}}\|_{L^{\infty}}
\lesssim(1+\|\nabla \eta\|_{L^{\infty}})\|\nabla u\|_{L^{\infty}}
\lesssim(1+\|\nabla\eta\|_{\underline{2},0}^{1/2}\|\nabla\eta\|_{\underline{1},1}^{1/4}\|\nabla\eta\|_{\underline{2},1}^{1/4})
\|\nabla u\|_{L^{\infty}},
\end{align*}
which imply \eqref{202108270845nm} and \eqref{202108270845n12m} under the assumption \eqref{202108261502} with sufficiently small $\delta$.

(2)
Applying the determinant expansion theorem, it is easy to see that
\begin{align}\label{202108271645}
J=\mm{det}(\nabla\eta+I)
=1+\mm{div}\eta+r_{\eta},
\end{align}
where we have defined that $r_{\eta}:=r_2^{\eta}+r_3^{\eta}$,
$$\begin{aligned}
&r_2^{\eta}:= \partial_1\eta_1\partial_2\eta_2+\partial_1\eta_1\partial_3\eta_3+\partial_2\eta_2\partial_3\eta_3
-\partial_2\eta_1\partial_1\eta_2-\partial_2\eta_3\partial_3\eta_2-
\partial_3\eta_1\partial_1\eta_3
,\\[0.5mm] &r_3^{\eta}:=\partial_1\eta_1(\partial_2\eta_2\partial_3\eta_3-\partial_2\eta_3\partial_3\eta_2)
+\partial_2\eta_1(\partial_1\eta_3\partial_3\eta_2-\partial_1\eta_2\partial_3\eta_3)
+\partial_3\eta_1(\partial_1\eta_2\partial_2\eta_3-\partial_1\eta_3\partial_2\eta_2).
\end{aligned}$$
Making use of \eqref{202108271400n}--\eqref{202108271400}, one easily deduces that
\begin{align}
&\label{2021sadf08271725}
\|r_2^{\eta}\|_{\underline{i},0}
\lesssim \|\nabla\eta\|_{\underline{2},0}^{1/2}\|\nabla\eta\|_{\underline{1},1}^{1/4}\|\nabla\eta\|_{\underline{2},1}^{1/4}
\|\nabla\eta\|_{\underline{i},0} \mbox{ for }0\leqslant i\leqslant 3,\\[0.5mm]
&\label{202108271710}
\|r_2^{\eta}\|_{\underline{j},1}
\lesssim \|\nabla\eta\|_{\underline{2},0}^{1/2}\|\nabla\eta\|_{\underline{1},1}^{1/4}\|\nabla\eta\|_{\underline{2},1}^{1/4}
\|\nabla\eta\|_{\underline{j},1}\mbox{ for }j=1 ,\ 2,\\[0.5mm]
&\label{202108271712}
\|r_2^{\eta}\|_{\underline{1},2}
\lesssim \|\nabla\eta\|_{\underline{2},0}^{1/2}\|\nabla\eta\|_{\underline{1},1}^{1/4}\|\nabla\eta\|_{\underline{2},1}^{1/4}
\|\nabla\eta\|_{\underline{1},2}.
  \end{align}

Let $1\leqslant$ $k$, $l$, $m$, $n$, $q$, $r\leqslant3$.
Following the derivations of \eqref{202108271400n}--\eqref{202108271400},
we can also estimate that
\begin{align*}
\|\partial_{k}\eta_{l}\partial_{m}\eta_{n}\partial_{q}\eta_{r}\|_{\underline{i},0}
&\lesssim\|\partial_{k}\eta_{l}\|_{\underline{2},0}^{1/2}\|\partial_{k}\eta_{l}\|_{\underline{1},1}^{1/4}
\|\partial_{k}\eta_{l}\|_{\underline{2},1}^{1/4}\|\partial_{m}\eta_{n}\partial_{q}\eta_{r}\|_{\underline{i},0}\nonumber \\
&\quad\quad+\|\partial_{m}\eta_{n}\partial_{q}\eta_{r}\|_{\underline{2},0}^{1/2}\|\partial_{m}\eta_{n}\partial_{q}\eta_{r}\|_{\underline{1},1}^{1/4}
\|\partial_{m}\eta_{n}\partial_{q}\eta_{r}\|_{\underline{2},1}^{1/4}\|\partial_{k}\eta_{l}\|_{\underline{i},0},\\[1mm]
\|\partial_3(\partial_{k}\eta_{l}\partial_{m}\eta_{n}\partial_{q}\eta_{r})\|_{\underline{j},0}
&\lesssim\|\partial_{k}\eta_{l}\partial_{3}(\partial_{m}\eta_{n}\partial_{q}\eta_{r})\|_{\underline{j},0}
+\|\partial_{3}\partial_{k}\eta_{l}(\partial_{m}\eta_{n}\partial_{q}\eta_{r})\|_{\underline{j},0}
\nonumber\\
&\lesssim
\|\partial_{k}\eta_{l}\|_{\underline{2},0}^{1/2}\|\partial_{k}\eta_{l}\|_{\underline{1},1}^{1/4}\|\partial_{k}\eta_{l}\|_{\underline{2},1}^{1/4}
\|\partial_{3}(\partial_{m}\eta_{n}\partial_{q}\eta_{r})\|_{\underline{j},0}\nonumber\\
&\quad\quad+
\|\partial_{k}\eta_{l}\|_{\underline{j},0}^{1/2}\|\partial_{k}\eta_{l}\|_{\underline{j},1}^{1/2}
\|\partial_{3}(\partial_{m}\eta_{n}\partial_{q}\eta_{r})\|_{\underline{1},0}^{3/4}
\|\partial_{3}(\partial_{m}\eta_{n}\partial_{q}\eta_{r})\|_{\underline{2},0}^{1/4}
\nonumber\\
&\quad\quad
+\|(\partial_{m}\eta_{n}\partial_{q}\eta_{r})\|_{\underline{2},0}^{1/2}\|(\partial_{m}\eta_{n}\partial_{q}\eta_{r})\|_{\underline{1},1}^{1/4}
\|(\partial_{m}\eta_{n}\partial_{q}\eta_{r})\|_{\underline{2},1}^{1/4}
\|\partial_{3}\partial_{k}\eta_{l}\|_{\underline{j},0}
\nonumber\\
&\quad\quad
+\|\partial_{3}\partial_{k}\eta_{l}\|_{\underline{1},0}^{3/4}\|\partial_{3}\partial_{k}\eta_{l}\|_{\underline{2},0}^{1/4}
\|(\partial_{m}\eta_{n}\partial_{q}\eta_{r})\|_{\underline{j},0}^{1/2}
\|(\partial_{m}\eta_{n}\partial_{q}\eta_{r})\|_{\underline{j},1}^{1/2} \nonumber
\end{align*}
and
\begin{align*}
&\|\partial_{3}^2(\partial_{k}\eta_{l}\partial_{m}\eta_{n}\partial_{q}\eta_{r})\|_{\underline{1},0}\nonumber\\
&\lesssim \|\partial_{3}^2\partial_{k}\eta_{l}(\partial_{m}\eta_{n}\partial_{q}\eta_{r})\|_{\underline{1},0}
+\|\partial_{3}\partial_{k}\eta_{l}\partial_{3}(\partial_{m}\eta_{n}\partial_{q}\eta_{r})\|_{\underline{1},0}
+\|\partial_{k}\eta_{l}\partial_{3}^2(\partial_{m}\eta_{n}\partial_{q}\eta_{r})\|_{\underline{1}0}\nonumber \\ &\lesssim\|\partial_{m}\eta_{n}\partial_{q}\eta_{r}\|_{\underline{2},0}^{1/2}\|\partial_{m}\eta_{n}\partial_{q}\eta_{r}\|_{\underline{1},1}^{1/4}
\|\partial_{m}\eta_{n}\partial_{q}\eta_{r}\|_{\underline{2},1}^{1/4}\|\partial_{3}^2\partial_{k}\eta_{l}\|_{\underline{1},0}\nonumber \\
&\quad\quad
+{\|\partial_{3}\partial_{k}\eta_{l}\|_{\underline{1},0}^{3/4}\|\partial_{3}\partial_{k}\eta_{l}\|_{\underline{2},0}^{1/4}}
\|\partial_{3}(\partial_{m}\eta_{n}\partial_{q}\eta_{r})\|_{\underline{1},0}^{1/2}
\|\partial_{3}(\partial_{m}\eta_{n}\partial_{q}\eta_{r})\|_{\underline{1},1}^{1/2}\nonumber \\
&\quad\quad
+{\|\partial_{3}(\partial_{m}\eta_{n}\partial_{q}\eta_{r})\|_{\underline{1},0}^{3/4}
\|\partial_{3}(\partial_{m}\eta_{n}\partial_{q}\eta_{r})\|_{\underline{2},0}^{1/4}}
\|\partial_{3}\partial_{k}\eta_{l}\|_{\underline{1},0}^{1/2}
\|\partial_{3}\partial_{k}\eta_{l}\|_{\underline{1},1}^{1/2}\nonumber \\
&\quad\quad
+\|\partial_{k}\eta_{l}\|_{\underline{2},0}^{1/2}\|\partial_{k}\eta_{l}\|_{\underline{1},1}^{1/4}
\|\partial_{k}\eta_{l}\|_{\underline{2},1}^{1/4}\|\partial_{3}^2(\partial_{m}\eta_{n}\partial_{q}\eta_{r})\|_{\underline{1},0},
\end{align*}
where  $0\leqslant i\leqslant 3$ and $j=1$, $2$.
Exploiting \eqref{202108261502}, \eqref{202108271400n}--\eqref{202108271400}, the interpolation inequality and the above three estimates, we obtain that
\begin{align}
&\label{2020109301202}
\|r_3^{\eta}\|_{\underline{i},0}
\lesssim \|\nabla\eta\|_{\underline{2},0}^{1/2}\|\nabla\eta\|_{\underline{1},1}^{1/4}\|\nabla\eta\|_{\underline{2},1}^{1/4}
\|\nabla\eta\|_{\underline{i},0} \mbox{ for }0\leqslant i\leqslant 3,\\[0.5mm]
&\label{2020109301205}
\|r_3^{\eta}\|_{\underline{j},1}
\lesssim \|\nabla\eta\|_{\underline{2},0}^{1/2}\|\nabla\eta\|_{\underline{1},1}^{1/4}\|\nabla\eta\|_{\underline{2},1}^{1/4}
\|\nabla\eta\|_{\underline{j},1} \mbox{ for }j=1,\ 2,\\[0.5mm]
&\label{20201093012saf02}
\|r_3^{\eta}\|_{\underline{1},2}
\lesssim \|\nabla\eta\|_{\underline{2},0}^{1/2}\|\nabla\eta\|_{\underline{1},1}^{1/4}\|\nabla\eta\|_{\underline{2},1}^{1/4}
\|\nabla\eta\|_{\underline{1},2}.
\end{align}
Combining \eqref{2021sadf08271725}--\eqref{202108271712}
with \eqref{2020109301202}--\eqref{20201093012saf02} then yields
\begin{align}
&\label{202109211000}
\|r_{\eta}\|_{\underline{i},0}\lesssim\|r_2^{\eta}\|_{\underline{i},0}+\|r_3^{\eta}\|_{\underline{i},0}
\lesssim \|\nabla\eta\|_{\underline{2},0}^{1/2}\|\nabla\eta\|_{\underline{1},1}^{1/4}\|\nabla\eta\|_{\underline{2},1}^{1/4}
\|\nabla\eta\|_{\underline{i},0} \mbox{ for }0\leqslant i\leqslant 3,\\[0.5mm]
&\label{202109211010}
\|r_{\eta}\|_{\underline{j},1}\lesssim\|r_2^{\eta}\|_{\underline{j},1}+\|r_3^{\eta}\|_{\underline{j},1}
\lesssim \|\nabla\eta\|_{\underline{2},0}^{1/2}\|\nabla\eta\|_{\underline{1},1}^{1/4}\|\nabla\eta\|_{\underline{2},1}^{1/4}
\|\nabla\eta\|_{\underline{j},1}\mbox{ for }j=1,\ 2,\\[0.5mm]
&\label{202108280850}
\|r_{\eta}\|_{\underline{1},2}\lesssim\|r_2^{\eta}\|_{\underline{1},2}+\|r_3^{\eta}\|_{\underline{1},2}
\lesssim
\|\nabla\eta\|_{\underline{2},0}^{1/2}\|\nabla\eta\|_{\underline{1},1}^{1/4}\|\nabla\eta\|_{\underline{2},1}^{1/4}
\|\nabla\eta\|_{\underline{1},2}.
 \end{align}
Hence,
in the light of \eqref{202108261502} and \eqref{202109211000}--\eqref{202108280850}, we arrive at
\begin{align}
&\label{202109071350}
\| \mm{div}\eta+r_{\eta} \|_{\underline{i},0}\lesssim \|\nabla\eta\|_{\underline{i},0} \mbox{ for } 0\leqslant i\leqslant 3,\\[0.5mm]
&\label{202109071355}
\| \mm{div}\eta+r_{\eta} \|_{\underline{j},1}\lesssim \|\nabla\eta\|_{\underline{j},1} \mbox{ for }j=1 ,\ 2,\\[0.5mm]
&\label{202108280925}
\| \mm{div}\eta+r_{\eta} \|_{\underline{1},2}\lesssim \|\nabla\eta\|_{\underline{1},2}.
\end{align}
In addition, similarly to the derivations of \eqref{202108270845nm} and \eqref{202108270845n12m}, we can infer that
\begin{align}
&
\|\partial_{t}r_{\eta}\|_{\underline{1},0}\lesssim
\|\nabla\eta\|_{\underline{2},0}^{1/2}\|\nabla\eta\|_{\underline{1},1}^{1/4}\|\nabla\eta\|_{\underline{2},1}^{1/4}
\|\nabla u\|_{\underline{1},0},\nonumber\\[0.5mm]
&
\|\partial_{t}r_{\eta}\|_{L^{\infty}}\lesssim
\|\nabla\eta\|_{\underline{2},0}^{1/2}\|\nabla\eta\|_{\underline{1},1}^{1/4}\|\nabla\eta\|_{\underline{2},1}^{1/4}
\|\nabla u\|_{L^{\infty}},\nonumber\\[0.5mm]
&\label{202408291512n}
\|\partial_t(\mm{div}\eta+r_{\eta})\|_{\underline{1},0}\lesssim\|\nabla u\|_{\underline{1},0},\\[0.5mm]
&\label{202408291515n}
\|\partial_t(\mm{div}\eta+r_{\eta})\|_{L^{\infty}}\lesssim\|\nabla u\|_{L^{\infty}}.
\end{align}

Next we are devoted to deriving \eqref{202109071535}--\eqref{202108270855n12nm} in sequence.
Using \eqref{202108261502},  \eqref{202109071350}--\eqref{202108280925} and \eqref{202108261424}, we obtain from \eqref{202108271645} that
\begin{align}
\label{202201093001606}
\|J-1\|_{L^\infty}\lesssim \| \mm{div}\eta+r_{\eta} \|_{L^\infty}\lesssim   \delta,
\end{align}
which yields \eqref{202109071535} for sufficiently small $\delta$.

In view of \eqref{202108271645} and the fact that $\eta_{t}=u$, we have
\begin{align}
&\label{202108280902nm}
J^{-1}-1
=(1-J^{-1})(\mm{div}\eta+r_{\eta})-(\mm{div}\eta+r_{\eta}),\\[1.5mm]
&\label{202108280902}
J^{-1}-1+\mm{div}\eta=(1-J^{-1})(\mm{div}\eta+r_{\eta})-r_{\eta},\\[1.5mm]
&\label{202108280912}
J_{t}^{-1}=(J^{-1}-1)_{t}
=(1-J^{-1})_{t}(\mm{div}\eta+r_{\eta})+(1-J^{-1})(\mm{div}\eta+r_{\eta})_{t}-(\mm{div}\eta+r_{\eta})_{t},\\[1.5mm]
&\label{202108280915}
(J^{-1}-1+\mm{div}\eta)_{t}
=(1-J^{-1})_{t}(\mm{div}\eta+r_{\eta})+(1-J^{-1})(\mm{div}\eta+r_{\eta})_{t}-\partial_{t}r_{\eta}.
\end{align}
Making use of \eqref{202109071350}--\eqref{202109071355} and \eqref{202108261428nm}, we infer that
\begin{align*}
&\|(1-J^{-1})(\mm{div}\eta+r_{\eta})\|_{\underline{s},0}\\
&\lesssim
\|(\mm{div}\eta+r_{\eta})\|_{\underline{2},0}^{1/2}
\|(\mm{div}\eta+r_{\eta})\|_{\underline{1},1}^{1/4}\|(\mm{div}\eta+r_{\eta})\|_{\underline{2},1}^{1/4}\|(1-J^{-1})\|_{\underline{s},0}\\[0.5mm]
&\lesssim\|\nabla\eta\|_{\underline{2},0}^{1/2}\|\nabla\eta\|_{\underline{1},1}^{1/4}\|\nabla\eta\|_{\underline{2},1}^{1/4}
\|(1-J^{-1})\|_{\underline{s},0} \mbox{ for } s=0,\ 1,
\end{align*}
which, together with \eqref{202108261502}, \eqref{202109071350} and \eqref{202108280902nm}, yields
\begin{align}
&\label{202408301024}
\|(1-J^{-1})\|_{\underline{s},0}\lesssim\|\nabla\eta\|_{\underline{s},0}\; \mbox{ for }s=0,\ 1.
\end{align}

In the same manner, we can utilize \eqref{202109071350}--\eqref{202109071355} and \eqref{202108261428nm} to derive that, for $1\leqslant n\leqslant 3$,
\begin{align*}
&\|\partial_{n}\big((1-J^{-1})(\mm{div}\eta+r_{\eta})\big)\|_{\underline{1},0}
\nonumber\\
&\lesssim\|\partial_{n}(1-J^{-1})(\mm{div}\eta+r_{\eta})\|_{\underline{1},0}
+\|(1-J^{-1})\partial_{n}(\mm{div}\eta+r_{\eta})\|_{\underline{1},0}
\nonumber\\
&\lesssim\|\partial_{n}(1-J^{-1})(\mm{div}\eta+r_{\eta})\|_{\underline{1},0}
+\|(\mm{div}\eta+r_{\eta})\partial_{n}(\mm{div}\eta+r_{\eta})\|_{\underline{1},0}\nonumber\\[0.5mm]
&\quad\quad
+\|(1-J^{-1})(\mm{div}\eta+r_{\eta})\partial_{n}(\mm{div}\eta+r_{\eta})\|_{\underline{1},0}
\nonumber\\
&\lesssim
\|(\mm{div}\eta+r_{\eta})\|_{\underline{2},0}^{1/2}
\|(\mm{div}\eta+r_{\eta})\|_{\underline{1},1}^{1/4}\|(\mm{div}\eta+r_{\eta})\|_{\underline{2},1}^{1/4}
\|\partial_{n}(1-J^{-1})\|_{\underline{1},0}\nonumber\\[0.5mm]
&\quad\;\;+
\|(\mm{div}\eta+r_{\eta})\|_{\underline{2},0}^{1/2}
\|(\mm{div}\eta+r_{\eta})\|_{\underline{1},1}^{1/4}\|(\mm{div}\eta+r_{\eta})\|_{\underline{2},1}^{1/4}
\|\partial_{n}(\mm{div}\eta+r_{\eta})\|_{\underline{1},0}\nonumber\\[0.5mm]
&\quad\;\;+\|(\mm{div}\eta+r_{\eta})\|_{\underline{2},0}^{1/2}
\|(\mm{div}\eta+r_{\eta})\|_{\underline{1},1}^{1/4}\|(\mm{div}\eta+r_{\eta})\|_{\underline{2},1}^{1/4}\nonumber\\[0.5mm]
&\quad\quad\times
{\|\partial_{n}(\mm{div}\eta+r_{\eta})\|_{\underline{1},0}^{1/2}\|\partial_{n}(\mm{div}\eta+r_{\eta})\|_{\underline{2},0}^{1/2}}
\|(1-J^{-1})\|_{\underline{1},0}^{1/2}\|(1-J^{-1})\|_{\underline{1},1}^{1/2}
\nonumber\\[0.5mm]
&\lesssim\|\nabla\eta\|_{\underline{2},0}^{1/2}\|\nabla\eta\|_{\underline{1},1}^{1/4}\|\nabla\eta\|_{\underline{2},1}^{1/4}
(\|(1-J^{-1})\|_{\underline{1},1}+\|\nabla\eta\|_{\underline{1},1})\nonumber\\[0.5mm]
&\quad\;\;
+\|\nabla\eta\|_{\underline{2},0}^{1/2}\|\nabla\eta\|_{\underline{1},1}^{1/4}\|\nabla\eta\|_{\underline{2},1}^{1/4}
{\|(1-J^{-1})\|_{\underline{1},0}^{1/2}\|\nabla\eta\|_{\underline{2},1}^{1/2}}
\|\nabla\eta\|_{\underline{1},1}^{1/2}\|(1-J^{-1})\|_{\underline{1},1}^{1/2},
\end{align*}
which, in combination with \eqref{202108261502}, \eqref{202109071350}--\eqref{202109071355}, \eqref{202408301024} and \eqref{202108280902nm}, gives rise to
\begin{align}
&\label{202408301027}
\|(1-J^{-1})\|_{\underline{1},1}\lesssim\|\nabla\eta\|_{\underline{1},1}.
\end{align}
In view of \eqref{202109071350}--\eqref{202109071355}, \eqref{202408301024}--\eqref{202408301027} and \eqref{202108261428nm} with $i=2$, we further deduce that
\begin{align*}
&\|(1-J^{-1})(\mm{div}\eta+r_{\eta})\|_{\underline{2},0}\nonumber\\[0.5mm]
&\lesssim \|(\mm{div}\eta+r_{\eta})\|_{\underline{2},0}^{1/2}\|(\mm{div}\eta+r_{\eta})\|_{\underline{1},1}^{1/4}
\|(\mm{div}\eta+r_{\eta})\|_{\underline{2},1}^{1/4}\|(1-J^{-1})\|_{\underline{2},0}\nonumber\\[0.5mm]
&\quad+\|(\mm{div}\eta+r_{\eta})\|_{\underline{2},0}^{1/2}\|(\mm{div}\eta+r_{\eta})\|_{\underline{3},0}^{1/2}
\|(1-J^{-1})\|_{\underline{1},0}^{1/2}\|(1-J^{-1})\|_{1}^{1/4}\|(1-J^{-1})\|_{\underline{1},1}^{1/4}
\nonumber\\[0.5mm]
&\lesssim \|\nabla\eta\|_{\underline{2},0}^{1/2}
\|\nabla\eta\|_{\underline{1},1}^{1/4}\|\nabla\eta\|_{\underline{2},1}^{1/4}\|(1-J^{-1})\|_{\underline{2},0}
+\|\nabla\eta\|_{\underline{2},0}^{1/2}\|\nabla\eta\|_{\underline{3},0}^{1/2}
\|\nabla\eta\|_{\underline{1},0}^{1/2}\|\nabla\eta\|_{1}^{1/4}\|\nabla\eta\|_{\underline{1},1}^{1/4},
\nonumber
\end{align*}
which, along with \eqref{202108261502}, \eqref{202109071355} and \eqref{202108280902nm}, implies that
\begin{align}
&\label{202408301045}
\|(1-J^{-1})\|_{\underline{2},0}\lesssim\|\nabla\eta\|_{\underline{2},0}.
\end{align}
Analogously, exploiting  \eqref{202109071350}--\eqref{202109071355}, \eqref{202408301024}--\eqref{202408301045} and \eqref{202108261428nm},
we can estimate that
\begin{align*}
&\|\partial_{n}\big((1-J^{-1})(\mm{div}\eta+r_{\eta})\big)\|_{\underline{2},0}
\nonumber\\
&\lesssim\|\partial_{n}(1-J^{-1})(\mm{div}\eta+r_{\eta})\|_{\underline{2},0}
+\|(1-J^{-1})\partial_{n}(\mm{div}\eta+r_{\eta})\|_{\underline{2},0}
\nonumber\\
&\lesssim
\|(\mm{div}\eta+r_{\eta})\|_{\underline{2},0}^{1/2}
\|(\mm{div}\eta+r_{\eta})\|_{\underline{1},1}^{1/4}\|(\mm{div}\eta+r_{\eta})\|_{\underline{2},1}^{1/4}
\|\partial_{n}(1-J^{-1})\|_{\underline{2},0}\nonumber\\[0.5mm]
&\quad\quad
+\|(\mm{div}\eta+r_{\eta})\|_{\underline{2},0}^{1/2}
\|(\mm{div}\eta+r_{\eta})\|_{\underline{2},1}^{1/2}
\|\partial_{n}(1-J^{-1})\|_{\underline{1},0}^{3/4}\|\partial_{n}(1-J^{-1})\|_{\underline{2},0}^{1/4}\nonumber\\[0.5mm]
&\quad\quad+
\|(1-J^{-1})\|_{\underline{2},0}^{1/2}
\|(1-J^{-1})\|_{\underline{1},1}^{1/4}\|(1-J^{-1})\|_{\underline{2},1}^{1/4}
\|\partial_{n}(\mm{div}\eta+r_{\eta})\|_{\underline{2},0}\nonumber\\[0.5mm]
&\quad\quad
+\|(1-J^{-1})\|_{\underline{2},0}^{1/2}
\|(1-J^{-1})\|_{\underline{2},1}^{1/2}
\|\partial_{n}(\mm{div}\eta+r_{\eta})\|_{\underline{1},0}^{3/4}\|\partial_{n}(\mm{div}\eta+r_{\eta})\|_{\underline{2},0}^{1/4}
\nonumber\\[0.5mm]
&\lesssim\|\nabla\eta\|_{\underline{2},0}^{1/2}\|\nabla\eta\|_{\underline{1},1}^{1/4}\|\nabla\eta\|_{\underline{2},1}^{1/4}
(\|(1-J^{-1})\|_{\underline{2},1}+\|\nabla\eta\|_{\underline{2},1}),
\end{align*}
which, along with \eqref{202108261502}, \eqref{202109071350}--\eqref{202109071355} and \eqref{202108280902nm},  yields
\begin{align}
\label{202201009292319}
\|(1-J^{-1})\|_{\underline{2},1}\lesssim\|\nabla\eta\|_{\underline{2},1}.
\end{align}
Consequently,  the estimate \eqref{202108270852} follows by combining \eqref{202408301027} with \eqref{202201009292319}.

Making use of \eqref{202108270852}, \eqref{202109071350}--\eqref{202109071355}, \eqref{202408301027}--\eqref{202408301045}
 and \eqref{202108261428nm}, we can have
\begin{align}
&\|(1-J^{-1})(\mm{div}\eta+r_{\eta})\|_{\underline{3},0}\nonumber\\
&\lesssim \|(1-J^{-1})\|_{\underline{3},0}\|(\mm{div}\eta+r_{\eta})\|_{\underline{2},0}^{1/2}
\|(\mm{div}\eta+r_{\eta})\|_{\underline{1},1}^{1/4}\|(\mm{div}\eta+r_{\eta})\|_{\underline{2},1}^{1/4}\nonumber\\
&\quad+\|(1-J^{-1})\|_{\underline{2},0}^{1/2}\|(1-J^{-1})\|_{\underline{1},1}^{1/4}\|(1-J^{-1})\|_{\underline{2},1}^{1/4}
\|(\mm{div}\eta+r_{\eta})\|_{\underline{3},0}
\nonumber\\
&\lesssim \|\nabla\eta\|_{\underline{2},0}^{1/2}\|\nabla\eta\|_{\underline{1},1}^{1/4}\|\nabla\eta\|_{\underline{2},1}^{1/4}
(\|(1-J^{-1})\|_{\underline{3},0}+\|\nabla\eta\|_{\underline{3},0}),
\nonumber
\end{align}
which, together with \eqref{202109071350}, \eqref{202408301024} and \eqref{202408301045}, implies \eqref{202109071532}
for sufficiently small $\delta$.

Additionally, exploiting \eqref{202108261428nm}, \eqref{202109071532}--\eqref{202108270852}, \eqref{202109071350}--\eqref{202109071355}, \eqref{202108280902nm} and the interpolation inequality, one infers that
\begin{align}\label{202408241725}
&\|\partial_{3}^2\big((1-J^{-1})(\mm{div}\eta+r_{\eta})\big)\|_{\underline{1},0}
\nonumber\\[0.5mm]
&\lesssim\|\partial_{3}^2(1-J^{-1})(\mm{div}\eta+r_{\eta})\|_{\underline{1},0}
+\|\partial_{3}(1-J^{-1})\partial_{3}(\mm{div}\eta+r_{\eta})\|_{\underline{1},0}
\nonumber \\
&\quad +\|(1-J^{-1})\partial_{3}^2(\mm{div}\eta+r_{\eta})\|_{\underline{1},0}\nonumber\\[0.5mm]
&\lesssim \|(\mm{div}\eta+r_{\eta})\|_{\underline{2},0}^{1/2}
\|(\mm{div}\eta+r_{\eta})\|_{\underline{1},1}^{1/4}\|(\mm{div}\eta+r_{\eta})
\|_{\underline{2},1}^{1/4}\|(1-J^{-1})\|_{\underline{1},2}
\nonumber\\[0.5mm]
&\quad +\|\partial_{3}(1-J^{-1})\|_{\underline{1},0}^{3/4}
\|\partial_{3}(1-J^{-1})\|_{\underline{2},0}^{1/4}\|\partial_{3}(\mm{div}\eta+r_{\eta})\|_{\underline{1},0}^{1/2}
\|\partial_{3}(\mm{div}\eta+r_{\eta})\|_{\underline{1},1}^{1/2}\nonumber\\[0.5mm]
&\quad +\|\partial_{3}(1-J^{-1})\|_{\underline{1},0}^{1/2}
\|\partial_{3}(1-J^{-1})\|_{\underline{1},1}^{1/2}\|\partial_{3}(\mm{div}\eta+r_{\eta})\|_{\underline{1},0}^{3/4}
\|\partial_{3}(\mm{div}\eta+r_{\eta})\|_{\underline{2},0}^{1/4}\nonumber\\[0.5mm]
&\quad +\|(1-J^{-1})\|_{\underline{2},0}^{1/2}
\|(1-J^{-1})\|_{\underline{1},1}^{1/4}\|(1-J^{-1})\|_{\underline{2},1}^{1/4}\|(\mm{div}\eta+r_{\eta})\|_{\underline{1},2}\nonumber\\[0.5mm]
&\lesssim \|\nabla\eta\|_{\underline{2},0}^{1/2}\|\nabla\eta\|_{\underline{1},1}^{1/4}\|\nabla\eta\|_{\underline{2},1}^{1/4}(\|\nabla\eta\|_{\underline{1},2}+\|(1-J^{-1})\|_{\underline{1},2}).
\end{align}
Combining \eqref{202108280925} with \eqref{202408241725} and \eqref{202108270852} then yields \eqref{202108270852nm} for sufficiently small $\delta$.

{Furthermore, utilizing \eqref{202109071532}--\eqref{202108270852nm}, \eqref{202109071350}--\eqref{202109071355}, \eqref{202408291512n}--\eqref{202408291515n} and \eqref{202408241430}--\eqref{202108261428nm}, we can deduce from \eqref{202108280912} that
\begin{align*}
&
\|J^{-1}_{t}\|_{\underline{1},0}\lesssim
\|\nabla\eta\|_{\underline{2},0}^{1/2}\|\nabla\eta\|_{\underline{1},1}^{1/4}
\|\nabla\eta\|_{\underline{2},1}^{1/4}(\|J^{-1}_{t}\|_{\underline{1},0}+\|\nabla u\|_{\underline{1},0})
+\|\nabla u\|_{\underline{1},0},\\[1mm]
&
\|J^{-1}_{t}\|_{L^{\infty}}\lesssim
\|\nabla\eta\|_{\underline{2},0}^{1/2}\|\nabla\eta\|_{\underline{1},1}^{1/4}
\|\nabla\eta\|_{\underline{2},1}^{1/4}(\|J^{-1}_{t}\|_{L^{\infty}}+\|\nabla u\|_{L^{\infty}})+\|\nabla u\|_{L^{\infty}}.
\end{align*}
Hence, we obtain \eqref{202108270852n12} and \eqref{202108270852n12nm} directly from the above two estimates
under the assumption \eqref{202108261502} with sufficiently small $\delta$.
}

Thanks to \eqref{202109071532}--\eqref{202108270852nm}, \eqref{202109071350}--\eqref{202108280925} and \eqref{202109211000}--\eqref{202108280850},
we can follow the derivation of \eqref{202108271400n}--\eqref{202108271400}  to deduce from \eqref{202108280902} that
\begin{align}
&\label{2021082sdfa70856}
 \|J^{-1}-1+\mm{div}\eta\|_{\underline{i},0}
\lesssim \|\nabla\eta\|_{\underline{2},0}^{1/2}\|\nabla\eta\|_{\underline{1},1}^{1/4}
\|\nabla\eta\|_{\underline{2},1}^{1/4}\|\nabla\eta\|_{\underline{i},0}
 \mbox{ for } 0\leqslant i\leqslant 3,\\[0.5mm]
&\label{2021082sadfa70854}
\|J^{-1}-1+\mm{div}\eta\|_{\underline{j},1}
\lesssim\|\nabla\eta\|_{\underline{2},0}^{1/2}\|\nabla\eta\|_{\underline{1},1}^{1/4}
\|\nabla\eta\|_{\underline{2},1}^{1/4}\|\nabla\eta\|_{\underline{3},0}  \mbox{ for }  j=1,\ 2,\\[0.5mm]
&\label{2021082sadfa70855}
\|J^{-1}-1+\mm{div}\eta\|_{\underline{2},1}
\lesssim\|\nabla\eta\|_{\underline{2},0}^{1/2}\|\nabla\eta\|_{\underline{1},1}^{1/4}
\|\nabla\eta\|_{\underline{2},1}^{1/4}\|\nabla\eta\|_{\underline{2},1}.
\end{align}
Similarly, we can  follow the derivations of \eqref{202108270852n12}--\eqref{202108270852n12nm} to deduce from \eqref{202108280915} that
\begin{align}
&\label{202408291522n}
\|(J^{-1}-1+\mm{div}\eta)_{t}\|_{\underline{1},0}\lesssim
\|\nabla\eta\|_{\underline{2},0}^{1/2}\|\nabla\eta\|_{\underline{1},1}^{1/4}
\|\nabla\eta\|_{\underline{2},1}^{1/4}\|\nabla u\|_{\underline{1},0},\\[1mm]
&\label{202408291525n}
\|(J^{-1}-1+\mm{div}\eta)_{t}\|_{L^{\infty}}\lesssim
\|\nabla\eta\|_{\underline{2},0}^{1/2}\|\nabla\eta\|_{\underline{1},1}^{1/4}
\|\nabla\eta\|_{\underline{2},1}^{1/4}\|\nabla u\|_{L^{\infty}}.
\end{align}

Next we turn to the estimates of $\mathcal{R}^\eta$.
From now on, we define that
$$P'':=\frac{\mm{d}^2}{\mm{d}z^2} P,\
\mathcal{P}:=\bar{\rho}\int_{0}^{\bar{\rho}(J^{-1}-1)}   P ''(\bar{\rho}+z)\mm{d}z$$
and
$$
\mathbb{S}(y,t):=
\begin{cases}
(0,\bar{\rho}(J^{-1}-1))&\mbox{for }J^{-1}-1>0;\\
(\bar{\rho}(J^{-1}-1),0)&\mbox{for }J^{-1}-1<0;\\
\{0\}&\mbox{for }J^{-1}-1=0.
\end{cases}$$
By \eqref{202201093001606}  and the condition $P(\cdot)\in C^4(\mathbb{R}^+)$, we easily see that, for sufficiently small
$\delta$,
\begin{align}
\nonumber
\bar{\rho}/2\leqslant \bar{\rho}J^{-1}\leqslant 2 \bar{\rho}
\end{align}
and
\begin{align}
\label{202202044}\sup_{(y,t)\in \Omega\times I_T} \left|\left. P'' (\bar{\rho}+z)
\right|_{z\in \mathbb{S}(y,t)}\right|
+\left\|\left.\frac{\mm{d}^{j+1}}{\mm{d}z^{j+1}} P (z)
\right|_{z=\bar{\rho}J^{-1}}\right\|_{L^\infty}\lesssim 1 ,
\end{align}
where $1\leqslant j\leqslant 3$.
Hence we  first have
\begin{align}
&\left\|\mathcal{P}\right\|_{0}\lesssim\left\|\left. P'' (z)
\right|_{z=\bar{\rho}J^{-1}}\right\|_{L^\infty}\|(J^{-1}-1)\|_{0} \lesssim
\|\nabla \eta\|_{0 }. \label{202010102105}
\end{align}

Moreover, noting that
$$\quad\partial_{m}\mathcal{P}
=\bar{\rho} \partial_{m} (J^{-1}-1) \left. P ''(z)
\right|_{z=\bar{\rho}J^{-1}} \mbox{ for }m=1,\ 2,\ 3,\ t,
$$
In light of \eqref{202109071532}--\eqref{202108270852nm} and \eqref{202202044}--\eqref{202010102105}, we can further derive that
\begin{align}
&\label{202109071550}
\|\mathcal{P} \|_{\underline{i},0}\lesssim \|\nabla\eta\|_{\underline{i},0} \mbox{ for } 0\leqslant i
\leqslant 3,\\[0.5mm]
&\label{202109071555}
\|\mathcal{P}\|_{\underline{j},1}\lesssim \|\nabla\eta\|_{\underline{j},1} \mbox{ for }j=1 ,\ 2,\\[0.5mm]
&\label{202108281225}
\|\mathcal{P}\|_{\underline{1},2}\lesssim \|\nabla\eta\|_{\underline{1},2}.
\end{align}
{
In the same manner,
\begin{align}
&\label{202408291633}
\left\|\mathcal{P}_{t}\right\|_{\underline{1},0}\lesssim \|\left. \partial_{t}(J^{-1}-1)P'' (z)
\right|_{z=\bar{\rho}J^{-1}} \|_{\underline{1},0}\lesssim \|J^{-1}_{t}\|_{\underline{1},0 },\\[0.5mm]
&\label{202408291633nm}
\left\|\mathcal{P}_{t}\right\|_{L^{\infty}}\lesssim \|\left. \partial_{t}(J^{-1}-1)P'' (z)
\right|_{z=\bar{\rho}J^{-1}} \|_{L^{\infty}}\lesssim \|J^{-1}_{t}\|_{L^{\infty}}.
\end{align}
}

Recalling the definition of $\mathcal{R}^{\eta}$, it holds that
\begin{align}\label{202408291702}
\partial_{m}\mathcal{R}^{\eta}
&=\partial_{m}(J^{-1}-1)\mathcal{P}
+(J^{-1}-1)\partial_{m}\mathcal{P}-\partial_{m}\int_{0}^{\bar{\rho}(J^{-1}-1)}z\frac{\mm{d}^2}{\mm{d}z^2} P (\bar{\rho}+z)\mm{d}z
\nonumber\\[1mm]
&=\partial_{m}(J^{-1}-1)\mathcal{P}
+(J^{-1}-1)\partial_{m}\mathcal{P}-\bar{\rho}^2(J^{-1}-1)P'' (\bar{\rho}J^{-1})\partial_{m}(J^{-1}-1)
\end{align}
for any $m=1$, $2$, $3$ and $t$.
Moreover, by \eqref{202202044} and \eqref{202109071532}--\eqref{202108270852nm}, one easily has
\begin{align}
\|\mathcal{R}^{\eta}\|_0
&\lesssim \|(J^{-1}-1)\|_{L^{\infty}}\|(J^{-1}-1)\|_{0}
\lesssim \|\nabla \eta\|_{\underline{2},0  }^{1/2}
\|\nabla \eta\|_{\underline{1}, 1}^{1/4}\|\nabla \eta\|_{\underline{2},1}^{1/4}\|\nabla \eta\|_{0}.
\label{da20202109301506}
\end{align}
Thanks to \eqref{202109071550}--\eqref{202108281225} and \eqref{202408291702}--\eqref{da20202109301506}, we easily derive that
\begin{align}
&\label{202108270asfd856}
 \|\mathcal{R}^{\eta}\|_{\underline{i},0}
 \lesssim \|\nabla\eta\|_{\underline{2},0}^{1/2}\|\nabla\eta\|_{\underline{1},1}^{1/4}\|\nabla\eta\|_{\underline{2},1}^{1/4}
\|\nabla\eta\|_{\underline{i},0}  \mbox{ for }  0\leqslant i\leqslant 3,\\[1mm]
&\label{202108sadf270854}
 \|\mathcal{R}^{\eta}\|_{\underline{j},1}
   \lesssim  \|\nabla\eta\|_{\underline{2},0}^{1/2}\|\nabla\eta\|_{\underline{1},1}^{1/4}\|\nabla\eta\|_{\underline{2},1}^{1/4}
   \|\nabla\eta\|_{\underline{j},1} \mbox{ for }  j=1,\ 2,\\[1mm]
&\label{202108sadf270855}
 \|\mathcal{R}^{\eta}\|_{\underline{1},2}
\lesssim
   \|\nabla\eta\|_{\underline{2},0}^{1/2}\|\nabla\eta\|_{\underline{1},1}^{1/4}\|\nabla\eta\|_{\underline{1},2}^{1/4}
   \|\nabla\eta\|_{\underline{2},1}.
\end{align}
Hence, we immediately obtain \eqref{202108270856}--\eqref{202108270854} by collecting \eqref{2021082sdfa70856}--\eqref{2021082sadfa70855} and \eqref{202108270asfd856}--\eqref{202108sadf270855} together.

Making use of \eqref{202108270852n12}--\eqref{202108270852n12nm}, \eqref{202202044}, \eqref{202408291633}--\eqref{202408291702} and \eqref{202108261428nm}, we easily obtain that
\begin{align}
&\label{202408270855}
 \|\mathcal{R}^{\eta}_{t}\|_{\underline{1},0}
   \lesssim  \|\nabla\eta\|_{\underline{2},0}^{1/2}\|\nabla\eta\|_{\underline{1},1}^{1/4}\|\nabla\eta\|_{\underline{2},1}^{1/4}
   \|\nabla u\|_{\underline{1},0},\\[1mm]
&\label{202408270855nm}
 \|\mathcal{R}^{\eta}_{t}\|_{L^{\infty}}
\lesssim
   \|\nabla\eta\|_{\underline{2},0}^{1/2}\|\nabla\eta\|_{\underline{1},1}^{1/4}\|\nabla\eta\|_{\underline{2},1}^{1/4}
   \|\nabla u\|_{L^{\infty}}.
\end{align}
Consequently, by collecting \eqref{202408291522n}--\eqref{202408291525n} and \eqref{202408270855}--\eqref{202408270855nm} together, the desired estimates \eqref{202108270855n12}--\eqref{202108270855n12nm} follow.
This completes the proof of Lemma \ref{lem:2108261504}.
\hfill$\Box$
\end{pf}

Now we further establish the estimates of the nonlinear terms.
\begin{lem}\label{2020201101328}
Let $0\leqslant i\leqslant3$. Under the condition \eqref{202108261502}  with sufficiently small $\delta$,
we have
\begin{align}
&\label{202108270922n12}
\|(\mathcal{N}^{\mu}_{u},\mathcal{N}^{\lambda}_{u})\|_{\underline{i},0}
\lesssim
\|\nabla\eta\|_{\underline{2},0}^{1/2}\|\nabla\eta\|_{\underline{1},1}^{1/4}\|\nabla\eta\|_{\underline{2},1}^{1/4}\|\nabla u\|_{\underline{i},0}\nonumber
\\
&\qquad \qquad\qquad\quad\ +\|\nabla\eta\|_{\underline{i},0} \|\nabla u\|_{\underline{2},0}^{1/2}\|\nabla u\|_{\underline{1},1}^{1/4}\|\nabla u\|_{\underline{1},2}^{1/4},\\[1mm]
&\label{202108270922}
\|(\mathcal{N}^{\mu}_{u},\mathcal{N}^{\lambda}_{u})\|_{\underline{1},j}
\lesssim
\begin{cases}
{\|\nabla \eta\|_{\underline{1},1}^{3/4}\|\nabla \eta\|_{\underline{2},1}^{1/4}}
\|\nabla u\|_{\underline{1},0}^{1/2}\|\nabla u\|_{\underline{1},1}^{1/2}&
\\
\;+{\|\nabla \eta\|_{\underline{1},1}^{1/2}\|\nabla \eta\|_{\underline{1},2}^{1/2}}
\|\nabla u\|_{\underline{1},0}^{3/4}\|\nabla u\|_{\underline{2},0}^{1/4}&\\
\;+\|\nabla \eta\|_{\underline{2},0}^{1/2}\|\nabla \eta\|_{\underline{1},1}^{1/4}\|\nabla \eta\|_{\underline{2},1}^{1/4}
\|\nabla u\|_{\underline{1},1} &\mbox{for }j=1;\\[1mm]
\|\nabla \eta\|_{\underline{2},0}^{1/2}\|\nabla \eta\|_{\underline{1},1}^{1/4}\|\nabla \eta\|_{\underline{2},1}^{1/4}
\|\nabla u\|_{\underline{1},2}
&\\
\;+\|\nabla \eta\|_{\underline{1},1}^{3/4}\|\nabla \eta\|_{\underline{2},1}^{1/4}
\|\nabla u\|_{\underline{1},1}^{1/2}\|\nabla u\|_{\underline{1},2}^{1/2}
\\
\;+\|\nabla \eta\|_{\underline{1},1}^{1/2}\|\nabla \eta\|_{\underline{2},1}^{1/2}
\|\nabla u\|_{\underline{1},1}^{3/4}\|\nabla u\|_{\underline{1},2}^{1/4}
\\
\; +
\|\nabla u\|_{\underline{2},0}^{1/2}\|\nabla u\|_{\underline{1},1}^{1/4}\|\nabla u\|_{\underline{2},1}^{1/4}
\|\nabla \eta\|_{\underline{1},2}&\mbox{for }j=2,
\end{cases}\\[1mm]
&\label{202408290922n12}
\|(\mathcal{N}^{\mu}_{u},\mathcal{N}^{\lambda}_{u})_{t}\|_{\underline{1},0}\lesssim
\|\nabla\eta\|_{\underline{2},0}^{1/2}\|\nabla\eta\|_{\underline{1},1}^{1/4}\|\nabla\eta\|_{\underline{2},1}^{1/4}\|\nabla u_{t}\|_{\underline{1},0}\nonumber \\
&\qquad \qquad\qquad\qquad +\|\nabla u\|_{\underline{1},0}\|\nabla u\|_{\underline{2},0}^{1/2}\|\nabla u\|_{\underline{1},1}^{1/4}\|\nabla u\|_{\underline{2},1}^{1/4},\\
&\label{202108270921}
\|\mathcal{N}^{P}_{\eta}\|_{\underline{i},0}
\lesssim
\|\nabla\eta\|_{\underline{2},0}^{1/2}\|\nabla\eta\|_{\underline{1},1}^{1/4}\|\nabla\eta\|_{\underline{2},1}^{1/4}\|\nabla \eta\|_{\underline{i},0}
,\\[1mm]
&\label{202108270932}
\|\mathcal{N}^{P}_{\eta}\|_{\underline{1},j}
\lesssim
\begin{cases}
{\|\nabla \eta\|_{\underline{1},1}^{1/2}\|\nabla \eta\|_{\underline{2},1}^{1/2}}
\|\nabla \eta\|_{\underline{1},0}^{1/2}\|\nabla \eta\|_{\underline{1},1}^{1/2}&
\\
\;+\|\nabla \eta\|_{\underline{2},0}^{1/2}\|\nabla \eta\|_{\underline{1},1}^{1/4}\|\nabla \eta\|_{\underline{2},1}^{1/4}
\|\nabla \eta\|_{\underline{1},1} &\mbox{for } j=1;\\[1mm]
{ \|\nabla \eta\|_{\underline{1},1}^{1/2}\|\nabla \eta\|_{\underline{2},1}^{1/2}}
\|\nabla \eta\|_{\underline{1},1}^{1/2}\|\nabla \eta\|_{\underline{1},2}^{1/2}&
\\
\;+\|\nabla \eta\|_{\underline{2},0}^{1/2}\|\nabla \eta\|_{\underline{1},1}^{1/4}
\|\nabla \eta\|_{\underline{2},1}^{1/4}
\|\nabla \eta\|_{\underline{1},2}
& \mbox{for } j=2,
\end{cases}\\[1mm]
&\label{202408290932n12}
\|\partial_{t}\mathcal{N}^{P}_{\eta}\|_{\underline{1},0}\lesssim
\|\nabla\eta\|_{\underline{2},0}^{1/2}\|\nabla\eta\|_{\underline{1},1}^{1/4}\|\nabla\eta\|_{\underline{2},1}^{1/4}\|\nabla u\|_{\underline{1},0}.
\end{align}
\end{lem}
\begin{pf}

(1) Let $\varphi$, $\psi\in \{\tilde{\mathcal{B}}_{ij}\}_{1\leqslant i,\ j\leqslant 3}$, and $f_{k}$ with $1\leqslant k\leqslant 5$ are defined as follows.
$$\begin{aligned}
&f_1:=\varphi,\ f_2:= J^{-1}-1 ,\  f_3:=( J^{-1}-1) \varphi  ,\\[1mm]
&f_4:=\varphi\psi\mbox{ and }f_5:=(J^{-1}-1)\varphi\psi.
\end{aligned}
$$
Following the derivations of \eqref{202109071530}--\eqref{202108270845n12m} by further using \eqref{202109071530}--\eqref{202108270845} and \eqref{202109071535}--\eqref{202108270852n12nm}, we can first have that
\begin{align}
&\label{202408290856}
 \|f_{k}\|_{\underline{i},0}
\lesssim
\|\nabla\eta\|_{\underline{i},0},
\\[1mm]
&\label{202408290854}
\|f_{k}\|_{\underline{j},1}
\lesssim
\|\nabla\eta\|_{\underline{j},1},
\\[1mm]
&\label{202408290855}
\|f_{k}\|_{\underline{1},2}
\lesssim
\|\nabla\eta\|_{\underline{1},2},\\[1mm]
 &\label{202409141119}
\|\partial_{t}f_{k} \|_{\underline{1},0}\lesssim \|\nabla u\|_{\underline{1},0},\\[0.5mm]
 &\label{202409141122}
\|\partial_{t}f_{k} \|_{L^{\infty}}\lesssim \|\nabla u\|_{L^{\infty}}
\end{align}
for any $1\leqslant k\leqslant 5$, $1\leqslant  i\leqslant 3$, and $j=1$, $2$.
Using \eqref{202408290856}--\eqref{202408290855} and \eqref{202108261428nm}, we can estimate that
\begin{align}
&
\|f_{k}\partial_{m} u_{n}\|_{\underline{i},0}
\lesssim
\|\nabla\eta\|_{\underline{2},0}^{1/2}\|\nabla\eta\|_{\underline{1},1}^{1/4}
\|\nabla\eta\|_{\underline{2},1}^{1/4}
\|\nabla u\|_{\underline{i},0}\nonumber \\
&\qquad \qquad \qquad \quad
+\|\nabla\eta\|_{\underline{i},0}\|\nabla u\|_{\underline{2},0}^{1/2}\|\nabla u\|_{\underline{1},1}^{1/4}\|\nabla u\|_{\underline{1},2}^{1/4},\label{2024dsaf193001516}  \\[1mm]
&\label{2024193001516}
\|\partial_3(f_{k}\partial_{m} u_{n})\|_{\underline{1},0}
\lesssim
\|\nabla \eta\|_{\underline{1},1}^{3/4}\|\nabla \eta\|_{\underline{2},1}^{1/4}
\|\nabla u\|_{\underline{1},0}^{1/2}\|\nabla u\|_{\underline{1},1}^{1/2}
\nonumber\\
&\qquad \qquad \qquad\quad\quad \   +\|\nabla \eta\|_{\underline{1},1}^{1/2}\|\nabla \eta\|_{\underline{1},2}^{1/2}
\|\nabla u\|_{\underline{1},0}^{3/4}\|\nabla u\|_{\underline{2},0}^{1/4}\nonumber\\
&\quad\quad\quad\quad\quad\quad\quad \quad \
+\|\nabla \eta\|_{\underline{2},0}^{1/2}\|\nabla \eta\|_{\underline{1},1}^{1/4}\|\nabla \eta\|_{\underline{2},1}^{1/4}
\|\nabla u\|_{\underline{1},1},\\[1mm]
&\label{202409061525}
\|\partial_3^2(f_{k}\partial_{m} u_{n})\|_{\underline{1},0}
\lesssim
\|\nabla \eta\|_{\underline{2},0}^{1/2}\|\nabla \eta\|_{\underline{1},1}^{1/4}\|\nabla \eta\|_{\underline{2},1}^{1/4}
\|\nabla u\|_{\underline{1},2}
\nonumber\\
&\quad\quad\quad\quad\quad\quad\quad\quad \ +
\|\nabla \eta\|_{\underline{1},1}^{3/4}\|\nabla \eta\|_{\underline{2},1}^{1/4}
\|\nabla u\|_{\underline{1},1}^{1/2}\|\nabla u\|_{\underline{1},2}^{1/2}\nonumber\\
&\quad\quad\quad\quad\quad\quad\quad\quad \
+\|\nabla \eta\|_{\underline{1},1}^{1/2}\|\nabla \eta\|_{\underline{1},2}^{1/2}
\|\nabla u\|_{\underline{1},1}^{3/4}\|\nabla u\|_{\underline{2},1}^{1/4}\nonumber\\
&\quad\quad\quad\quad\quad\quad\quad\quad \ +\|\nabla u\|_{\underline{2},0}^{1/2}\|\nabla u\|_{\underline{1},1}^{1/4}\|\nabla u\|_{\underline{2},1}^{1/4}
\|\nabla \eta\|_{\underline{1},2}.
\end{align}

Recalling the definitions of $\mathcal{N}_{u}^{\mu}$ and $\mathcal{N}_{u}^{\lambda}$ in \eqref{20201101020111209}, we see that
 \begin{align}
&\mbox{each component of }\mathcal{N}_{u}^{\mu}\mbox{ and }\mathcal{N}_{u}^{\lambda}\mbox{ is a linear combinations of }\{f_k\partial_m u_n\}. \label{202020100015016}
\end{align}
Thus we immediately derive \eqref{202108270922n12}--\eqref{202108270922} from \eqref{2024dsaf193001516}--\eqref{202409061525} and the fact stated in \eqref{202020100015016}.

Furthermore, making use of \eqref{202408290856}--\eqref{202409141122} and \eqref{202108261428nm}, we can estimate that
\begin{align}
\|(f_{k}\partial_{m} u_{n})_{t}\|_{\underline{1},0}
\lesssim&\|\nabla\eta\|_{\underline{2},0}^{1/2}\|\nabla\eta\|_{\underline{1},1}^{1/4}\|\nabla\eta\|_{\underline{2},1}^{1/4}\|\nabla u_{t}\|_{\underline{1},0}
\nonumber \\
&+\|\partial_{t}f_{k}\|_{\underline{1},0}\|\partial_{m} u_{n}\|_{L^{\infty}}
+\|\partial_{t}f_{k}\|_{L^{\infty}}\|\partial_{m} u_{n}\|_{\underline{1},0}\nonumber\\[0.5mm]
\lesssim&\|\nabla\eta\|_{\underline{2},0}^{1/2}\|\nabla\eta\|_{\underline{1},1}^{1/4}\|\nabla\eta\|_{\underline{2},1}^{1/4}\|\nabla u_{t}\|_{\underline{1},0}
+\|\nabla u\|_{\underline{1},0}\|\nabla u\|_{L^{\infty}}\nonumber,
\end{align}
which, together with \eqref{202020100015016}, yields \eqref{202408290922n12} by using \eqref{202108261424}.

In addition, noting that
\begin{align}
&\mbox{each component of }\mathcal{N}_{\eta}^{P} \mbox{ is a linear combinations of }
\tilde{\mathcal{N}}_{\eta}^{3},\ f_1\tilde{\mathcal{N}}_{\eta}^3\mbox{ and }f_1\partial_m \eta_n.
\label{2020210930158222}
\end{align}
where $\tilde{\mathcal{N}}_{\eta}^3$ is defined by \eqref{201903120816}.
Therefore, following the derivations of \eqref{202108270922n12}--\eqref{202408290922n12} by further utilizing  \eqref{202109071530}--\eqref{202108270845n12m}, \eqref{202108270856}--\eqref{202108270855n12nm} and the interpolation inequality,
the desired estimates \eqref{202108270921}--\eqref{202408290932n12}   easily follow.
This completes the proof.
\hfill$\Box$
\end{pf}

\section{Proof of Theorem \ref{thm:2109n}} \label{20201010052127}

This section is devoted to the proof of Theorem \ref{thm:2109n}. The key step in the proof is to \emph{a priori} derive the  energy inequality \eqref{202109090810} for the solution $(\eta,u)$ of the problem \eqref{20210824nnm} under the \emph{a priori} assumptions of \eqref{prio1} and \eqref{prio2}, where the constant $K$ will be determined later by \eqref{20241101031319} and the smallness of $\delta$ depends on $\bar{\rho}$, $\mu$, $\lambda$ and $\Omega$, but not on $\kappa$.
Obviously, the condition \eqref{202108261502} and the preliminary estimates in Lemmas \ref{lem:2108261504}--\ref{2020201101328}  hold automatically under the assumptions \eqref{prio1} and \eqref{prio2} with sufficiently small $\delta$.

\subsection{Basic energy estimates}

Now we are focus on the derivation of basic energy estimates for $(\eta, u)$, which
include the $y_{\mm{h}}$-derivative estimates of $(\eta, u)$ in Lemmas \ref{lem:22083001}--\ref{lem:24083001}, the estimate of temporal derivative of $u$ in Lemma \ref{12281500}, and the normal derivative estimates of $(\eta,u)$ in Lemmas \ref{lem:dfifessimM}--\ref{lem:dfifessim2M}. Next, we shall establish these estimates
in sequence.

\begin{lem}\label{lem:22083001}
Let $\alpha$ satisfy $0\leqslant |\alpha|\leqslant 3$. Under the assumptions \eqref{prio1} and \eqref{prio2} with sufficiently small $\delta$, it holds that
\begin{align}
&\label{202108301219}
\frac{1}{2}\frac{\mm{d}}{\mm{d}t} \left(\|\partial_{\mm{h}}^{\alpha}(\sqrt{\bar{\rho}}u,\sqrt{P'(\bar{\rho})\bar{\rho}} \mm{div}\eta,\sqrt{\kappa}\nabla\eta)\|_0^2\right)
+ \|\partial_{\mm{h}}^{\alpha}(\sqrt{\mu/2}\mathbb{D}u,\sqrt{\lambda}\mm{div}u)\|_0^2\nonumber\\[0.5mm]
&\lesssim (\kappa^{-1/8}+\kappa^{-5/8})\sqrt{\mathcal{E}}\mathcal{D}.
\end{align}
\end{lem}
\begin{pf}
Applying $\partial_{\mm{h}}^{\alpha}$  to \eqref{20210824nnmm} yields
\begin{equation}\label{20240824nnmm}
\begin{cases}
\partial_{\mm{h}}^{\alpha}\eta_t=\partial_{\mm{h}}^{\alpha}u&\mbox{in } \Omega ,\\[0.5mm]
\bar{\rho} \partial_{\mm{h}}^{\alpha}u_t
-\mm{div}\partial_{\mm{h}}^{\alpha}({P}'(\bar{\rho})\bar{\rho}\mm{div}\eta \mathbb{I}+\mu\mathbb{D}u+\lambda\mm{div}u\mathbb{I}
+\kappa\nabla\eta^\top)
=\mm{div}\partial_{\mm{h}}^{\alpha}({\mathcal{N}}^{\mu}_{u}+{\mathcal{N}}^{\lambda}_{u}-{\mathcal{N}}^{P}_{\eta})
&\mbox{in } \Omega ,\\[0.5mm]
\llbracket  \partial_{\mm{h}}^{\alpha}\eta\rrbracket=\llbracket \partial_{\mm{h}}^{\alpha} u\rrbracket=0,
& \mbox{on } \Sigma,\\[0.5mm]
 \llbracket
\partial_{\mm{h}}^{\alpha}({P}'(\bar{\rho})\bar{\rho}\mm{div}\eta \mathbb{I}+
\mu\mathbb{D}u+ \lambda\mm{div}u\mathbb{I}+\kappa\nabla\eta )\rrbracket \mathbf{e}^3
=\llbracket\partial_{\mm{h}}^{\alpha}({\mathcal{N}}^{P}_{\eta}-{\mathcal{N}}^{\mu}_{u}-{\mathcal{N}}^{\lambda}_{u})\rrbracket \mathbf{e}^3
& \mbox{on } \Sigma,\\[0.5mm]
(\partial_{\mm{h}}^{\alpha}\eta, \partial_{\mm{h}}^{\alpha}u)=(0,0) & \mbox{on }\Sigma_{-}^{+}.
\end{cases}
\end{equation}
Taking the inner product of \eqref{20240824nnmm}$_2$ and $\partial_{\mm{h}}^{\alpha}u$ in $L^2$, we get
\begin{align}\label{202108301410}
&\frac{1}{2}\frac{\mm{d}}{\mm{d}t}\|\sqrt{\bar{\rho}}\partial_{\mm{h}}^{\alpha}u\|_0^2
+\int\mm{div}\partial_{\mm{h}}^{\alpha}\left({P}'(\bar{\rho})\bar{\rho}\mm{div}\eta \mathbb{I}+\mu\mathbb{D}u+\lambda\mm{div}u\mathbb{I}
+\kappa\nabla\eta^\top\right)
\cdot\partial_{\mm{h}}^{\alpha}u\mm{d}y\nonumber\\
&-\int\mm{div}\partial_{\mm{h}}^{\alpha}({\mathcal{N}}^{\mu}_{u}+{\mathcal{N}}^{\lambda}_{u}-{\mathcal{N}}^{P}_{\eta})
\cdot\partial_{\mm{h}}^{\alpha}u\mm{d}y=0.
\end{align}
Exploiting the integration by parts and the boundary conditions \eqref{20240824nnmm}$_3$--\eqref{20240824nnmm}$_5$, we have
\begin{align}
&\int\mm{div}\partial_{\mm{h}}^{\alpha}\left({P}'(\bar{\rho})\bar{\rho}\mm{div}\eta \mathbb{I}+\mu\mathbb{D}u+\lambda\mm{div}u\mathbb{I}
+\kappa\nabla\eta^\top\right)
\cdot\partial_{\mm{h}}^{\alpha}u\mm{d}y\nonumber\\
&-\int\mm{div}\partial_{\mm{h}}^{\alpha}({\mathcal{N}}^{\mu}_{u}+{\mathcal{N}}^{\lambda}_{u}-{\mathcal{N}}^{P}_{\eta})
\cdot\partial_{\mm{h}}^{\alpha}u\mm{d}y\nonumber\\
&= \frac{1}{2}\frac{\mm{d}}{\mm{d}t}\|\partial_{\mm{h}}^{\alpha}(\sqrt{P'(\bar{\rho})\bar{\rho}}\mm{div}\eta,
\sqrt{\kappa}\nabla \eta)\|_0^2
+ \|\partial_{\mm{h}}^{\alpha}(\sqrt{\mu/2}\mathbb{D}u,\sqrt{\lambda} \mm{div}u)\|_0^2+\sum_{i=1}^3I_{i}  ,
\label{202108301420}
\end{align}
where we have defined that
\begin{align}
&I_{1}:=\int\partial_{\mm{h}}^{\alpha}{\mathcal{N}}^{\mu}_{u}:\nabla\partial_{\mm{h}}^{\alpha}u^{\top}\mm{d}y,
\ I_{2}:=\int\partial_{\mm{h}}^{\alpha}{\mathcal{N}}^{\lambda}_{u}\nabla \partial_{\mm{h}}^{\alpha}\mm{div}u\mm{d}y,
\mbox{ and } I_{3}:=-\int\partial_{\mm{h}}^{\alpha}{\mathcal{N}}^{P}_{\eta} \partial_{\mm{h}}^{\alpha}\mm{div}u\mm{d}y.\nonumber
\end{align}
Using the integration by parts, \eqref{202108270922} and \eqref{202108270921},  we obtain
\begin{align}\label{202108301442}
 &|I_{1}+I_{2}|\nonumber\\
&\lesssim
 \|({\mathcal{N}}^{1}_{u},\mathcal{N} ^{2}_{u})\|_{\underline{3},0}
\|\nabla u\|_{\underline{3},0} \nonumber\\
&\lesssim
\left(
\|\nabla\eta\|_{\underline{2},0}^{1/2}\|\nabla \eta\|_{\underline{1},1}^{1/4}\|\nabla \eta\|_{\underline{2},1}^{1/4}\|\nabla u\|_{\underline{3},0}+\|\nabla\eta\|_{\underline{3},0}\|\nabla u\|_{\underline{2},0}^{1/2}\|\nabla u\|_{\underline{1},1}^{1/4}\|\nabla u\|_{\underline{2},1}^{1/4} \right)
\|\nabla u\|_{\underline{3},0}\nonumber\\
&\lesssim
\kappa^{-1/8}\sqrt{\mathcal{E}}\mathcal{D}
\end{align}
and
\begin{align}\label{202108301445}
| I_{3}| &\lesssim
\|{\mathcal{N}}^{3 }_{\eta}\|_{\underline{3},0}\|\nabla u\|_{\underline{3},0}  \lesssim
\|\nabla\eta\|_{\underline{2},0}^{1/2}\|\nabla \eta\|_{\underline{1},1}^{1/4}\|\nabla \eta\|_{\underline{2},1}^{1/8}\|\nabla\eta\|_{\underline{3},0 }
\|\nabla u\|_{\underline{3},0} \lesssim\kappa^{-5/8}\sqrt{\mathcal{E}}\mathcal{D}.
\end{align}
Plugging \eqref{202108301420}--\eqref{202108301445} into \eqref{202108301410}, then yields \eqref{202108301219}.
This completes the proof.
\hfill$\Box$
\end{pf}

\begin{lem}\label{lem:24083001}
Let $\alpha$ satisfy $0\leqslant |\alpha|\leqslant 3$. Under the assumptions \eqref{prio1} and \eqref{prio2} with sufficiently small $\delta$,  it holds that
\begin{align}
&\label{202408301219n}
\frac{\mm{d}}{\mm{d}t} \left(2\int \bar{\rho}\partial_{\mm{h}}^{\alpha}u\cdot\partial_{\mm{h}}^{\alpha}\eta\mm{d}y
+\|\partial_{\mm{h}}^{\alpha}(\sqrt{\mu/2} \mathbb{D}\eta,\sqrt{\lambda}\mm{div}\eta)\|_0^2\right)
+\|\partial_{\mm{h}}^{\alpha}(\sqrt{P'(\bar{\rho})\bar{\rho}}\mm{div}\eta,\sqrt{\kappa}\nabla \eta)\|_0^2\nonumber\\[0.5mm]
&\lesssim \|\sqrt{\bar{\rho}}\partial_{\mm{h}}^{\alpha}u\|_{0}^2
+\kappa^{-1/2}(\kappa^{-1/8}+\kappa^{-5/8})\sqrt{\mathcal{E}}\mathcal{D}.
\end{align}
\end{lem}
\begin{pf}
Taking the inner product of \eqref{20240824nnmm}$_2$ and $\partial_{\mm{h}}^{\alpha}\eta$ in $L^2$ yields
\begin{align}\label{202108301606}
&\frac{\mm{d}}{\mm{d}t}\int \bar{\rho}\partial_{\mm{h}}^{\alpha}u\cdot\partial_{\mm{h}}^{\alpha}\eta\mm{d}y
+\int\mm{div}\partial_{\mm{h}}^{\alpha}\left({P}'(\bar{\rho})\bar{\rho}\mm{div}\eta \mathbb{I}+\mu\mathbb{D}u+\lambda\mm{div}u\mathbb{I}
+\kappa\nabla\eta^\top\right)
\cdot\partial_{\mm{h}}^{\alpha}u\mm{d}y\nonumber\\
&-\int\mm{div}\partial_{\mm{h}}^{\alpha}({\mathcal{N}}^{\mu}_{u}+{\mathcal{N}}^{\lambda}_{u}-{\mathcal{N}}^{P}_{\eta})
\cdot\partial_{\mm{h}}^{\alpha}u\mm{d}y=\|\sqrt{\bar{\rho}}\partial_{\mm{h}}^{\alpha}u\|_0^2.
\end{align}
Similarly to \eqref{202108301420}, exploiting the integration by parts and the boundary conditions \eqref{20240824nnmm}$_3$--\eqref{20240824nnmm}$_5$, we can have
\begin{align}
&\int\mm{div}\partial_{\mm{h}}^{\alpha}\left({P}'(\bar{\rho})\bar{\rho}\mm{div}\eta \mathbb{I}+\mu\mathbb{D}u+\lambda\mm{div}u\mathbb{I}
+\kappa\nabla\eta^\top\right)
\cdot\partial_{\mm{h}}^{\alpha}u\mm{d}y\nonumber\\
&-\int\mm{div}\partial_{\mm{h}}^{\alpha}({\mathcal{N}}^{\mu}_{u}+{\mathcal{N}}^{\lambda}_{u}-{\mathcal{N}}^{P}_{\eta})
\cdot\partial_{\mm{h}}^{\alpha}u\mm{d}y\nonumber\\
&= \frac{1}{2}\frac{\mm{d}}{\mm{d}t}\|\partial_{\mm{h}}^{\alpha}(\sqrt{\mu/2}\mathbb{D}u,\sqrt{\lambda} \mm{div}u)\|_0^2
+\|\partial_{\mm{h}}^{\alpha}(\sqrt{P'(\bar{\rho})\bar{\rho}}\mm{div}\eta,
\sqrt{\kappa}\nabla \eta)\|_0^2
+\sum_{i=4}^6I_{i}  ,
\label{202408301420}
\end{align}
where we have defined that
\begin{align}
&I_{4}:=\int\partial_{\mm{h}}^{\alpha}{\mathcal{N}}^{\mu}_{u}:\nabla\partial_{\mm{h}}^{\alpha}\eta^{\top}\mm{d}y,
\  I_{5}:=\int\partial_{\mm{h}}^{\alpha}{\mathcal{N}}^{\lambda}_{u} \partial_{\mm{h}}^{\alpha}\mm{div}\eta\mm{d}y\mbox{ and }I_{6}:=-\int\partial_{\mm{h}}^{\alpha}{\mathcal{N}}^{P}_{\eta} \partial_{\mm{h}}^{\alpha}\mm{div}\eta\mm{d}y.\nonumber
\end{align}
Following the arguments of \eqref{202108301442}--\eqref{202108301445}, we can also  estimate that
\begin{align}\label{202108301619}
\left|I_{4}+I_{5}+I_{6}\right|  \lesssim \kappa^{-1/2}(\kappa^{-1/8}+\kappa^{-5/8})\sqrt{\mathcal{E}}\mathcal{D}.
\end{align}
Putting \eqref{202408301420}--\eqref{202108301619} into \eqref{202108301606}, then \eqref{202408301219n} follows.
This completes the proof.
\hfill$\Box$
\end{pf}

\begin{lem}\label{12281500}
Under the assumptions \eqref{prio1} and \eqref{prio2} with sufficiently small $\delta$, it holds that
\begin{align}
&\kappa^{-2}\frac{\mathrm{d}}{\mathrm{d}t}\left(\|\sqrt{\bar{\rho}} u_{t}\|_{\underline{1},0}^2+
\|(\sqrt{P'(\bar{\rho})\bar{\rho}}\mm{div}u,\sqrt{\kappa}\nabla u)\|_{\underline{1},0}^2\right) + {c}{\kappa^{-2}}\|u_{t}\|_{\underline{1},1}^2\nonumber \\[1mm]
& \label{12281503}
\lesssim\left(1+\kappa^{-1}\right)\kappa^{-1/8}\sqrt{\mathcal{E}}\mathcal{D}.
\end{align}
\end{lem}
\begin{pf}
Applying $\partial_{\mm{h}}^{\alpha}\partial_{t}$ with $|\alpha|\leqslant1$ to \eqref{20210824nnmm},
we have
\begin{equation}\label{20240817nnmm}
\begin{cases}
\bar{\rho} \partial_{\mm{h}}^{\alpha}u_{tt}
-\mm{div}\partial_{\mm{h}}^{\alpha}\left({P}'(\bar{\rho})\bar{\rho}\mm{div}u \mathbb{I}+\mu\mathbb{D}u_{t}+\lambda\mm{div}u_{t}\mathbb{I}
+\kappa\nabla u^\top\right)&
\\[0.5mm]
=\mm{div}\partial_{\mm{h}}^{\alpha}\partial_{t}(\mathcal{N}^{\mu}_{u}+{\mathcal{N}}^{\lambda}_{u}-{\mathcal{N}}^{P}_{\eta})
&\mbox{in } \Omega ,\\[0.5mm]
\llbracket \partial_{\mm{h}}^{\alpha} u\rrbracket=\llbracket \partial_{\mm{h}}^{\alpha} u_{t}\rrbracket=0,
& \mbox{on } \Sigma,\\[0.5mm]
\llbracket\partial_{\mm{h}}^{\alpha}(
{P}'(\bar{\rho})\bar{\rho}\mm{div}u \mathbb{I}+\mu\mathbb{D}u_{t}+ \lambda\mm{div}u_{t}\mathbb{I}+\kappa\nabla u) \rrbracket \mathbf{e}^3
=\llbracket\partial_{\mm{h}}^{\alpha}\partial_{t}({\mathcal{N}}^{P}_{\eta}-{\mathcal{N}}^{\mu}_{u}
-{\mathcal{N}}^{\lambda}_{u})\rrbracket \mathbf{e}^3
& \mbox{on } \Sigma,\\[0.5mm]
(\partial_{\mm{h}}^{\alpha} u, \partial_{\mm{h}}^{\alpha} u_{t})=(0,0) & \mbox{on }\Sigma_{-}^{+}.
\end{cases}
\end{equation}
Taking the inner product of \eqref{20240817nnmm}$_1$ with $\kappa^{-2}\partial_{\mm{h}}^{\alpha}u_{t}$ in $L^2$,
and then integrating by parts over $\Omega$, we have
\begin{align} \label{12281503n1}
&{\kappa^{-2}}\frac{\mathrm{d}}{\mathrm{d}t}\left(\|\sqrt{\bar{\rho}} \partial_{\mm{h}}^{\alpha}u_{t}\|_{0}^2+
\|\partial_{\mm{h}}^{\alpha}(\sqrt{P'(\bar{\rho})\bar{\rho}}\mm{div}u,\sqrt{\kappa}\nabla u)\|_{0}^2\right) +{\kappa^{-2}}\|\partial_{\mm{h}}^{\alpha}(\sqrt{\mu/2}\mathbb{D} (u_{t}),\sqrt{\lambda}\mm{div}u_{t})\|_{1}^2\nonumber\\[1mm]
&=- {\kappa^{-2}}\int\left(\partial_{\mm{h}}^{\alpha}\partial_t\mathcal{N}^{\mu}_{u}:\nabla\partial_{\mm{h}}^{\alpha} u_{t}^{\top}
+\partial_{\mm{h}}^{\alpha}(\partial_t\mathcal{N}^{\lambda}_{u}
-\partial_t\mathcal{N}^{P}_{\eta})\partial_{\mm{h}}^{\alpha}\mm{div}u_{t}\right)\mm{d}y.
\end{align}
Thanks to \eqref{202408290922n12} and \eqref{202408290932n12}, we can estimate that
\begin{align*}
&-{\kappa^{-2}}\int\left(\partial_{\mm{h}}^{\alpha}\partial_t\mathcal{N}^{\mu}_{u}:\nabla \partial_{\mm{h}}^{\alpha}u_{t}^{\top}
+\partial_{\mm{h}}^{\alpha}(\partial_t\mathcal{N}^{\lambda}_{u}-\partial_t\mathcal{N}^{P}_{\eta})
\mm{div}\partial_{\mm{h}}^{\alpha}u_{t}\right)\mm{d}y\\
&\lesssim{\kappa^{-2}}\|\nabla u_{t}\|_{\underline{1},0}
\bigg(
\|\nabla\eta\|_{\underline{2},0}^{1/2}\|\nabla\eta\|_{\underline{1},1}^{1/4}\|\nabla\eta\|_{\underline{2},1}^{1/4}\|\nabla u_{t}\|_{\underline{1},0}
+\|\nabla u\|_{\underline{1},0}\|\nabla u\|_{\underline{2},0}^{1/2}
\|\nabla u\|_{\underline{1},1}^{3/8}\|\nabla u\|_{\underline{2},1}^{1/8}\bigg)\\[0.5mm]
&\quad+{\kappa^{-2}}\|\nabla\partial_{\mm{h}}^{\alpha} u_{t}\|_{0}
\|\nabla\eta\|_{\underline{2},0}^{1/2}\|\nabla\eta\|_{\underline{1},1}^{1/4}\|\nabla\eta\|_{\underline{2},1}^{1/4}\|\nabla u\|_{\underline{1},0}
\\[0.5mm]
&\lesssim{\kappa^{-2}}\|\nabla\partial_{\mm{h}}^{\alpha} u_{t}\|_{0}
\bigg({\kappa^{7/8}}\|{\kappa}^{-1/2}\nabla u\|_{\underline{1},0}\|\nabla u\|_{\underline{2},0}^{1/2}
\|{\kappa}^{-1/2}\nabla u\|_{\underline{1},1}^{1/4}\|\kappa^{-1}\nabla u\|_{\underline{2},1}^{1/4}\\
&\quad\quad\quad\quad\quad\quad\quad\;\;+\|\nabla\eta\|_{\underline{2},0}^{1/2}
\|\nabla\eta\|_{\underline{1},1}^{1/4}\|\nabla\eta\|_{\underline{2},1}^{1/4}(\|u_{t}\|_{\underline{1},1}+\|u\|_{\underline{1},1})\bigg)\\[0.5mm]
&\lesssim\left(1+\kappa^{-1}\right)\kappa^{-1/8}\sqrt{\mathcal{E}}\mathcal{D}.
\end{align*}
Inserting the above estimate into \eqref{12281503n1},  then summing up over $\alpha$ from $|\alpha|=0$ to $|\alpha|=1$, and finally using Korn's inequality yields \eqref{12281503} immediately. This completes the proof.
\hfill$\Box$
\end{pf}

\begin{lem}\label{lem:dfifessimM}
Under the assumptions of \eqref{prio1} and \eqref{prio2} with sufficiently small $\delta$, it holds that
\begin{align}
&\label{201903142019}
  \frac{\mm{d}}{\mm{d}t}\mathcal{H}_{\kappa} (\eta) + c\| (\eta, \kappa^{-1}u)\|_{\underline{1},3}^2 \nonumber\\[1mm]
&\lesssim
(1+\kappa^{-4})\|(\sqrt{\kappa}\eta,u)\|_{\underline{3},1}^2 +(1+\kappa^{-2})\left(\kappa^{-2}\|u_t\|_{\underline{1},1}^2
+(\kappa^{-1/4}+\kappa^{-5/2})\mathcal{E}\mathcal{D}\right),
\end{align}
 where the energy functional $\mathcal{H}_{\kappa} (\eta)$ satisfies
\begin{align}\label{Hetaprofom}
 {\kappa}^{-1} \|\partial_3^2\eta \|_{\underline{1},1}^2 \lesssim {\mathcal{H}}_{\kappa} (\eta)
\lesssim \left(1+{\kappa}^{-2}\right){\kappa}^{-1} \|\partial_3^2\eta \|_{\underline{1},1}^2.
\end{align}
\end{lem}
\begin{rem}
By virtue of \eqref{Hetaprofom}, we can see that
\begin{align*}
\kappa^{-1} \|\eta \|_{\underline{1},3}^2
\lesssim {\mathcal{H}}_{\kappa} (\eta)+ (1+ {\kappa}^{-2})\kappa\|\eta\|_{\underline{3},1}^2
\lesssim (1+ {\kappa}^{-2}) \big(\kappa^{-1} \|\eta \|_{\underline{1},3}^2+\kappa\|\eta\|_{\underline{3},1}^2\big).
\end{align*}
\end{rem}
\begin{pf}
The first two components and the third component of \eqref{20210824nnmm}$_2$ read as
\begin{align}
& \label{Stokesequson1137}
-\mu \partial_3^2 u_\mm{h}-\kappa\partial_3^2\eta_\mm{h}=\mathcal{K}_{\mm{h}} ,\\[1mm]
& \label{Stokesequson1}
-(\mu +\lambda)\partial_3^2 u_3
-(P'(\bar{\rho})\bar{\rho}+\kappa)  \partial_3^2\eta_3= \mathcal{K}_{3} ,
\end{align}
where we have defined that
\begin{align}
&\label{Stokesequson1137n}
\mathcal{K}_{\mm{h}}:= \mu  \Delta_\mm{h} u_\mm{h} + {\lambda} \nabla_\mm{h} \mm{div} u
 +\kappa\Delta_\mm{h}\eta_\mm{h} +P'(\bar{\rho})\bar{\rho}\nabla_\mm{h}\mm{div}\eta-\bar{\rho} \partial_tu_{\mm{h}}\nonumber\\[1mm]
&\quad\quad\;
+(\mm{div}\partial_{\mm{h}}^{\alpha}({\mathcal{N}}^{\mu}_{u}+{\mathcal{N}}^{\lambda}_{u}-{\mathcal{N}}^{P}_{\eta}))_\mm{h},\\[1mm]
 & \label{Stokesequson1n}
\mathcal{K}_{3}:= \mu \Delta_\mm{h} u_3
+\lambda\partial_3\mm{div}_{\mm{h}}u_{\mm{h}}
+\kappa\Delta_\mm{h} \eta_3+P'(\bar{\rho})\bar{\rho}\partial_3\mm{div}_{\mm{h}}\eta_{\mm{h}} -\bar{\rho}\partial_t u_3\nonumber\\[1mm]
&\quad\quad\;
+ (\mm{div}\partial_{\mm{h}}^{\alpha}({\mathcal{N}}^{\mu}_{u}+{\mathcal{N}}^{\lambda}_{u}-{\mathcal{N}}^{P}_{\eta}))_3 ,
\end{align}
and
$f_\mm{h}$, $f_3 $ represnet the first two components and the third component of $f:=(f_\mm{h},f_3)^\top$, resp.
Noting that the order of $\partial_3$ in the linear parts on the right-hand side of \eqref{Stokesequson1137} and \eqref{Stokesequson1} is lower than that of $\partial_3$ on the left-hand side,
this fact clearly provides that the $y_3$-derivative estimates of $\eta$ can be converted to the $y_{\mm{h}}$-derivative
estimates.

Let $0\leqslant i$, $l\leqslant 1$, $0\leqslant j\leqslant i$ and  $0\leqslant k\leqslant i-j$.
Applying $ \partial_{\mm{h}}^{j+k+l}\partial_3^{i-j-k}$ to \eqref{Stokesequson1137}--\eqref{Stokesequson1} yields that
\begin{align}
& \label{Stesen1137}
-\partial_{\mm{h}}^{j+k+l}\partial_3^{i-j-k}(\mu \partial_3^2 u_\mm{h}+\kappa\partial_3^2\eta_\mm{h})=\partial_{\mm{h}}^{j+k+l}\partial_3^{i-j-k}\mathcal{K}_{\mm{h}} ,\\[2mm]
& \label{Stosequson1}
-\partial_{\mm{h}}^{j+k+l}\partial_3^{i-j-k}((\mu +\lambda)\partial_3^2 u_3
+(P'(\bar{\rho})\bar{\rho}+\kappa)  \partial_3^2\eta_3)=\partial_{\mm{h}}^{j+k+l}\partial_3^{i-j-k} \mathcal{K}_{3} .
\end{align}
Taking the inner products of $-\kappa^{-1}\partial_{\mm{h}}^{j+k+l}\partial_3^{i-j-k+2} \eta_{\mm{h}}$, resp. $-\kappa^{-1}\partial_{\mm{h}}^{j+k+l}\partial_3^{i-j-k+2} \eta_3$ with \eqref{Stesen1137} resp.
\eqref{Stosequson1} in $L^2$,
we find by $\eta_t=u$ that
\begin{align} & {\kappa}^{-1}\int
\partial_{\mm{h}}^{j+k+l}\partial_3^{i-j-k+2}
\left( \mu \partial_t \eta_\mm{h}
+\kappa\eta_\mm{h}\right)
\cdot\partial_{\mm{h}}^{j+k+l}\partial_3^{i-j-k+2} \eta_{\mm{h}}\mm{d}y\nonumber \\
&= -{\kappa}^{-1}\int
\partial_{\mm{h}}^{j+k+l}\partial_3^{i-j-k}\mathcal{K}_{\mm{h}}
\cdot\partial_{\mm{h}}^{j+k+l}\partial_3^{i-j-k+2} \eta_{\mm{h}}  \mm{d}y \label{201912021540}
\end{align}
and
\begin{align}\label{201912021542}
&{\kappa}^{-1}\int
 \partial_{\mm{h}}^{j+k+l}\partial_3^{i-j-k}\left((\mu+\lambda) \partial_3^2
\partial_t\eta_3 +(P'(\bar{\rho})\bar{\rho} +\kappa)\partial_3^2\eta_3\right)
\partial_{\mm{h}}^{j+k+l}\partial_3^{i-j-k+2}\eta_3
\mm{d}y\nonumber\\
& =-{\kappa}^{-1}\int \partial_{\mm{h}}^{j+k+l}\partial_3^{i-j-k}\mathcal{K}_{3}
\partial_{\mm{h}}^{j+k+l}\partial_3^{i-j-k+2}\eta_3\mm{d}y.
\end{align}
Using Young's inequality,
we can further deduce from \eqref{201912021540}--\eqref{201912021542} that
\begin{equation}
\label{omdm12dfs2MJnm}\begin{aligned}
&{\kappa}^{-1}\frac{\mm{d}}{\mm{d}t}\left(\mu \| \partial_3^2 \eta_\mm{h}\|_{j+l,i-j}^2+(\mu +\lambda)\| \partial_3^2\eta_3
\|_{j+l,i-j}^2\right)
+\frac{1}{2}\|\partial_3^2 \eta_\mm{h} \|_{j+l,i-j}^2\\[1mm]
&+(1+\kappa^{-1}P'(\bar{\rho})\bar{\rho})
\left\|
\partial_3^2\eta_3\right\|_{j+l,i-j}^2
\lesssim {\kappa^{-2}}\left\|\mathcal{K}  \right\|_{j+l,i-j}^2.
 \end{aligned}   \end{equation}

Similarly,  taking the inner products of  $-\kappa^{-2}\partial_{\mm{h}}^{j+k+l}\partial_3^{i-j-k+2} u_{\mm{h}}$, resp.
$-\kappa^{-2}\partial_{\mm{h}}^{j+k+l}\partial_3^{i-j-k+2} u_3$ with \eqref{Stesen1137}, resp.
\eqref{Stosequson1} in $L^2$, and then
following the derivation of \eqref{omdm12dfs2MJnm}, we can obtain that
 \begin{align}
 &{\kappa}^{-1} \frac{\mm{d}}{\mm{d}t}
 \left(\left\|\partial_3^2 \eta_{\mm{h}}
 \right\|_{j+l,i-j}^2+ \left(1+\kappa^{-1}P'(\bar{\rho})\bar{\rho}\right)\left\|  \partial_3^2\eta_3
\right\|_{j+l,i-j}^2 \right)\nonumber \\[1mm]
&+{c}{\kappa^{-2}}\|(\sqrt{\mu}\partial_3^2u_{\mm{h}},
\sqrt{\left(\mu+\lambda\right)}\partial_3^2u_3 )\|_{j+l,i-j}^2 \lesssim {\kappa^{-2}}\left\|\mathcal{K}  \right\|_{j+l,i-j}^2.
\label{201910021625m}
 \end{align}
Thus, combining \eqref{omdm12dfs2MJnm} with \eqref{201910021625m}, we ontain
\begin{align}     \label{201701231114}
\begin{aligned}
  \frac{\mm{d}}{\mm{d}t} {\mathcal{H}}_{i,j,\kappa}(\eta) + c  \|\partial_3^2 (\eta ,\kappa^{-1}u ) \|_{j+l,i-j}^2
  \lesssim   {\kappa^{-2}} \left\|\mathcal{K} \right\|_{j+l,i-j}^2,
\end{aligned}
\end{align}
where we have defined that
 $$ {\mathcal{H}}_{i,j,\kappa}(\eta):= {\kappa}^{-1}\left\|\left(\sqrt{(\mu+1)}\partial_3^2 \eta_{\mm{h}},
 \sqrt{(1+\mu+\lambda+ \kappa^{-1}P'(\bar{\rho})\bar{\rho})}\partial_3^2 \eta_3\right) \right\|_{j+l,i-j}^2.$$

In addition, we have
\begin{equation*} \left\|\mathcal{K} \right\|_{j+l,i-j}^2
\lesssim  \| (\sqrt{\kappa}\eta, \sqrt{P'(\bar{\rho})\bar{\rho}}\eta,  u)\|_{j+l+1,i-j+1}^2+\| ( u_t,\mm{div}(\mathcal{N}^{\mu}_{u}+{\mathcal{N}}^{\lambda}_{u}-{\mathcal{N}}^{P}_{\eta})) \|_{j+l,i-j}^2,
\end{equation*}
from which and \eqref{201701231114} we arrive at
\begin{align}
&\frac{\mm{d}}{\mm{d}t}{\mathcal{H}}_{\kappa} (\eta)+   \| (\eta ,\kappa^{-1}u ) \|_{\underline{1},3}^2\nonumber \\
&\lesssim \| (\eta ,\kappa^{-1}u ) \|_{\underline{3},1}^2+
  \left(1+{\kappa}^{-1}\right)^2\left( {\kappa^{-2}}\|(\sqrt{\kappa}\eta,\eta,u) \|_{\underline{3}, 1}^2
   +{\kappa^{-2}}\|(u_t,\mm{div}(\mathcal{N}^{\mu}_{u}+{\mathcal{N}}^{\lambda}_{u}-{\mathcal{N}}^{P}_{\eta})) \|_{\underline{1},1}^2\right)\nonumber \\[1mm]
&\lesssim
  \left(1+ {\kappa}^{-4}\right)\|(\sqrt{\kappa}\eta,u)\|_{\underline{3}, 1}^2
   +\left(1+ {\kappa}^{-2}\right)\left( {\kappa^{-2}}\|u_t\|_{\underline{1},1}^2
   +{\kappa^{-2}}\|(\mathcal{N}^{\mu}_{u}+{\mathcal{N}}^{\lambda}_{u}-{\mathcal{N}}^{P}_{\eta}) \|_{\underline{1},2}^2\right),\label{omdm12dfs2nm}
\end{align}
where we have defined that
$${\mathcal{H}}_{\kappa}(\eta):=\sum_{l=0}^1\left({h}_{1,0} {\mathcal{H}}_{1,0,\kappa}(\eta)
+\left(1+{\kappa}^{-1}\right)^2{h}_{1,1} {\mathcal{H}}_{1,1,\kappa}(\eta)+{h}_{0,0} {\mathcal{H}}_{0,0,\kappa}(\eta)\right)$$ for some positive constants ${h}_{i,j}$ depending on the domain and other known physical parameters, but not on $\kappa$.

Moreover, utilizing \eqref{202108270922}--\eqref{202108270932} with $j=2$, we can estimate that
\begin{align*}
&{\kappa^{-2}}\|(\mathcal{N}^{\mu}_{u}+{\mathcal{N}}^{\lambda}_{u}-{\mathcal{N}}^{P}_{\eta}) \|_{\underline{1},2}^2\\[0.5mm]
&\lesssim\kappa^{-2}\bigg(
\|\nabla \eta\|_{\underline{2},0}\|\nabla \eta\|_{\underline{1},1}^{1/2}\|\nabla \eta\|_{\underline{2},1}^{1/2}
\|\nabla u\|_{\underline{1},2}^2
+ \|\nabla \eta\|_{\underline{1},1}^{3/2}\|\nabla \eta\|_{\underline{2},1}^{1/2}
\|\nabla u\|_{\underline{1},1}\|\nabla u\|_{\underline{1},2}\\[0.5mm]
&\quad\
+ \|\nabla \eta\|_{\underline{1},1}\|\nabla \eta\|_{\underline{2},1}
\|\nabla u\|_{\underline{1},1}^{3/2}\|\nabla u\|_{\underline{1},2}^{1/2}
+\|\nabla u\|_{\underline{2},0}\|\nabla u\|_{\underline{1},1}^{1/2}\|\nabla u\|_{\underline{2},1}^{1/2}
\|\nabla \eta\|_{\underline{1},2}^2\\[0.5mm]
&\quad\
+{ \|\nabla \eta\|_{\underline{1},1}\|\nabla \eta\|_{\underline{2},1}}
\|\nabla \eta\|_{\underline{1},1}\|\nabla \eta\|_{\underline{1},2}
+\|\nabla \eta\|_{\underline{2},0}\|\nabla \eta\|_{\underline{1},1}^{1/2}\|\nabla \eta\|_{\underline{2},1}^{1/2}
\|\nabla \eta\|_{\underline{1},2}^2\bigg) \\[0.5mm]
&\lesssim(\kappa^{-1/4}+\kappa^{-5/2})\mathcal{E}\mathcal{D}.
\end{align*}
Consequently, plugging the above estimate into \eqref{omdm12dfs2nm},
we deduce \eqref{201903142019} from \eqref{omdm12dfs2nm}.
In addition, ${\mathcal{H}}_{\kappa} (\eta)$ obviously satisfies \eqref{Hetaprofom}.
This completes the proof.
\hfill$\Box$
\end{pf}
\begin{lem}\label{lem:dfifessim2M}
Under the assumptions of \eqref{prio1} and \eqref{prio2} with
sufficiently small $\delta$, it holds that
\begin{align}
\label{202403142019}
\kappa^{-2}\|u\|_{\underline{1},2}^2
&\lesssim\kappa^{-2}\|u_{t}\|_{\underline{1},0}^2+(1+\kappa^{-2})\|\eta\|_{\underline{1},2}^2
+{ (\kappa^{-1/4}+\kappa^{-9/4})\mathcal{E}^2}.
\end{align}
\end{lem}
\begin{pf}
Let us consider the following stratified Lam\'e problem
\begin{equation}\label{20210824nm12}
\begin{cases}
-\mu\Delta u-(\mu+\lambda)\nabla\mm{div}u=F^1\quad
&\mbox{in } \Omega ,\\[0.5mm]
\llbracket  u\rrbracket=0& \mbox{on } \Sigma,\\[0.5mm]
\llbracket\mu\mathbb{D}u+ \lambda\mm{div}u\mathbb{I}\rrbracket \mathbf{e}^3=\llbracket F^2\rrbracket \mathbf{e}^3
& \mbox{on } \Sigma,\\[0.5mm]
 u=0 & \mbox{on }\Sigma_{-}^{+},
\end{cases}
\end{equation}
where we have defined that
\begin{align*}
&F^1:=-\bar{\rho} u_t
+\mm{div}\left({P}'(\bar{\rho})\bar{\rho}\mm{div}\eta \mathbb{I}
+\kappa\nabla\eta+{\mathcal{N}}^{\mu}_{u}+{\mathcal{N}}^{\lambda}_{u}-{\mathcal{N}}^{P}_{\eta}\right),\\[0.5mm]
&F^2:=-({P}'(\bar{\rho})\bar{\rho}\mm{div}\eta \mathbb{I}
+\kappa\nabla\eta+{\mathcal{N}}^{\mu}_{u}+{\mathcal{N}}^{\lambda}_{u})+{\mathcal{N}}^{P}_{\eta}.
\end{align*}

Exploiting the regularity of solutions of stratified Lam\'e problem in \eqref{Ellipticestimate},  we obtain that
\begin{align}\label{202403142019n}
\kappa^{-2}\|u\|_{\underline{1},2}^2&\lesssim\kappa^{-2}\left(\|F^{1}\|_{\underline{1},0}^2+|F^2|_{3/2}^2\right)
\nonumber\\[0.5mm]
&\lesssim\kappa^{-2}\|u_{t}\|_{\underline{1},0}^2+(1+\kappa^{-2})\|\eta\|_{\underline{1},2}^2
+\kappa^{-2}\|({\mathcal{N}}^{\mu}_{u},{\mathcal{N}}^{\lambda}_{u}, {\mathcal{N}}^{P}_{\eta})\|_{\underline{1},1}^2.
\end{align}
Making use of \eqref{202108270922} and \eqref{202108270932} with $j=1$, we can estimate that
\begin{align}
&\kappa^{-2}\|({\mathcal{N}}^{\mu}_{u},{\mathcal{N}}^{\lambda}_{u}, {\mathcal{N}}^{P}_{\eta})\|_{\underline{1},1}^2\nonumber\\[0.5mm]
&\lesssim\kappa^{-2}\big(
{ \|\nabla \eta\|_{\underline{1},1}^{3/2}\|\nabla \eta\|_{\underline{2},1}^{1/2}}
\|\nabla u\|_{\underline{1},0}\|\nabla u\|_{\underline{1},1}
+{ \|\nabla \eta\|_{\underline{1},1}\|\nabla \eta\|_{\underline{1},2}}
\|\nabla u\|_{\underline{1},0}^{3/2}\|\nabla u\|_{\underline{2},0}^{1/2}\nonumber\\[1mm]
&\quad+\|\nabla \eta\|_{\underline{2},0}\|\nabla \eta\|_{\underline{1},1}^{1/2}\|\nabla \eta\|_{\underline{2},1}^{1/2}
\|\nabla u\|_{\underline{1},1}^2
+{ \|\nabla \eta\|_{\underline{1},1}\|\nabla \eta\|_{\underline{2},1}}
\|\nabla \eta\|_{\underline{1},0}\|\nabla \eta\|_{\underline{1},1}\nonumber\\[1mm]
&\quad
+\|\nabla \eta\|_{\underline{2},0}\|\nabla \eta\|_{\underline{1},1}^{1/2}\|\nabla \eta\|_{\underline{2},1}^{1/2}
\|\nabla \eta\|_{\underline{1},1}^2\big)\nonumber\\[1mm]
&\lesssim
(\kappa^{-1/4}+\kappa^{-9/4})\mathcal{E}^2\nonumber.
\end{align}
Putting the above estimate into \eqref{202403142019n} and using Young's inequality then yields \eqref{202403142019} immediately.
This completes the proof.
\hfill$\Box$
\end{pf}

\subsection{Stability estimates}
With Lemmas \ref{lem:22083001}--\ref{lem:dfifessim2M} in hand, we are now in a position to establish the desired  \emph{a priori} estimate \eqref{pr} under the assumptions \eqref{prio1} and \eqref{prio2} for sufficiently small $\delta$.

We first define
\begin{align}
\bar{\mathcal{E}}:=&C_1\|(\sqrt{\bar{\rho}}u,\sqrt{P'(\bar{\rho})\bar{\rho}}\mm{div}\eta,
\sqrt{\kappa}\nabla\eta)\|_{\underline{3},0}^2\nonumber\\
& +2\sum_{|\alpha|=0}^3\int\bar{\rho}\partial_{\mm{h}}^{\alpha}u\cdot\partial_{\mm{h}}^{\alpha}\eta\mm{d}y+
\|(\sqrt{\mu/2}\mathbb{D}\eta,\sqrt{\lambda}\mm{div}\eta)\|_{\underline{3},0}^2,\nonumber\\
 \bar{\mathcal{D}}:=&C_1\|   u\|_{\underline{3},1}^2 +\|\mm{div}\eta \|_{\underline{3},0}^2+  {\kappa}\|\eta\|_{\underline{3},1}^2\nonumber
\end{align}
for some positive constant $C_1$.
Using Korn's inequality \eqref{psfxdxx} and Poinc\'are's inequality, we deduce from Lemmas \ref{lem:22083001}--\ref{lem:24083001} that
\begin{align}\label{202109021542}
\frac{1}{2}\frac{\mm{d}}{\mm{d}t}\bar{\mathcal{E}}+c\bar{\mathcal{D}}
 \lesssim(1+\kappa^{-1/2})(\kappa^{-1/8}+\kappa^{-5/8})\sqrt{\mathcal{E}}\mathcal{D},
\end{align}
and that $\bar{\mathcal{E}}$ satisfies
$$
\| u\|_{\underline{3},0}^2+(1+\kappa) \|\eta\|_{\underline{3},1}^2\lesssim\bar{\mathcal{E}}
\lesssim \| u\|_{\underline{3},0}^2+(1+\kappa) \|\eta\|_{\underline{3},1}^2.
$$
Hence, we further deduce from Lemmas \ref{12281500}--\ref{lem:dfifessimM} and \eqref{202109021542} that
\begin{align}
  \frac{\mm{d}}{\mm{d}t}\tilde{\mathcal{E}}
+c \tilde{\mathcal{D}}  \lesssim  (1+\kappa^{-9/2})(\kappa^{-1/8}+\kappa^{-5/2})(\sqrt{\mathcal{E}}+\mathcal{E})\mathcal{D},\label{202109021724}
\end{align}
where we have defined that
\begin{align}
\tilde{\mathcal{E}}:=&
\mathcal{H}_{\kappa} (\eta)+{C_2(1+\kappa^{-4})}{\kappa^{-2}}\left(\|\sqrt{\bar{\rho}} u_{t}\|_{\underline{1},0}^2+
\|(\sqrt{P'(\bar{\rho})\bar{\rho}}\mm{div}u,\sqrt{\kappa}\nabla u)\|_{\underline{1},0}^2\right)
+C_2^2 (1+\kappa^{-4})\bar{\mathcal{E}},\nonumber\\[1mm]
\tilde{\mathcal{D}}:=& \| (\eta, \kappa^{-1}u)\|_{\underline{1},3}^2
+ {C_2(1+\kappa^{-4})}{\kappa^{-2}}\|u_{t}\|_{\underline{1},1}^2
+C_2^2 (1+\kappa^{-4})\bar{\mathcal{D}}
 \nonumber
\end{align}
for some positive constant $C_2$.

Using the interpolation inequality, it is easy to see that
 \begin{align*}
&\|\eta \|_{\underline{1},2}^2\lesssim\| \kappa^{-1/2}\eta \|_{\underline{1},3}\|\kappa^{1/2}\eta\|_{\underline{1},1},\ {\kappa}^{1/2}\|\eta \|_{\underline{1},2}^2\lesssim\|\eta \|_{\underline{1},3}\|\kappa^{1/2}\eta\|_{\underline{1},1},\\[2mm]
&\mbox{ and } {\kappa}^{-1}\| u \|_{\underline{1},2}^2\lesssim\kappa^{-1}\|u \|_{\underline{1},3}\| u\|_{\underline{1},1}.
\end{align*}
These three inequalities, in combination with \eqref{Hetaprofom}, yields
\begin{align}
\label{202201011001400}
\mathcal{E}-\kappa^{-2}\|u\|_{\underline{1},2}^2\lesssim \tilde{\mathcal{E}}\lesssim (1+\kappa^{-4})(\mathcal{E}-\kappa^{-2}\|u\|_{\underline{1},2}^2)
\; \mbox{ and }\; \mathcal{D}\lesssim \tilde{\mathcal{D}}.
\end{align}

Therefore, it follows from \eqref{202109021724} and $\mathcal{D}\lesssim \tilde{\mathcal{D}}$ that
there exists a constant $\tilde{\delta}_1\in(0,1)$ such that, for any $\delta\leqslant \tilde{\delta}_1$,
\begin{align}
\label{202109030732}
\frac{\mm{d}}{\mm{d}t}\tilde{\mathcal{E}}
+c{\mathcal{D}}\leqslant0.
\end{align}
Integrating the above inequality over $(0,t)$, and then using \eqref{202201011001400}, we get
\begin{align*}
 \mathcal{E}(t)-\kappa^{-2}\|u\|_{\underline{1},2}^2
+\int_0^{t} \mathcal{D}(\tau)\mm{d}\tau\lesssim (1+\kappa^{-4})\mathcal{E}(0),
\end{align*}
which, together with \eqref{202403142019}, yields that
there exists a constant $\tilde{\delta}_2\in(0,1)$ such that, for any $\delta\leqslant \tilde{\delta}_2$,
\begin{align}\label{202409030732}
 {\mathcal{E}}(t)
+\int_0^{t} {\mathcal{D}}(\tau)\mm{d}\tau\lesssim (1+\kappa^{-6})\mathcal{E}(0),
\end{align}
where we have defined $f(0):=f(t)\big|_{t=0}$.

Furthermore, exploiting \eqref{202108270922} and \eqref{202108270932} with $i=1$ for $t=0$,  we deduce from \eqref{20210824nnmm}$_2$ that
there exists a constant $\tilde{\delta}_3\in(0,1)$ such that, for any $\delta\leqslant \tilde{\delta}_3$,
\begin{align}\label{202409030932}
 \kappa^{-2}\|u_{t}|_{t=0}\|_{\underline{1},0}^2
&\lesssim {\kappa^{-2}}\big((1+\kappa^2)\|\eta^0\|_{\underline{1},2}^2+\|u^0\|_{\underline{1},2}^2
+\|({\mathcal{N}}^{\mu}_{u}+{\mathcal{N}}^{\lambda}_{u}-{\mathcal{N}}^{P}_{\eta})(0)
\|_{\underline{1},1}^2\big)\nonumber\\[1mm]
&\lesssim(1+{\kappa}^{-2}){E}^0.
\end{align}
Substituting the above estimate into \eqref{202409030732} yields that, for any $\delta\leqslant \delta_1:=\min\{\tilde{\delta}_1,\tilde{\delta}_2,\tilde{\delta}_3\}$,
\begin{align}\label{20240903073n2}
\mathcal{E}(t)+\int_0^{t} {\mathcal{D}}(\tau)\mm{d}\tau
\lesssim (1+{\kappa}^{-8}){E}^0\leqslant c_2(1+{\kappa}^{-8}){E}^0\mbox{ for any }t\in [0,T],
\end{align}
where the constant $c_2\geqslant1$.

Now, if we take
\begin{align}
K^2=4c_2(1+{\kappa}^{-8}){E}^0 \label{20241101031319},
\end{align}
it then follows from \eqref{20240903073n2} that
\begin{align}\label{202109030803}
 {\mathcal{E}}(t)
+\int_0^{t} {\mathcal{D}}(\tau)\mm{d}\tau\leqslant K^2/4
\mbox{ for any }t\in [0,T].
\end{align}

Additionally, by \eqref{202201011001400}, one has
\begin{align}
 \tilde{\mathcal{E}}\leqslant c_3(1+\kappa^{-1}){\mathcal{D}}\quad\mbox{for some constant $c_3>0$.}\nonumber
\end{align}
Thanks to \eqref{202409030932} and the above relation, we derive from \eqref{202201011001400}--\eqref{202109030732} that,
for any $t\in [0,T]$,
\begin{align}\label{202409031000}
& \mathcal{E} (t)-\kappa^{-2}\|u\|_{\underline{1},2}^2
 + \int_{0}^{t} {\mathcal{D}}(\tau)e^{\frac{(\tau-t)}{c_3 (1+\kappa^{-1})}}\mm{d}\tau
\lesssim e^{-\frac{t}{c_3 (1+\kappa^{-1})}}(1+{\kappa}^{-6}){E}^0.
\end{align}
Consequently, combining \eqref{202403142019} with \eqref{202109030803}--\eqref{202409031000} and \eqref{prio2} (with $\delta\leqslant\delta_1$),
we finally obtain that
\begin{align}\label{202409031200}
& \mathcal{E} (t)
 + \int_{0}^{t} {\mathcal{D}}(\tau)e^{\frac{(\tau-t)}{c_3 (1+\kappa^{-1})}}\mm{d}\tau
\lesssim e^{-\frac{t}{c_3 (1+\kappa^{-1})}}(1+{\kappa}^{-8}){E}^0.
\end{align}
This completes the derivation of \emph{a priori} stability estimate. More precisely, we have the following conclusion of the \emph{a priori} stability estimate.
\begin{pro}
\label{202020101202031308}
Let $(\eta,u)$ be the solution of the initial-boundary value problem \eqref{20210824nnm}.  If $(\eta,u)$ satisfies \eqref{prio1}--\eqref{prio2} for $K$ defined by \eqref{20241101031319} and for sufficiently small $\delta$ such that $\delta\leqslant\delta_1$, then $(\eta,u)$ satisfies \eqref{202109030803} and the exponential stability estimate \eqref{202409031200}
for any $t\in [0,T]$.
\end{pro}

In addition, we introduce a local well-posedness for the initial-boundary value problem
\eqref{20210824nnm}.
\begin{pro}\label{pro:210903}
Let $(\eta^0,u^0)\in H^{1,3}_{ 0}\times H^{1,2}_{0}$ satisfy the compatibility condition \eqref{202409091654nn} and $$\|(\nabla\eta^0,u^0)\|_{\underline{1},2}\leqslant B,$$
where $B$ is a positive constant.
Then there is a constant $\delta_2\in (0,1]$ such that, if $\eta^0$ further satisfies
\begin{align}
\|\nabla \eta^0\|_{L^{\infty}} \leqslant \delta_2, \label{201912261426}
\end{align}
 there exist a local existence time $T>0$ (depending possibly on $B$, $\mu$, $\lambda$, $\kappa$ and $\delta_2$) and a unique local strong solution
$(\eta, u)\in C^0([0,T],  H^{1,3}_{0}\times H^{1,2}_{0})$
to the initial-boundary value problem   \eqref{20210824nnm}, satisfying
$0<\displaystyle\inf_{(y,t)\in \mathbb{R}^2\times [0,T]} \det(\nabla \eta+I)$ and $\displaystyle\sup_{t\in [0,T]}\|\nabla \eta\|_{L^{\infty}}\leqslant 2\delta_2$.
\end{pro}
\begin{pf}
The local and global well-posedness results of  stratified compressible viscous fluids have been established by Jang--Tice--Wang  in \cite{JJTIIAWangYJC,JJTIWYJTC}.
Thus, following the standard iteration method in \cite{JJTIIAWangYJC}, we can easily obtain the  well-posedness result
of the initial-boundary value problem  \eqref{20210824nnm} stated in Proposition \ref{pro:210903}.
\hfill $\Box$
\end{pf}

With Propositions \ref{202020101202031308} and \ref{pro:210903} in hand and the homeomorphism mapping theorem in Lemma \ref{pro:040s1nasfxdxx}, we easily establish Theorem \ref{thm:2109n}, please refer to \cite[Theorem 2.1]{JFJSGSfdsa} or \cite[Theorem 1.4]{JFJSsdafsafJMFMOSERT} for the proof.

\section{Proof of Theorem \ref{thm:210902}}\label{asymptotic}

Let  $(\eta,u)$ be the solution to the problem \eqref{20210824nnm} established in Theorem \ref{thm:2109n}.
By \eqref{202109091654}, we find that
\begin{align}\label{202408201622}
\|\nabla \eta^0\|_{L^{\infty}}\lesssim1,\ \|\nabla \eta^0\|_{\underline{3},0}\|\nabla \eta^0\|_{\underline{1},1}\lesssim1
 \mbox{ and }\|\nabla \eta^0\|_{\underline{3},0}^{1/2}\|\nabla \eta^0\|_{\underline{1},1}^{1/4}\|\nabla \eta^0\|_{\underline{2},1}^{1/4}\lesssim1,
\end{align}
where $(\eta^0,u^0)$ is the initial data of $(\eta,u)$ given in Theorem \ref{thm:2109n}.

Since the initial data $(\eta^0,u^0)$  does not generally satisfy the compatibility jump condition \eqref{20210824nnmmn}$_4$ of the linear problem \eqref{20210824nnmmn}, one has to modify the initial data $(\eta^0,u^0)$ so that the modified initial data, denoted by $(\tilde{\eta}^0,\tilde{u}^0)$, can be used as the initial data for the linear problem \eqref{20210824nnmmn}.

To this purpose, we consider the following stratified Lam\'e problem:
\begin{equation}\label{20210824n1mn}
\begin{cases}
-\mu\Delta u^{\mathrm{r}}-(\mu+\lambda)\nabla \mm{div}^{\mathrm{r}}=0&\mbox{in }  \Omega,\\[0.5mm]
\llbracket  u^{\mathrm{r}}\rrbracket=0,
& \mbox{on } \Sigma,\\
\llbracket \mu\mathbb{D}u^{\mathrm{r}}+\lambda\mm{div}u^{\mathrm{r}}\mathbb{I}\rrbracket \mathbf{e}^3
=-\llbracket{\mathcal{N}}^{P}_{\eta}(0)-{\mathcal{N}}^{\mu}_{u}(0)-{\mathcal{N}}^{\lambda}_{u}(0)\rrbracket \mathbf{e}^3
& \mbox{on } \Sigma,\\[0.5mm]
u^{\mathrm{r}}=0 & \mbox{on }\Sigma_{-}^{+}.
\end{cases}
\end{equation}
In view of the existence theory of stratified Lam\'e problem (see Lemma \ref{prosfxdxx}),  there exists a solution $u^{\mathrm{r}}\in {H}_{0}^{1,2}$ to \eqref{20210824n1mn}; moreover, by \eqref{Ellipticestimate}, the solution enjoys that
 \begin{align}\label{202408191959n}
\|u^{\mathrm{r}}\|_{\underline{1},2}^2
\lesssim& \left|\llbracket{\mathcal{N}}^{P}_{\eta}(0)-{\mathcal{N}}^{\mu}_{u}(0)-{\mathcal{N}}^{\lambda}_{u}(0)\rrbracket\right|_{1+1/2}^2.
\end{align}
Making use of  \eqref{202109071530}--\eqref{202108270837},
\eqref{202109071532}--\eqref{202108270852}, \eqref{202108270856}--\eqref{202108270854} for $t=0$, \eqref{202408201622} and
\eqref{embed3}--\eqref{202108261402}, we can deduce from \eqref{202408191959n} that the solution $u^{\mathrm{r}}$ enjoys the estimate \eqref{202408191959}.

Let $\tilde{\eta}^0=\eta^0$ and $\tilde{u}^0=u^0+u^{\mathrm{r}}$, where $u^{\mathrm{r}}$ is constructed above.
Then $(\tilde{\eta}^0,\tilde{u}^0)$ satisfies the jump condition
\begin{align}\label{20240825n1mn}
\llbracket {P}'(\bar{\rho})\bar{\rho}\mm{div}\tilde{\eta}^0 \mathbb{I}+\mu\mathbb{D}\tilde{u}^0+\lambda\mm{div}\tilde{u}^0\mathbb{I}+\kappa\nabla\tilde{\eta}^0\rrbracket \mathbf{e}^3=0
 \mbox{ on } \Sigma.
\end{align}
Following the argument of the well-posedness result in \eqref{pro:210903}, we easily prove that there exists a unique global solution $(\eta^{\mm{l}},u^{\mm{l}})$ to
the linearized problem \eqref{20210824nnmmn} with initial data $(\eta^{\mm{l}},u^{\mm{l}})|_{t=0}=(\tilde{\eta}^0,\tilde{u}^0)$.
{Moreover, the solution enjoys the same regularity as well as the solution $(\eta,u)$ in Theorem \ref{thm:2109n}.
In particular, following the derivations of \eqref{12281503} by using \eqref{202408201622}, we can deduce from \eqref{20210824nnmmn} and \eqref{202408191959} that
\begin{align}
&\label{202409141532nm}
(\kappa^{-2}\|u_t^{\mm{l}}\|_{\underline{1},0}^2+\kappa^{-1}\|u^{\mm{l}}\|_{\underline{1},1}^2)
+\kappa^{-2}\int_0^t\|u_{t}^{\mm{l}}(\tau)\|_{\underline{1},1}^2\mm{d}\tau\lesssim
(1+\kappa^{-2}){E}^0.
\end{align}
}

Subtracting \eqref{20210824nnmmn} from \eqref{20210824nnmm}, then we obtain that
\begin{equation}\label{20210824nnmmk}
\begin{cases}
\eta^{\mm{d}}_t=u^{\mm{d}}&\mbox{in } \Omega ,\\[0.5mm]
\bar{\rho} u^{\mm{d}}_t
-\mm{div}\big({P}'(\bar{\rho})\bar{\rho}\mm{div}\eta^{\mm{d}} \mathbb{I}+\mu\mathbb{D}u^{\mm{d}}+\lambda\mm{div}u^{\mm{d}}\mathbb{I}
+\kappa(\nabla \eta^{\mm{d}})^\top\big)&\\[0.5mm]
=\mm{div}\big({\mathcal{N}}^{\mu}_{u}+{\mathcal{N}}^{\lambda}_{u}-{\mathcal{N}}^{P}_{\eta}\big) &\mbox{in } \Omega ,\\[0.5mm]
\llbracket  \eta^{\mm{d}}\rrbracket=\llbracket  u^{\mm{d}}\rrbracket=0,
& \mbox{on } \Sigma,\\[0.5mm]
\llbracket
{P}'(\bar{\rho})\bar{\rho}\mm{div}\eta^{\mm{d}}+\mu\mathbb{D}u^{\mm{d}}+
\lambda\mm{div}u^{\mm{d}}
   \mathbb{I}+\kappa\nabla\eta^{\mm{d}} \rrbracket \mathbf{e}^3
=\llbracket{\mathcal{N}}^{P}_{\eta}-{\mathcal{N}}^{\mu}_{u}-{\mathcal{N}}^{\lambda}_{u}\rrbracket
& \mbox{on } \Sigma,\\[0.5mm]
(\eta^{\mm{d}}, u^{\mm{d}})=(0,0) & \mbox{on }\Sigma_{-}^{+},\\[0.5mm]
(\eta^{\mm{d}}, u^{\mm{d}})|_{t=0}=(0,u^{\mm{r}})& \mbox{in }\Omega.
\end{cases}
\end{equation}
Following the arguments  of Lemmas \ref{lem:22083001}--\ref{lem:24083001} with slight modifications, we can deduce from \eqref{20210824nnmmk} that
\begin{align}\label{202109021542n}
\frac{1}{2}\frac{\mm{d}}{\mm{d}t}\bar{\mathcal{E}}^{\mm{d}}
+c\bar{\mathcal{D}}^{\mm{d}}
&\lesssim
\bigg(
\|\nabla\eta\|_{\underline{2},0}^{1/2}\|\nabla \eta\|_{\underline{1},1}^{1/4}\|\nabla \eta\|_{\underline{2},1}^{1/4}\|\nabla u\|_{\underline{3},0}+\|\nabla\eta\|_{\underline{3},0}\|\nabla u\|_{\underline{2},0}^{1/2}\|\nabla u\|_{\underline{1},1}^{1/4}\|\nabla u\|_{\underline{2},1}^{1/4} \nonumber\\
&\quad\quad+\|\nabla\eta\|_{\underline{2},0}^{1/2}\|\nabla \eta\|_{\underline{1},1}^{1/4}\|\nabla \eta\|_{\underline{2},1}^{1/4}\|\nabla\eta\|_{\underline{3},0 }\bigg)
\|(\nabla \eta^{\mm{d}},\nabla u^{\mm{d}})\|_{\underline{3},0}\nonumber\\
&\lesssim(1+\kappa^{-1/2})(\kappa^{-1/8}+\kappa^{-5/8})\sqrt{\mathcal{E}\mathcal{D}}
\|(\sqrt{\kappa}\nabla \eta^{\mm{d}},\nabla u^{\mm{d}})\|_{\underline{3},0},
\end{align}
where we have defined that
\begin{align}
\bar{\mathcal{E}}^{\mm{d}}:=& C_3 \|(\sqrt{\bar{\rho}}u^{\mm{d}},\sqrt{P'(\bar{\rho})\bar{\rho}}\mm{div}\eta^{\mm{d}},
\sqrt{\kappa}\nabla\eta^{\mm{d}})\|_{\underline{3},0}^2\nonumber\\
& +2\sum_{|\alpha|=0}^3\int \bar{\rho}\partial_{\mm{h}}^{\alpha}u^{\mm{d}}\cdot\partial_{\mm{h}}^{\alpha}\eta^{\mm{d}}\mm{d}y+
\|(\sqrt{\mu/2}\mathbb{D}\eta^{\mm{d}},\sqrt{\lambda}\mm{div}\eta^{\mm{d}})\|_{\underline{3},0}^2,\nonumber\\[1mm]
\bar{\mathcal{D}}^{\mm{d}}:=& C_3\|   u^{\mm{d}} \|_{\underline{3},1}^2
+\| \mm{div}\eta^{\mm{d}}\|_{\underline{3},0}^2+{\kappa}\|\eta^{\mm{d}}\|_{\underline{3},1}^2\nonumber
\end{align}
for some suitably large $C_3$,
and that $\bar{\mathcal{E}}^{\mm{d}}$ satisfies
\begin{align*}
\|(\mm{div}\eta^{\mm{d}},u^{\mm{d}})\|_{\underline{3},0}^2+(1+\kappa) \|\eta^{\mm{d}}\|_{\underline{3},1}^2\lesssim\bar{\mathcal{E}}^{\mm{d}}
\lesssim \|(\mm{div}\eta^{\mm{d}},u^{\mm{d}})\|_{\underline{3},0}^2+(1+\kappa) \|\eta^{\mm{d}}\|_{\underline{3},1}^2,
\end{align*}
which implies that
\begin{align}
\bar{\mathcal{E}}^{\mm{d}}\lesssim(1+\kappa^{-1})\bar{\mathcal{D}}^{\mm{d}}.
\label{202021101101141}
\end{align}
Furthermore, using Young's inequality, we can deduce from \eqref{202109021542n} that
\begin{align}\label{202409021542n}
\frac{\mm{d}}{\mm{d}t}\bar{\mathcal{E}}^{\mm{d}}
+c\bar{\mathcal{D}}^{\mm{d}}
&\lesssim(1+\kappa^{-1})(\kappa^{-1/4}+\kappa^{-5/4}){\mathcal{E}\mathcal{D}}.
\end{align}

Thanks to \eqref{202408191959}, we can see that
\begin{align}\label{202409051542n}
\bar{\mathcal{E}}^{\mm{d}}(0)\lesssim\|u^{\mathrm{r}}\|_{\underline{1},2}^2
\lesssim\kappa^{-1/2}{E}^0\|(\eta^0,u^0)\|_{\underline{1},2}^2.
\end{align}
Consequently, in the light of \eqref{202021101101141}--\eqref{202409051542n}, along with \eqref{202109090810}, we then obtain the desired error estimate \eqref{202109061353}.

Additionally, taking the inner product of \eqref{20210824nnmmk} with $\partial_{\mm{h}}^{\alpha}\partial_{3}^2\eta^{\mm{d}}$ for $|\alpha|\leqslant1$ in $L^2$,
and following the derivation of \eqref{20210824nnmmn} with slight modifications, we obtain that
\begin{align}
&\label{202409141542}
\frac{\mm{d}}{\mm{d}t}\overline{\|\partial_{3}^2\eta^{\mm{d}}\|_{\underline{1},0}^2}
+ \|\sqrt{\kappa}\partial_{3}^2\eta^{\mm{d}}\|_{\underline{1},0}^2 \nonumber\\[0.5mm]
&\lesssim
(1+\kappa^{-2})\|(\sqrt{\kappa}\eta^{\mm{d}},u^{\mm{d}})\|_{\underline{2},1}^2 +\kappa^{-1}\|u_t^{\mm{d}}\|_{\underline{1},0}^2
+\kappa^{-1}\|({\mathcal{N}}^{\mu}_{u},{\mathcal{N}}^{\lambda}_{u}, {\mathcal{N}}^{P}_{\eta})\|_{\underline{1},1}^2\nonumber\\[1mm]
&\lesssim
(1+\kappa^{-2})\|(\sqrt{\kappa}\eta^{\mm{d}},u^{\mm{d}})\|_{\underline{2},1}^2 +\kappa^{-1}\|u_t^{\mm{d}}\|_{\underline{1},0}^2
+(\kappa^{-1/4}+\kappa^{-5/4})\mathcal{E}\mathcal{D},
\end{align}
where $\overline{\|\partial_{3}^2\eta^{\mm{d}}\|_{\underline{1},0}^2}$ and $\|\partial_{3}^2\eta^{\mm{d}}\|_{\underline{1},0}^2$ are equivalent to each other.

Taking the inner product of \eqref{20210824nnmmk}$_2$ with ${\kappa}^{-1}\partial_{\mm{h}}^{\alpha}u^{\mm{d}}_t$ in $L^2$, and then using
the integration by parts, and \eqref{20210824nnmmk}$_3$--\eqref{20210824nnmmk}$_5$, we infer that
\begin{align}
&\label{202409141545}
{\kappa}^{-1}\frac{\mm{d}}{\mm{d}t}\left(\int\partial_{\mm{h}}^{\alpha}\left( {P}'(\bar{\rho})\bar{\rho}\mm{div}\eta^{\mm{d}} \mathbb{I}+\kappa\nabla \eta^{\mm{d}}\right):\nabla \partial_{\mm{h}}^{\alpha}u^{\mm{d}}\mm{d}y
+\|\partial_{\mm{h}}^{\alpha}(\sqrt{\mu/2} \mathbb{D}u^{\mm{d}},\sqrt{\lambda}\mm{div}u^{\mm{d}})\|_0^2/2\right)\nonumber\\[1mm]
&
+{\kappa}^{-1}\|\sqrt{\bar{\rho}}\partial_{\mm{h}}^{\alpha}u^{\mm{d}}_t\|_{0}^2\nonumber\\[0.5mm]
&={\kappa}^{-1}\left(\|\partial_{\mm{h}}^{\alpha}(\sqrt{P'(\bar{\rho})\bar{\rho}}\mm{div}u^{\mm{d}},\sqrt{\kappa}\nabla u^{\mm{d}})\|_{0}^2\right)\nonumber\\
&\quad-{\kappa}^{-1}\left(\int\partial_{\mm{h}}^{\alpha}\mathcal{N}^{\mu}_{u}:\nabla(\partial_{\mm{h}}^{\alpha}u^{\mm{d}}_t)^{\top}\mm{d}y
+\int\partial_{\mm{h}}^{\alpha}\mathcal{N}^{\lambda}_{u} \mm{div}\partial_{\mm{h}}^{\alpha}u^{\mm{d}}_t\mm{d}y
-\int\partial_{\mm{h}}^{\alpha}\mathcal{N}^{P}_{\eta}\mm{div}\partial_{\mm{h}}^{\alpha}u^{\mm{d}}_t\mm{d}y\right).
\end{align}
Moreover, by virtue of \eqref{202108270922n12} and \eqref{202108270921} with $i=1$, we can estimate that
\begin{align}\label{202409141745}
&-{\kappa}^{-1}\left(\int\partial_{\mm{h}}^{\alpha}\mathcal{N}^{\mu}_{u}:\nabla(\partial_{\mm{h}}^{\alpha}u^{\mm{d}}_t)^{\top}\mm{d}y
+\int\partial_{\mm{h}}^{\alpha}\mathcal{N}^{\lambda}_{u} \mm{div}\partial_{\mm{h}}^{\alpha}u^{\mm{d}}_t\mm{d}y
-\int\partial_{\mm{h}}^{\alpha}\mathcal{N}^{P}_{\eta}\mm{div}\partial_{\mm{h}}^{\alpha}u^{\mm{d}}_t\mm{d}y\right)\nonumber\\[1mm]
&\lesssim\|\nabla\eta\|_{\underline{2},0}^{1/2}\|\nabla\eta\|_{\underline{1},1}^{1/4}\|\nabla\eta\|_{\underline{2},1}^{1/4}
(\|\nabla \eta\|_{\underline{1},0}+\|\nabla u\|_{\underline{1},0})
\|{\kappa}^{-1}u_t^{\mm{d}}\|_{\underline{1},1}\nonumber\\[1mm]
&
\lesssim(1+\kappa^{-1/2})\kappa^{-1/8}\sqrt{\mathcal{E}\mathcal{D}}\|{\kappa}^{-1}(u_t,u_t^{\mm{l}})\|_{\underline{1},1}.
\end{align}
Finally, making use of \eqref{202109090810}, \eqref{202109061353}, \eqref{202408191959}, \eqref{202409141532nm} and Young's inequality, we immediately obtain \eqref{202409061553} from \eqref{202409141542}--\eqref{202409141745}.
This completes the proof of Theorem \ref{thm:210902}.

\appendix
\section{Analysis tools}\label{sec:09}
\renewcommand\thesection{A}
This appendix is devoted to listing some mathematical analysis tools that have been used in the previous sections.
It should be remarked that in this appendix we still adopt the simplified mathematical notations introduced in Section \ref{202411061605}. In particular, we still use the notation $a\lesssim b$ to mean that $a\leqslant c b$ for some constant $c $, which  may  depend on the domain and the other given parameters in the lemmas below.

\begin{lem}\label{lem:sadf0826}
 Interpolation inequalities in 1D case (see \cite[Theorem]{NLOE}): For any given interval $(a,b)$,
\begin{align}
&\label{202108sdfa261404}\|\varphi\|_{L^{\infty}(\mathbb{R})}\lesssim
\|\varphi\|_{L^2(\mathbb{R})}^{1/2}\|\varphi'\|_{L^2(\mathbb{R})}^{1/2},\\[1mm]
&\label{2021082sadf61417}
\|\varphi\|_{L^{\infty}(a,b)}\lesssim
\|\varphi\|_{L^2(a,b)}^{1/2}\|\varphi\|_{H^1(a,b)}^{1/2}.
\end{align}
\end{lem}
\begin{lem}\label{lem:0826}
Embedding inequalities   in  $H^{s}(\mathbb{R}^2)$ (see \cite[Theorem 1.66]{BJAFP315200} and \cite[Lemma A.1]{JFJHJS2023}):
For any given constant $s>1$,
\begin{align}
&\label{202108261404}|\varphi|_{L^{\infty}}\lesssim|\varphi|_{s},\\[1mm]
&\label{embed3}
|\phi|_{L^{\infty}}\lesssim\|\phi\|_{W^{1,4}(\mathbb{R}^2)}\lesssim |\phi|_{3/2}.
\end{align}
\end{lem}

 \begin{lem} Trace estimate (see \cite[Lemma A.6]{JFJHJS2023}):
For any given  non-negative integer $i$,
\begin{align}
&\label{202108261406}
|\varphi|_{i+1/2}\lesssim\|\varphi\|_{\underline{i+1},0}^{1/2}\|\partial_3\varphi\|_{\underline{i},0}^{1/2}
+\|\varphi\|_{\underline{i+1},0}.
\end{align}
\end{lem}
\begin{lem} Product estimates in $H^{s}(\mathbb{R}^2)$ (see \cite[Corollary 2.86]{BJAFP315200}): For any given $s>0$,
\begin{align}
&\label{202108261400}
|\varphi\phi|_{s}\lesssim |
\varphi|_{L^{\infty}}|\phi|_{s}+|\varphi|_{s}|\phi|_{L^{\infty}},
\end{align}
in particular, by \eqref{202108261404}, we further have, for any given $s>1$,
\begin{align}\label{202108261402}|\varphi\phi|_{s}\lesssim |\varphi|_{s}|\phi|_{s}.  \end{align}
\end{lem}
\begin{lem}\label{lesdfafm:0826}
\begin{enumerate}
\item[(1)]  Interpolation inequalities in Sobolev spaces defined in $\Omega$:
Let $i$ and $j$ be given nonnegative integers and satisfy $0\leqslant j \leqslant i$, then
\begin{align}
&\label{202108261405}
\|f\|_{j}\lesssim \|f\|_{0}^{1-j/i}\|f\|_{i}^{j/i}, \\[1mm]
&\label{202108261424}
\|f\|_{L^{\infty}}
 \lesssim \|f\|_{\underline{2},0}^{3/8}\|f\|_{\underline{2},0}^{1/8} \|f\|_{\underline{1},1}^{3/8}\|f\|_{\underline{2},1}^{1/8}. \end{align}
\item[(2)] Product estimates in  Sobolev spaces defined in $\Omega$:
\begin{align}
&\label{202408241430}
\|\varphi\phi\|_{0}\lesssim\begin{cases}
\|\varphi\|_{\underline{1},0}^{3/8}\|\varphi\|_{\underline{2},0}^{1/8}\|\varphi\|_{\underline{1},1}^{3/8}
\|\varphi\|_{\underline{2},1}^{1/8}\|\phi\|_{0};\\[0.5mm]
 {\|\varphi\|_{\underline{1},0}^{3/4}\|\varphi\|_{\underline{2},0}^{1/4}}\|\phi\|_{0}^{1/2}\|\phi\|_{1}^{1/2}
;\\[0.5mm]
 {\|\varphi\|_{\underline{1},0}^{1/2}\|\varphi\|_{1}^{1/4}\|\varphi\|_{\underline{1},1}^{1/4}}
\|\phi\|_{0}^{1/2}\|\phi\|_{\underline{1},0}^{1/2},
\end{cases}\\[1mm]
&\label{202108261428nm}
\left\|\varphi\phi\right\|_{\underline{i},0}
\lesssim
\begin{cases}
\!\!\!\left.\begin{array}{l}
 {\|\varphi\|_{\underline{2},0}^{1/2}\|\varphi\|_{\underline{1},1}^{1/4}\|\varphi\|_{\underline{2},1}^{1/4}}
\|\phi\|_{\underline{i},0} \\
{\|\varphi\|_{\underline{1},0}^{1/2}\|\varphi\|_{\underline{2},0}^{1/2}}
\|\phi\|_{\underline{i},0}^{1/2}\|\phi\|_{\underline{i},1}^{1/2}\\
\|\varphi\|_{L^{\infty}}
\|\phi\|_{\underline{i},0}+\|\varphi\|_{\underline{i},0}\|\phi\|_{L^{\infty}}\\
 \|\varphi\|_{\underline{i},0}^{1/2}\|\varphi\|_{\underline{i},1}^{1/2}{\|\phi\|_{\underline{1},0}^{3/4}
\|\phi\|_{\underline{2},0}^{1/4}}
+\|\varphi\|_{\underline{1},0}^{3/4}\|\varphi\|_{\underline{2},0}^{1/4}{\|\phi\|_{\underline{i},0}^{1/2}
\|\phi\|_{\underline{i},1}^{1/2}}
\end{array}\right\}  \mbox{for }i=0,\ 1;\\[2mm]
\!\!\!\left.\begin{array}{l}
{\|\varphi\|_{\underline{2},0}^{1/2}\|\varphi\|_{\underline{1},1}^{1/4}\|\varphi\|_{\underline{2},1}^{1/4}}
\|\phi\|_{\underline{i},0}
+\|\varphi\|_{\underline{i},0}^{1/2}\|\varphi\|_{\underline{i},1}^{1/2}{\|\phi\|_{\underline{1},0}^{3/4}
\|\phi\|_{\underline{2},0}^{1/4}}
\\
{\|\varphi\|_{\underline{2},0}^{1/2}\|\varphi\|_{\underline{1},1}^{1/4}\|\varphi\|_{\underline{2},1}^{1/4}}
\|\phi\|_{\underline{i},0}
+\|\varphi\|_{\underline{i},0}^{1/2}\|\varphi\|_{\underline{i+1},0}^{1/2}
{\|\phi\|_{\underline{1},0}^{1/2}\|\phi\|_{1}^{1/4}\|\phi\|_{\underline{1},1}^{1/4}}
 \end{array}\right\}  \mbox{for }1\leqslant i\leqslant2;\\[1mm]
{\|\varphi\|_{\underline{2},0}^{1/2}\|\varphi\|_{\underline{1},1}^{1/4}\|\varphi\|_{\underline{2},1}^{1/4}}
\|\phi\|_{\underline{i},0}
+\|\varphi\|_{\underline{i},0}{\|\phi\|_{\underline{2},0}^{1/2}
\|\phi\|_{\underline{1},1}^{1/4}\|\phi\|_{\underline{2},1}^{1/4}}\; \mbox{ for }1\leqslant i\leqslant3.
\end{cases}\quad
\end{align}
\end{enumerate}
\end{lem}
\begin{pf}
The estimate \eqref{202108261405} can be found in \cite[Theroem 5.2]{ARAJJFF}.
By \eqref{202108sdfa261404}--\eqref{2021082sadf61417} and the density theorem of Sobolev spaces (see \cite[Theorem 3.22]{ARAJJFF}),
we easily obtain the estimates  \eqref{202108261424}--\eqref{202408241430}. Similarly, using
\eqref{202408241430} and the interpolation inequality, we can get the estimate \eqref{202108261428nm},
please refer to \cite[Lemma 3.1]{JFJSGSfdsa} for details.
\hfill$\Box$
\end{pf}

\begin{lem} Korn's inequality (see \cite[Proposition A.8]{JJTIWYJTC}):
\begin{align}\label{psfxdxx}
\|w\|_1 \lesssim  \|\mathbb{D}w\|_0 \quad\mbox{ for any }w\in H_0^1.
\end{align}
\end{lem}

\begin{lem}\label{pro:040s1nasfxdxx}
Homeomorphism mapping theorem (see \cite[Lemma 4.2]{JFJSOMITIN}): There is a positive constant $\delta_3$, such that for any $\varphi\in H^{1,3}$ satisfying
$\|\nabla \varphi\|_{L^{\infty}}\leqslant \delta_3$, we have (after possibly being redefined on a set of measure zero) $1/2\leqslant \det(\nabla \varphi+I)\leqslant 3/2$ and
\begin{align}
&   \psi : \Omega\to \Omega \mbox{ is a }C^1\mbox{-homeomorphism mapping},
\end{align}
where $\psi:=\varphi+y$.
\end{lem}

\begin{lem}\label{prosfxdxx}
 Stratified elliptic theory (see \cite[Lemma A.10]{JJTIWYJTC}): Let $i\geqslant 0$, $\mathcal{F}\in H^{i }$ and $\mathcal{G}\in H^{i+1/2}$, then there exists a unique solution $u\in H^{i+2}$
to the following  stratified Lam\'e problem:
\begin{equation*} \begin{cases}
\mu\Delta w+(\mu+\lambda)\nabla \mm{div}w= \mathcal{F}&\mbox{in }\Omega, \\[1mm]
\llbracket w  \rrbracket=0&\mbox{on }\Sigma, \\[1mm]
 \llbracket\mu\mathbb{D}w+ \lambda\mm{div}w\mathbb{I}\rrbracket \mathbf{e}^3= \mathcal{G}&\mbox{on }\Sigma, \\[1mm]
w=0 &\mbox{on }\Sigma_{-}^{+};
\end{cases} \end{equation*}
moreover, the solution enjoys
\begin{equation}
\label{Ellipticestimate}
\|u\|_{i+2}\lesssim
\|\mathcal{F}\|_{i}+|\mathcal{G}|_{i+1/2}.
\end{equation}
\end{lem}
\vspace{5mm} \noindent\textbf{Acknowledgements.}
The research of Fei Jiang was supported by  NSFC (Nos. 12022102 and 12231016), the Natural Science Foundation of Fujian Province of China (Nos. 2022J01105 and 2024J011011), the Central Guidance on Local Science and Technology Development Fund of Fujian Province (No. 2023L3003) and Fujian Alliance Of Mathematics (No. 2025SXLMMS01), and the research of Youyi Zhao was supported
by NSFC (No. 12401289), the Natural Science Foundation of Fujian Province of China (No. 2024J08029)
and the Research Foundation of Fuzhou University (No. XRC-24050).


\vspace{5mm}
\noindent\textbf{Conflict of Interest.}
On behalf of all authors, the corresponding author states that there is no conflict of interest.

\renewcommand\refname{References}
\renewenvironment{thebibliography}[1]{%
\section*{\refname}
\list{{\arabic{enumi}}}{\def\makelabel##1{\hss{##1}}\topsep=0mm
\parsep=0mm
\partopsep=0mm\itemsep=0mm
\labelsep=1ex\itemindent=0mm
\settowidth\labelwidth{\small[#1]}%
\leftmargin\labelwidth \advance\leftmargin\labelsep
\advance\leftmargin -\itemindent
\usecounter{enumi}}\small
\def\newblock{\ }
\sloppy\clubpenalty4000\widowpenalty4000
\sfcode`\.=1000\relax}{\endlist}
\bibliographystyle{model1b-num-names}

\begin{thebibliography}{53}
\expandafter\ifx\csname natexlab\endcsname\relax\def\natexlab#1{#1}\fi
\providecommand{\bibinfo}[2]{#2}
\ifx\xfnm\relax \def\xfnm[#1]{\unskip,\space#1}\fi
\bibitem[{Adams and Fournier(2005)}]{ARAJJFF}
\bibinfo{author}{R.A. Adams}, \bibinfo{author}{J.J.F. Fournier},
  \bibinfo{title}{{Sobolev Space}}, \bibinfo{publisher}{Academic Press: New
  York}, \bibinfo{year}{2005}.
\bibitem[{Bahouri et~al.(2011)Bahouri, Chemin and Danchin}]{BJAFP315200}
\bibinfo{author}{H.~Bahouri}, \bibinfo{author}{J.Y. Chemin},
  \bibinfo{author}{R.~Danchin}, \bibinfo{title}{Fourier Analysis and Nonlinear
  Partial Differential Equations}, \bibinfo{publisher}{Springer},
  \bibinfo{year}{2011}.
\bibitem[{Bejaoui and Majdoub(2013)}]{MR2990047}
\bibinfo{author}{O.~Bejaoui}, \bibinfo{author}{M.~Majdoub},
  \bibinfo{title}{Global weak solutions for some {O}ldroyd models},
  \bibinfo{journal}{J. Differential Equations} \bibinfo{volume}{254}
  (\bibinfo{year}{2013}) \bibinfo{pages}{660--685}.
\bibitem[{Castro et~al.(2019)Castro, C{\'o}rdoba, Fefferman, Gancedo and
  G{\'o}mez-Serrano}]{castro2019splash}
\bibinfo{author}{A.~Castro}, \bibinfo{author}{D.~C{\'o}rdoba},
  \bibinfo{author}{C.L. Fefferman}, \bibinfo{author}{F.~Gancedo},
  \bibinfo{author}{J.~G{\'o}mez-Serrano}, \bibinfo{title}{{Splash singularities
  for the free boundary Navier--Stokes equations}}, \bibinfo{journal}{Ann. PDE}
  \bibinfo{volume}{5} (\bibinfo{year}{2019}) \bibinfo{pages}{1--117}.
\bibitem[{Chen et~al.(2017)Chen, Hu and Wang}]{CCRMHJLWDH}
\bibinfo{author}{R.M. Chen}, \bibinfo{author}{J.L. Hu}, \bibinfo{author}{D.H.
  Wang}, \bibinfo{title}{Linear stability of compressible vortex sheets in
  two-dimensional elastodynamics}, \bibinfo{journal}{Adv. Math.}
  \bibinfo{volume}{311} (\bibinfo{year}{2017}) \bibinfo{pages}{18--60}.
\bibitem[{Chen et~al.(2020)Chen, Hu and Wang}]{chen2020linear}
\bibinfo{author}{R.M. Chen}, \bibinfo{author}{J.L. Hu}, \bibinfo{author}{D.H.
  Wang}, \bibinfo{title}{Linear stability of compressible vortex sheets in 2{D}
  elastodynamics: variable coefficients}, \bibinfo{journal}{Math. Ann.}
  \bibinfo{volume}{376} (\bibinfo{year}{2020}) \bibinfo{pages}{863--912}.
\bibitem[{Chen and Zhang(2006)}]{chen2006global}
\bibinfo{author}{Y.~Chen}, \bibinfo{author}{P.~Zhang}, \bibinfo{title}{The
  global existence of small solutions to the incompressible viscoelastic fluid
  system in 2 and 3 space dimensions}, \bibinfo{journal}{Comm. Partial
  Differential Equations} \bibinfo{volume}{31} (\bibinfo{year}{2006})
  \bibinfo{pages}{1793--1810}.
\bibitem[{Di~Iorio et~al.(2020)Di~Iorio, Marcati and Spirito}]{di2020splash}
\bibinfo{author}{E.~Di~Iorio}, \bibinfo{author}{P.~Marcati},
  \bibinfo{author}{S.~Spirito}, \bibinfo{title}{Splash singularity for a
  free-boundary incompressible viscoelastic fluid model},
  \bibinfo{journal}{Adv. Math.} \bibinfo{volume}{368} (\bibinfo{year}{2020})
  \bibinfo{pages}{107124}.
\bibitem[{Gu and Lei(2022)}]{GL2022}
\bibinfo{author}{X.M. Gu}, \bibinfo{author}{Z.~Lei}, \bibinfo{title}{{Local
  well-posedness of free-boundary incompressible elastodynamics with surface
  tension via vanishing viscosity limit}}, \bibinfo{journal}{Arch. Ration.
  Mech. Anal.} \bibinfo{volume}{245} (\bibinfo{year}{2022})
  \bibinfo{pages}{1285--1338}.
\bibitem[{Gu and Mei(2023)}]{GuMei2023}
\bibinfo{author}{X.M. Gu}, \bibinfo{author}{Y.~Mei}, \bibinfo{title}{{Vanishing
  viscosity limits for the free boundary problem of compressible viscoelastic
  fluids with surface tension}}, \bibinfo{journal}{Sci. China Math.}
  \bibinfo{volume}{66} (\bibinfo{year}{2023}) \bibinfo{pages}{1263--1300}.
\bibitem[{Hu(2018)}]{HXPGETDCMFJDE}
\bibinfo{author}{X.P. Hu}, \bibinfo{title}{{Global existence of weak solutions
  to two dimensional compressible viscoelastic flows}}, \bibinfo{journal}{J.
  Differential Equations} \bibinfo{volume}{265} (\bibinfo{year}{2018})
  \bibinfo{pages}{3130--3167}.
\bibitem[{Hu and Lin(2016)}]{HXPLFHG}
\bibinfo{author}{X.P. Hu}, \bibinfo{author}{F.H. Lin}, \bibinfo{title}{Global
  solution to two dimensional incompressible viscoelastic fluid with
  discontinuous data}, \bibinfo{journal}{Comm. Pure Appl. Math.}
  \bibinfo{volume}{69} (\bibinfo{year}{2016}) \bibinfo{pages}{372--404}.
\bibitem[{Hu and Wang(2010)}]{hu2010local}
\bibinfo{author}{X.P. Hu}, \bibinfo{author}{D.H. Wang}, \bibinfo{title}{Local
  strong solution to the compressible viscoelastic flow with large data},
  \bibinfo{journal}{J. Differential Equations} \bibinfo{volume}{249}
  (\bibinfo{year}{2010}) \bibinfo{pages}{1179--1198}.
\bibitem[{Hu and Wang(2011)}]{HXWDHTTIVPTCVF2}
\bibinfo{author}{X.P. Hu}, \bibinfo{author}{D.H. Wang}, \bibinfo{title}{{Global
  existence for the multi-dimensional compressible viscoelastic flows}},
  \bibinfo{journal}{J. Differential Equations} \bibinfo{volume}{250}
  (\bibinfo{year}{2011}) \bibinfo{pages}{1200--1231}.
\bibitem[{Hu and Wang(2012)}]{hu2012strong}
\bibinfo{author}{X.P. Hu}, \bibinfo{author}{D.H. Wang}, \bibinfo{title}{Strong
  solutions to the three-dimensional compressible viscoelastic fluids},
  \bibinfo{journal}{J. Differential Equations} \bibinfo{volume}{252}
  (\bibinfo{year}{2012}) \bibinfo{pages}{4027--4067}.
\bibitem[{Hu and Wang(2015)}]{HXWDHTTIVPTCVF1}
\bibinfo{author}{X.P. Hu}, \bibinfo{author}{D.H. Wang}, \bibinfo{title}{{The
  initial-boundary value problem for the compressible viscoelastic flows}},
  \bibinfo{journal}{Discrete Contin. Dyn. Syst.} \bibinfo{volume}{35}
  (\bibinfo{year}{2015}) \bibinfo{pages}{917--934}.
\bibitem[{Hu and Wu(2013)}]{HXPWGC}
\bibinfo{author}{X.P. Hu}, \bibinfo{author}{G.C. Wu}, \bibinfo{title}{{Global
  existence and optimal decay rates for three-dimensional compressible
  viscoelastic flows}}, \bibinfo{journal}{SIAM J. Math. Anal.}
  \bibinfo{volume}{45} (\bibinfo{year}{2013}) \bibinfo{pages}{2815--2833}.
\bibitem[{Hu and Zhao(2020)}]{MR4064197}
\bibinfo{author}{X.P. Hu}, \bibinfo{author}{W.B. Zhao}, \bibinfo{title}{Global
  existence of compressible dissipative elastodynamics systems with zero shear
  viscosity in two dimensions}, \bibinfo{journal}{Arch. Ration. Mech. Anal.}
  \bibinfo{volume}{235} (\bibinfo{year}{2020}) \bibinfo{pages}{1177--1243}.
\bibitem[{Huang et~al.(2022)Huang, Wang, Wen and Zi}]{MR4335131}
\bibinfo{author}{J.R. Huang}, \bibinfo{author}{Y.H. Wang},
  \bibinfo{author}{H.Y. Wen}, \bibinfo{author}{R.Z. Zi},
  \bibinfo{title}{Optimal time-decay estimates for an {O}ldroyd-{B} model with
  zero viscosity}, \bibinfo{journal}{J. Differential Equations}
  \bibinfo{volume}{306} (\bibinfo{year}{2022}) \bibinfo{pages}{456--491}.
\bibitem[{Ishigaki(2020)}]{MR4152207}
\bibinfo{author}{Y.~Ishigaki}, \bibinfo{title}{Diffusion wave phenomena and
  {$L^p$} decay estimates of solutions of compressible viscoelastic system},
  \bibinfo{journal}{J. Differential Equations} \bibinfo{volume}{269}
  (\bibinfo{year}{2020}) \bibinfo{pages}{11195--11230}.
\bibitem[{Jang et~al.(2016{\natexlab{a}})Jang, Tice and Wang}]{JJTIIAWangYJC}
\bibinfo{author}{J.~Jang}, \bibinfo{author}{I.~Tice}, \bibinfo{author}{Y.J.
  Wang}, \bibinfo{title}{{The compressible viscous surface-internal wave
  problem: local well-posedness}}, \bibinfo{journal}{SIAM J. Math. Anal.}
  \bibinfo{volume}{48} (\bibinfo{year}{2016}{\natexlab{a}})
  \bibinfo{pages}{2602--2673}.
\bibitem[{Jang et~al.(2016{\natexlab{b}})Jang, Tice and Wang}]{JJTIWYJTC}
\bibinfo{author}{J.~Jang}, \bibinfo{author}{I.~Tice}, \bibinfo{author}{Y.J.
  Wang}, \bibinfo{title}{{The compressible viscous surface-internal wave
  problem: stability and vanishing surface tension limit}},
  \bibinfo{journal}{Comm. Math. Phys.} \bibinfo{volume}{343}
  (\bibinfo{year}{2016}{\natexlab{b}}) \bibinfo{pages}{1039--1113}.
\bibitem[{Jiang et~al.(2023)Jiang, Jiang and Jiang}]{JFJHJS2023}
\bibinfo{author}{F.~Jiang}, \bibinfo{author}{H.~Jiang},
  \bibinfo{author}{S.~Jiang}, \bibinfo{title}{{Rayleigh--Taylor instability in
  stratified compressible fluids with/without the interfacial surface
  tension}}, \bibinfo{journal}{ariXiv:submit/5130422 [math.AP] 23 Sep 2023}
  (\bibinfo{year}{2023}).
\bibitem[{Jiang and Jiang(2019{\natexlab{a}})}]{JFJSCVPDE1}
\bibinfo{author}{F.~Jiang}, \bibinfo{author}{S.~Jiang},
  \bibinfo{title}{{Nonlinear stability and instability in the Rayleigh--Taylor
  problem of stratified compressible MHD fluids}}, \bibinfo{journal}{Calc. Var.
  Partial Differential Equations} \bibinfo{volume}{58}
  (\bibinfo{year}{2019}{\natexlab{a}}) \bibinfo{pages}{Art. 29, 61 pp.}
\bibitem[{Jiang and Jiang(2019{\natexlab{b}})}]{JFJSOMITIN}
\bibinfo{author}{F.~Jiang}, \bibinfo{author}{S.~Jiang}, \bibinfo{title}{{On
  magnetic inhibition theory in non-resistive magnetohydrodynamic fluids}},
  \bibinfo{journal}{Arch. Rational Mech. Anal.} \bibinfo{volume}{233}
  (\bibinfo{year}{2019}{\natexlab{b}}) \bibinfo{pages}{749--798}.
\bibitem[{Jiang and Jiang(2021{\natexlab{a}})}]{JFJSGSfdsa}
\bibinfo{author}{F.~Jiang}, \bibinfo{author}{S.~Jiang}, \bibinfo{title}{{
  Asymptotic behaviors of global solutions to the two-dimensional non-resistive
  MHD equations with large initial perturbations}}, \bibinfo{journal}{Adv.
  Math.} \bibinfo{volume}{393} (\bibinfo{year}{2021}{\natexlab{a}})
  \bibinfo{pages}{108084}.
\bibitem[{Jiang and Jiang(2021{\natexlab{b}})}]{JFJSsdafsafJMFMOSERT}
\bibinfo{author}{F.~Jiang}, \bibinfo{author}{S.~Jiang}, \bibinfo{title}{{
  Strong solutions of the equations for viscoelastic fluids in some classes of
  large data}}, \bibinfo{journal}{J. Differential Equations}
  \bibinfo{volume}{282} (\bibinfo{year}{2021}{\natexlab{b}})
  \bibinfo{pages}{148--183}.
\bibitem[{Jiang and Jiang(2024)}]{JFJSGS}
\bibinfo{author}{F.~Jiang}, \bibinfo{author}{S.~Jiang},
  \bibinfo{title}{{Non-formation of singularities on the free-boundary of an
  incompressible viscoelastic fluid with relatively larger elasticity
  coefficient}},
  \bibinfo{journal}{https://pan.baidu.com/s/1sRQ3zXdrnmM1Ary-1Gdfzw?pwd=1234}
  (\bibinfo{year}{2024}).
\bibitem[{Jiang et~al.(2017)Jiang, Jiang and Wu}]{JFJWGCOSdd}
\bibinfo{author}{F.~Jiang}, \bibinfo{author}{S.~Jiang}, \bibinfo{author}{G.C.
  Wu}, \bibinfo{title}{{ On stabilizing effect of elasticity in the
  Rayleigh--Taylor problem of stratified viscoelastic fluids}},
  \bibinfo{journal}{J. Funct. Anal.} \bibinfo{volume}{272}
  (\bibinfo{year}{2017}) \bibinfo{pages}{3763--3824}.
\bibitem[{Jiang and Liu(2020)}]{jiang2020}
\bibinfo{author}{F.~Jiang}, \bibinfo{author}{M.M. Liu},
  \bibinfo{title}{Nonlinear stability of the viscoelastic {B}\'enard problem},
  \bibinfo{journal}{Nonlinearity} \bibinfo{volume}{33} (\bibinfo{year}{2020})
  \bibinfo{pages}{1677--1704}.
\bibitem[{Jiang et~al.(2016)Jiang, Wu and Zhong}]{FJWGCZXOE}
\bibinfo{author}{F.~Jiang}, \bibinfo{author}{G.C. Wu},
  \bibinfo{author}{X.~Zhong}, \bibinfo{title}{{ On exponential stability of
  gravity driven viscoelastic flows}}, \bibinfo{journal}{J. Differential
  Equations} \bibinfo{volume}{260} (\bibinfo{year}{2016})
  \bibinfo{pages}{7498--7534}.
\bibitem[{Lei(2016)}]{LZSTCZY2}
\bibinfo{author}{Z.~Lei}, \bibinfo{title}{{Global well-posedness of
  incompressible elastodynamics in two dimensions}}, \bibinfo{journal}{Comm.
  Pure Appl. Math.} \bibinfo{volume}{69} (\bibinfo{year}{2016})
  \bibinfo{pages}{2072--2106}.
\bibitem[{Lei et~al.(2018)Lei, Liu and Zhou}]{LZLCZYGARMA}
\bibinfo{author}{Z.~Lei}, \bibinfo{author}{C.~Liu}, \bibinfo{author}{Y.~Zhou},
  \bibinfo{title}{Global solutions for incompressible viscoelastic fluids},
  \bibinfo{journal}{Arch. Ration. Mech. Anal.} \bibinfo{volume}{188}
  (\bibinfo{year}{2018}) \bibinfo{pages}{371--398}.
\bibitem[{Lei and Zhou(2005)}]{LZZYGE}
\bibinfo{author}{Z.~Lei}, \bibinfo{author}{Y.~Zhou}, \bibinfo{title}{Global
  existence of classical solutions for the two-dimensional {O}ldroyd model via
  the incompressible limit}, \bibinfo{journal}{SIAM Journal on Mathematical
  Analysis} \bibinfo{volume}{37} (\bibinfo{year}{2005})
  \bibinfo{pages}{797--814}.
\bibitem[{Li et~al.(2019)Li, Wang and Zhang}]{LHWWZZF}
\bibinfo{author}{H.~Li}, \bibinfo{author}{W.~Wang}, \bibinfo{author}{Z.F.
  Zhang}, \bibinfo{title}{Well-posedness of the free boundary problem in
  incompressible elastodynamics}, \bibinfo{journal}{J. Differential Equations}
  \bibinfo{volume}{267} (\bibinfo{year}{2019}) \bibinfo{pages}{6604--6643}.
\bibitem[{Lin(2012)}]{LFHSM}
\bibinfo{author}{F.H. Lin}, \bibinfo{title}{Some analytical issues for elastic
  complex fluids}, \bibinfo{journal}{Comm. Pure Appl. Math.}
  \bibinfo{volume}{65} (\bibinfo{year}{2012}) \bibinfo{pages}{893--919}.
\bibitem[{Lin et~al.(2005)Lin, Liu and Zhang}]{LFHLCZPO2}
\bibinfo{author}{F.H. Lin}, \bibinfo{author}{C.~Liu},
  \bibinfo{author}{P.~Zhang}, \bibinfo{title}{On hydrodynamics of viscoelastic
  fluids}, \bibinfo{journal}{Comm. Pure Appl. Math.} \bibinfo{volume}{58}
  (\bibinfo{year}{2005}) \bibinfo{pages}{1437--1471}.
\bibitem[{Lin and Zhang(2008)}]{LFZPGCon2}
\bibinfo{author}{F.H. Lin}, \bibinfo{author}{P.~Zhang}, \bibinfo{title}{On the
  initial-boundary value problem of the incompressible viscoelastic fluid
  system}, \bibinfo{journal}{Comm. Pure Appl. Math.} \bibinfo{volume}{61}
  (\bibinfo{year}{2008}) \bibinfo{pages}{539--558}.
\bibitem[{Lions and Masmoudi(2000)}]{LPLMNCAMN}
\bibinfo{author}{P.L. Lions}, \bibinfo{author}{N.~Masmoudi},
  \bibinfo{title}{Global solutions for some {O}ldroyd models of non-{N}ewtonian
  flows}, \bibinfo{journal}{Chin. Ann. Math. Ser. B} \bibinfo{volume}{21}
  (\bibinfo{year}{2000}) \bibinfo{pages}{131--146}.
\bibitem[{Masmoudi(2013)}]{MNIMG}
\bibinfo{author}{N.~Masmoudi}, \bibinfo{title}{Global existence of weak
  solutions to the fene dumbbell model of polymeric flows},
  \bibinfo{journal}{Invent. Math.} \bibinfo{volume}{191} (\bibinfo{year}{2013})
  \bibinfo{pages}{427--500}.
\bibitem[{Nirenberg(1959)}]{NLOE}
\bibinfo{author}{L.~Nirenberg}, \bibinfo{title}{{On elliptic partial
  differential equations}}, \bibinfo{journal}{Estratto dagli Annali della
  Scuola Normale Superiore di Pisa Serie III} \bibinfo{volume}{XIII. Fasc. II}
  (\bibinfo{year}{1959}).
\bibitem[{Pan and Xu(2019)}]{MR3927502}
\bibinfo{author}{X.H. Pan}, \bibinfo{author}{J.~Xu}, \bibinfo{title}{Global
  existence and optimal decay estimates of the compressible viscoelastic flows
  in {$L^p$} critical spaces}, \bibinfo{journal}{Discrete Contin. Dyn. Syst.}
  \bibinfo{volume}{39} (\bibinfo{year}{2019}) \bibinfo{pages}{2021--2057}.
\bibitem[{Pan et~al.(2022)Pan, Xu and Zhu}]{MR4435915}
\bibinfo{author}{X.H. Pan}, \bibinfo{author}{J.~Xu}, \bibinfo{author}{Y.~Zhu},
  \bibinfo{title}{Global existence in critical spaces for non-{N}ewtonian
  compressible viscoelastic flows}, \bibinfo{journal}{J. Differential
  Equations} \bibinfo{volume}{331} (\bibinfo{year}{2022})
  \bibinfo{pages}{162--191}.
\bibitem[{Qian and Zhang(2010)}]{QJZZPWZZF}
\bibinfo{author}{J.Z. Qian}, \bibinfo{author}{Z.F. Zhang},
  \bibinfo{title}{{Global well-posedness for compressible viscoelastic fluids
  near equilibrium}}, \bibinfo{journal}{Arch. Ration. Mech. Anal.}
  \bibinfo{volume}{198} (\bibinfo{year}{2010}) \bibinfo{pages}{835--868}.
\bibitem[{Sideris and Thomases(2005)}]{STCBT1}
\bibinfo{author}{T.~Sideris}, \bibinfo{author}{B.~Thomases},
  \bibinfo{title}{{Global existence for three-dimensional incompressible
  isotropic elastodynamics via the incompressible limit}},
  \bibinfo{journal}{Comm. Pure Appl. Math.} \bibinfo{volume}{58}
  (\bibinfo{year}{2005}) \bibinfo{pages}{750--788}.
\bibitem[{Sideris and Thomases(2007)}]{STCBT2}
\bibinfo{author}{T.~Sideris}, \bibinfo{author}{B.~Thomases},
  \bibinfo{title}{{Global existence for three-dimensional incompressible
  isotropic elastodynamics}}, \bibinfo{journal}{Comm. Pure Appl. Math.}
  \bibinfo{volume}{60} (\bibinfo{year}{2007}) \bibinfo{pages}{1707--1730}.
\bibitem[{Trakhinin(2018)}]{TrakhininJDE}
\bibinfo{author}{Y.~Trakhinin}, \bibinfo{title}{{Well-posedness of the free
  boundary problem in compressible elastodynamics}}, \bibinfo{journal}{J.
  Differential Equations} \bibinfo{volume}{264} (\bibinfo{year}{2018})
  \bibinfo{pages}{1661--1715}.
\bibitem[{Wang and Wen(2020)}]{MR4059368}
\bibinfo{author}{W.J. Wang}, \bibinfo{author}{H.Y. Wen}, \bibinfo{title}{The
  {C}auchy problem for an {O}ldroyd-{B} model in three dimensions},
  \bibinfo{journal}{Math. Models Methods Appl. Sci.} \bibinfo{volume}{30}
  (\bibinfo{year}{2020}) \bibinfo{pages}{139--179}.
\bibitem[{Wang(2017)}]{WXCG}
\bibinfo{author}{X.C. Wang}, \bibinfo{title}{{Global existence for the 2D
  incompressible isotropic elastodynamics for small initial data}},
  \bibinfo{journal}{Ann. Henri Poincar\'e} \bibinfo{volume}{18}
  (\bibinfo{year}{2017}) \bibinfo{pages}{1213--1267}.
\bibitem[{Xu et~al.(2013)Xu, Zhang and Zhang}]{XLZPZZFGAR}
\bibinfo{author}{L.~Xu}, \bibinfo{author}{P.~Zhang}, \bibinfo{author}{Z.F.
  Zhang}, \bibinfo{title}{{Global solvability of a free boundary
  three-dimensional incompressible viscoelastic fluid system with surface
  tension}}, \bibinfo{journal}{Arch. Ration. Mech. Anal.} \bibinfo{volume}{208}
  (\bibinfo{year}{2013}) \bibinfo{pages}{753--803}.
\bibitem[{Zhang(2022)}]{ZJY2022}
\bibinfo{author}{J.Y. Zhang}, \bibinfo{title}{{Local well-posedness and
  incompressible limit of the free-boundary problem in compressible
  elastodynamics}}, \bibinfo{journal}{Arch. Ration. Mech. Anal.}
  \bibinfo{volume}{244} (\bibinfo{year}{2022}) \bibinfo{pages}{599--697}.
\bibitem[{Zhao(2024)}]{MR4762618}
\bibinfo{author}{Y.Y. Zhao}, \bibinfo{title}{On the inhibition of
  {R}ayleigh--{T}aylor instability by surface tension in stratified
  incompressible viscous fluids under {L}agrangian coordinates},
  \bibinfo{journal}{Discrete Contin. Dyn. Syst.} \bibinfo{volume}{44}
  (\bibinfo{year}{2024}) \bibinfo{pages}{2815--2845}.
\bibitem[{Zhu(2022)}]{MR4350195}
\bibinfo{author}{Y.~Zhu}, \bibinfo{title}{Global classical solutions of 3{D}
  compressible viscoelastic system near equilibrium}, \bibinfo{journal}{Calc.
  Var. Partial Differential Equations} \bibinfo{volume}{61}
  (\bibinfo{year}{2022}) \bibinfo{pages}{Paper No. 21, 22}.

\end{thebibliography}

\end{document}